\documentclass[a4paper]{amsart}
\usepackage{color}
\usepackage[latin1]{inputenc}
\usepackage[T1]{fontenc}
\usepackage{amsfonts}
\usepackage{amssymb}
\usepackage{amsmath}
\usepackage{amsthm}
\usepackage{stmaryrd}
\usepackage{enumerate}
\usepackage{multirow}
\usepackage{graphicx}

\usepackage{graphicx,type1cm,eso-pic,color}
\usepackage{pstricks}
\usepackage{float}
\newtheorem{theorem}{Theorem}[section]
\newtheorem{lemma}[theorem]{Lemma}
\newtheorem{proposition}[theorem]{Proposition}

\newtheorem{assumption}[theorem]{Assumption}
\newtheorem{remark}[theorem]{Remark}
\newtheorem{example}[theorem]{Example}
\begin{document}
\setlength\arraycolsep{2pt}
\title{Nonparametric Estimation in SDE Models Involving an Explanatory Process}
\author{Fabienne COMTE$^{\dag}$}
\email{fabienne.comte@u-paris.fr}
\author{Nicolas MARIE$^{\diamond}$}
\email{nmarie@parisnanterre.fr}
\date{}
\maketitle
\noindent
$^{\dag}$Universit\'e Paris Cit\'e, CNRS, MAP5, F-75006 Paris, France.\\
$^{\diamond}$Universit\'e Paris Nanterre, CNRS, Modal'X, 92001 Nanterre, France.
%


%
\begin{abstract}
This paper deals with the process $X = (X_t)_{t\in [0,T]}$ defined by the stochastic differential equation (SDE) $dX_t = (a(X_t) + b(Y_t))dt +\sigma(X_t)dW_1(t)$, where $W_1$ is a Brownian motion and $Y$ is an exogenous process. The first task - of probabilistic nature - is to properly define the model, to prove the existence and uniqueness of the solution of such an equation, and then to establish the existence and a suitable control of a density with respect to the Lebesgue measure of the distribution of $(X_t,Y_t)$ ($t > 0$). In the second part of the paper, a risk bound and a rate of convergence in specific Sobolev spaces are established for a copies-based projection least squares estimator of the $\mathbb R^2$-valued function $(a,b)$. Moreover, a model selection procedure making the adequate bias-variance compromise both in theory and practice is investigated.
\end{abstract}
\noindent
{\bf MSC2020 subject classifications:} 60H10; 62G05.\\
{\bf Keywords:} Diffusion processes; Explanatory process; Nonparametric estimation; Projection least squares estimator.
%


%
\section{Introduction}\label{section_introduction}
Let $(\Omega,\mathcal F,\mathbb F,\mathbb P)$ be the filtered probability space induced by the $2$-dimensional Brownian motion $W$ - with independent components $W_1$ and $W_2$ - defined on $[0,T]$ ($T > 0$). Let $\mathbb F_{\ell}$ be the natural filtration of $W_{\ell}$ for $\ell = 1,2$. In this paper, we consider the stochastic differential equation
\begin{equation}\label{main_equation}
X_t = x_0 +\int_{0}^{t}a(X_s)ds +
\int_{0}^{t}b(Y_s)ds +
\int_{0}^{t}\sigma(X_s)dW_1(s)
\textrm{ $;$ }t\in [0,T],
\end{equation}
under the following assumption:
%


%
\begin{assumption}\label{main_assumption}
$x_0\in\mathbb R$, $a$, $b$ and $\sigma$ are Lipschitz continuous from $\mathbb R$ into itself, and $Y$ is a $\mathbb F_2$-adapted explanatory continuous process.
\end{assumption}
In other words, we intend to add in the dynamics of the process $X$ the influence of an exogenous explanatory process $Y$. It is not difficult to assert that, under Assumption \ref{main_assumption}, Equation (\ref{main_equation}) admits a unique strong solution.

Our aim is to build and study nonparametric estimators of the functions $a(.)$ and $b(.)$ in the copies-based estimation setting: this means that we observe $N$ independent and identically distributed $\mathbb R^2$-valued processes $(X^i,Y^i)$ ($i = 1,\dots,N$), where $X^i :=\mathcal I(Y^i,W_{1}^{i})$, $W^1,\dots,W^N$ are $N$ independent copies of $W$, and $\mathcal I$ is the solution map for Equation (\ref{main_equation}).

Because of our statistical question, probabilistic results are required. Indeed, if $Y$ is a diffusion process, then well-known results on multidimensional SDE ensure the existence and a suitable control of the density of $(X_t,Y_t)$ ($t > 0$) (see e.g. Menozzi et al. \cite{MPZ21}). However, in such models, the explanatory process may not be a diffusion, and have a quite different nature. This is why we first establish conditions ensuring useful results on the distribution of $(X_t,Y_t)$. Thanks to Malliavin calculus-based techniques (see Nualart \cite{NUALART06}, Chapter 2), we can prove that, when $Y$ is only a regular enough $\mathbb F_2$-adapted process, relevant properties are available. Moreover, thanks again to the Malliavin calculus, a preliminary result (see Proposition \ref{technical_result_empirical_norm}) is established to give conditions ensuring that the map
\begin{displaymath}
(\tau,\nu)\longmapsto
\sqrt{\frac{1}{N(T - t_0)}\sum_{i = 1}^{N}\int_{t_0}^{T}(\tau(X_{s}^{i}) +\nu(Y_{s}^{i}))^2ds}
\textrm{ $;$ }t_0\in (0,T),
\end{displaymath}
involved in the definition of the objective function $\gamma_N$ associated to the estimator of $(a,b)$ (see (\ref{objective_function})), is actually an empirical norm on the domain of $\gamma_N$.

As the model is new in the sense we stated above, there are no previous results in the corresponding statistical setting, for the estimation of the couple of functions $(a,b)$.

We mention that Equation (\ref{main_equation}) can be seen as an extension of the following functional regression model with functional response (FRMFR):
\begin{equation}\label{FRMFR}
X_t = x_0 +\int_{0}^{T}\beta(s,t)b(Y_s)ds +\varepsilon_t,
\end{equation}
where $\beta(.,t) :=\mathbf 1_{[0,t]}(.)$ and $\varepsilon$ is a centered second order process. In addition to the relationship between $Y$ and $X$ already taken into account in (\ref{FRMFR}), Equation (\ref{main_equation}) models the own dynamics of $X$ thanks to a drift term, and to a multiplicative noise instead of $\varepsilon$. On the FRMFR, the reader can refer to Wang et al. \cite{WCM16} (see Section 3.2, p. 271-272, and references therein) for linear models, and to Febrero-Bande et al. \cite{FODA25} (see Section 3 and references therein) for non-linear ones. In \cite{FODA25}, Section 7.2, the authors compare several FRMFR to model the relationship between the electricity price $X_t$ at time $t\in [0,24]$ (one day) and the electricity demand. In this application, thanks to Equation (\ref{main_equation}), one could take into account the own dynamics of the electricity price $X$ - as in SDE-based mathematical finance (see e.g. Lamberton and Lapeyre \cite{LL07}, Chapter 4) - in addition to its relationship with the demand $Y$. Examples can also be found in the medical field, when reading Zhu et al. \cite{ZST11}, who argue that "with the advent of many high-throughput technologies, functional data are routinely collected". Their paper aims at using a diffusion model to estimate the mean function of a continuous time process. More precisely, they deal with statistical inference in the model $P_t = U_t +\varepsilon_t$, where $P$ is the Prostate Specific Antigen (PSA) level during a prostate cancer, and $U$, describing the population mean PSA process, is such that $d^{m - 1}U_t/dt^{m - 1}$ is a diffusion process.

For SDE models with no explanatory process (i.e. Model (\ref{main_equation}) with $b(.) = 0$), several nonparametric copies-based estimators of $a$ have been investigated in the literature over the last decade. For instance, Comte and Genon-Catalot \cite{CGC20} (resp. Comte and Marie \cite{CM23}) deals with a projection least squares estimator of the drift function computed from independent (resp. correlated) continuous-time observations of a diffusion process. Denis et al. \cite{DDM21} also deals with such an estimator, but computed from independent discrete-time observations of the diffusion. Moreover, Marie and Rosier \cite{MR23} deals with continuous-time and discrete-time versions of a copies-based Nadaraya-Watson estimator of the drift function of a diffusion process. On such an estimator for interacting particle systems, the reader can refer to Della Maestra and Hoffmann \cite{DMH22}, while Amorino et al. \cite{ABPPZ25} investigates the estimation of the interaction function for a class of McKean-Vlasov stochastic differential equations when both $T$ and $N$ grow to infinity. In most of these references, the authors establish a risk bound on an adaptive version of the estimator under consideration.

Let us point that in our specific model, an identifiability constraint is needed. This leads to the condition $\int b = 0$, which is taken into account in the definition of our estimator of $(a,b)$ (see (\ref{projection_LS_estimator})). The present paper deals with a projection least squares estimator of $(a,b)$ defined as a constrained minimizer of the aforementioned objective function $\gamma_N$, which is computed from the observations $(X^i,Y^i)$ ($i = 1,\dots,N$). For this estimator, a risk bound and a rate of convergence in specific Sobolev spaces are established for a fixed model. It is noteworthy that the two functions are simultaneously handled. Then, a risk bound is established on an adaptive estimator defined thanks to a Birg\'e-Massart type model selection procedure under the identifiability constraint, still for both functions simultaneously.

The plan of the paper is the following. Section \ref{section_preliminaries} provides all the probabilistic results required to define and investigate the projection least squares estimator of $(a,b)$ on the theoretical side. Risk bounds and rates of convergence for a fixed model are established in Section \ref{section_projection_LSE}, followed by a risk bound on an adaptive estimator of $(a,b)$ defined thanks to a model selection procedure. Finally, Section \ref{section_numerical_experiments} provides some numerical experiments. The results of Section \ref{section_preliminaries} (resp. Section \ref{section_projection_LSE}) are proved in Appendix \ref{section_proofs_probabilistic_results} (resp. Appendix \ref{section_proofs_statistical_results}).
\\
\\
\textbf{Notations and basic definitions:}
\begin{itemize}
 \item The usual inner product in $\mathbb R^d$ ($d\in\mathbb N^*$) is denoted by $\langle .,.\rangle_{2,\mathbb R^d}$. Moreover, for $p\geqslant 1$,
 \begin{displaymath}
 \|\mathbf x\|_{p,\mathbb R^d} :=
 \left(\sum_{\ell = 1}^{d}|\mathbf x_{\ell}|^p\right)^{\frac{1}{p}}
 \textrm{$;$ }\forall\mathbf x\in\mathbb R^d.
 \end{displaymath}
 \item For any $n,d\in\mathbb N^*$, the space of $n\times d$ real matrices is denoted by $\mathcal M_{n,d}(\mathbb R)$, and equipped with the operator norm $\|.\|_{\rm op}$ defined by
 \begin{displaymath}
 \|M\|_{\rm op} :=
 \sup_{\mathbf x\in\mathbb R^d}\frac{\|M\mathbf x\|_{2,\mathbb R^n}}{\|\mathbf x\|_{2,\mathbb R^d}}
 \textrm{ $;$ }\forall M\in\mathcal M_{n,d}(\mathbb R).
 \end{displaymath}
 Moreover, for every $M\in\mathcal M_{n,d}(\mathbb R)$, the transpose (resp. the trace) of $M$ is denoted by $M^*$ (resp. ${\rm Tr}(M)$).
 \item The usual inner product in $\mathbb L^2(\mathbb R;\mathbb R^d)$ (resp. $\mathbb L^2([0,T];\mathbb R^d)$) is denoted by $\langle .,.\rangle$ (resp. $\langle .,.\rangle_T$).
 \item In the sequel,
 \begin{displaymath}
 \mathcal C :=\left\{H
 \textrm{, $\mathbb F$-adapted continuous process} :
 \mathbb E\left(\sup_{t\in [0,T]}H_{t}^{2}\right) <\infty\right\}
 \end{displaymath}
 is equipped with the norm $\|.\|_{\mathcal C}$ defined by
 \begin{displaymath}
 \|H\|_{\mathcal C} :=
 \mathbb E\left(\sup_{t\in [0,T]}H_{t}^{2}\right)^{\frac{1}{2}}
 \textrm{$;$ }\forall H\in\mathcal C.
 \end{displaymath}
 The normed vector space $(\mathcal C,\|.\|_{\mathcal C})$ is complete.
 \item The Malliavin derivative of order $k\in\mathbb N^*$ - with respect to $W = (W_1,W_2)$ - is denoted by
 \begin{displaymath}
 \mathbf D^{(k)} = (\mathbf D^{(k),\ell_1,\dots,\ell_k})_{\ell_1,\dots,\ell_k\in\{1,2\}}.
 \end{displaymath}
 Moreover, for every $p\geqslant 2$, the domain of $\mathbf D^{(k)}$ in $\mathbb L^p(\Omega)$ is denoted by $\mathbb D^{k,p}$. Finally, for the sake of readability, $\mathbf D^{(1)}$ (resp. $\mathbf D^{(1),\ell}$ with $\ell\in\{1,2\}$) is denoted by $\mathbf D$ (resp. $\mathbf D^{\ell}$).
 \item Consider
 \begin{displaymath}
 \mathcal C_p :=\left\{H
 \textrm{, $\mathbb F$-adapted continuous process} :
 \mathbb E\left(\sup_{t\in [0,T]}|H_t|^p\right) <\infty\right\}
 \textrm{ $;$ }p\geqslant 1,
 \end{displaymath}
 and then
 \begin{displaymath}
 \mathbb H^{\infty} :=
 \bigcap_{p\geqslant 1}
 \left\{H\in\mathcal C_p :
 \forall k\in\mathbb N^*,
 \sup_{s_1,\dots,s_k\in [0,T]}\mathbb E\left(
 \sup_{t\in [s_1\vee\dots\vee s_k,T]}
 \|\mathbf D_{s_1,\dots,s_k}^{(k)}H_t\|_{p,\mathbb R^{2^k}}^{p}\right)
 <\infty\right\}.
 \end{displaymath}
 \item For every measurable function $f$ from $[0,T]$ into $\mathbb R^d$, $f\not\equiv 0$ means that $\lambda\{t\in [0,T] : f(t)\neq 0\} > 0$.
\end{itemize}
%


%
\section{Probabilistic preliminaries}\label{section_preliminaries}
This section deals with conditions that $(X,Y)$ needs to satisfy for our statistical purposes. 
%


%
\subsection{First properties of $(X_t,Y_t)$ and its distribution}\label{section_first_properties_XY}
From now on, we assume that Assumption \ref{main_assumption} holds. In particular, it ensures the following result.
%


%
\begin{proposition}\label{existence_uniqueness_solution} Under Assumption \ref{main_assumption}, Equation (\ref{main_equation}) has a unique (strong) solution in $\mathcal C$.
\end{proposition}
\noindent
Proposition \ref{existence_uniqueness_solution} is a straightforward consequence of Nualart \cite{NUALART06}, Lemma 2.2.1. In the sequel, $a$, $b$, $\sigma$ and $Y$ fulfill the two following assumptions.
%


%
\begin{assumption}\label{assumption_drift_volatility}
The functions $a$, $b$ and $\sigma$ are continuously differentiable from $\mathbb R$ into itself, $a'$, $b'$, $\sigma$ and $\sigma'$ are bounded, and $\inf_{\mathbb R}|\sigma| > 0$.
\end{assumption}
%


%
\begin{assumption}\label{assumption_Y}
First,
\begin{displaymath}
\mathbb E\left(\sup_{t\in [0,T]}Y_{t}^{2}\right) +
\sup_{s\in [0,T]}\mathbb E\left(\sup_{t\in [s,T]}
(\mathbf D_{s}^{2}Y_t)^2\right) <\infty.
\end{displaymath}
Moreover, for every $t\in (0,T]$,
\begin{enumerate}
 \item $\mathbf D^2Y_t\not\equiv 0$, and
 \item the distribution of $Y_t$ has a positive and continuously differentiable density $f_{Y,t}$ with respect to the Lebesgue measure on $\mathbb R$.
\end{enumerate}
\end{assumption}
\noindent
First, in order to ensure that the empirical norm in Section \ref{section_projection_LSE} is positive definite, let us establish the following technical result.
%


%
\begin{proposition}\label{technical_result_empirical_norm}
Consider $\mathcal S_1,\mathcal S_2\subset C^1(\mathbb R)$. Under Assumptions \ref{assumption_drift_volatility} and \ref{assumption_Y}, for every $\tau\in\mathcal S_1$, and every
\begin{equation}\label{technical_result_empirical_norm_1}
\nu\in\mathcal S_2
\quad\textrm{such that}\quad\int_{-\infty}^{\infty}\nu(y)dy = 0,
\end{equation}
if $\tau(X_t) +\nu(Y_t) = 0$ for every $t\in (0,T]$, then $\nu(.) = 0$, and $\tau(X_t) = 0$ for every $t\in (0,T]$.
\end{proposition}
\noindent
Now, for our statistical purposes, $(X,Y)$ needs to fulfill the following assumption.
%


%
\begin{assumption}\label{assumption_XY}
For every $t\in (0,T]$, the distribution of $(X_t,Y_t)$ has a density $f_t$ with respect to the Lebesgue measure on $\mathbb R^2$. Moreover, there exists $t_0\in [0,T)$ such that, for every $(x,y)\in\mathbb R^2$, the map $s\mapsto f_s(x,y)$ belongs to $\mathbb L^1([t_0,T])$.
\end{assumption}
\noindent
Of course, by assuming that $Y$ is a diffusion process, well-known results on multidimensional stochastic differential equations allow to establish that $Y$ and $(X,Y)$ fulfill Assumptions \ref{assumption_Y} and \ref{assumption_XY} (see Proposition \ref{sufficient_condition_AC_diffusion}). However, since $Y$ doesn't need to be a diffusion process in Section \ref{section_projection_LSE}, the following proposition provides a general sufficient condition on $a$, $b$, $\sigma$ and $Y$ for $(X,Y)$ to fulfill Assumption \ref{assumption_XY}.
%


%
\begin{proposition}\label{general_sufficient_condition_AC}
Consider $t\in (0,T]$, and assume that $a$, $b$ and $\sigma$ fulfill Assumption \ref{assumption_drift_volatility}.
\begin{enumerate}
 \item Under Assumption \ref{assumption_Y}.(1), the distribution of $(X_t,Y_t)$ has a density with respect to the Lebesgue measure on $\mathbb R^2$.
 \item Assume that $\sigma$ is bounded, and that $a,b,\sigma\in C^{\infty}(\mathbb R)$ with all derivatives bounded. Assume also that $Y\in\mathbb H^{\infty}$, and that $1/\|\mathbf D^2Y_t\|_T$ belongs to $\mathbb L^p(\Omega)$ for every $p\geqslant 1$. Then, the distribution of $(X_t,Y_t)$ has a smooth (i.e. infinitely differentiable) density $f_t$ with respect to the Lebesgue measure on $\mathbb R^2$, and there exist $\mathfrak c_{\ref{general_sufficient_condition_AC}} > 0$ and $(m,\alpha)\in\mathbb N^*\times (2,\infty)$, depending on $T$ but not on $t$, such that for every $(x,y)\in\mathbb R^2$,
 \begin{displaymath}
 f_t(x,y)\leqslant
 \frac{\mathfrak c_{\ref{general_sufficient_condition_AC}}}{t^m}
 \mathbb E\left(\frac{1}{\|\mathbf D^2Y_t\|_{T}^{4\alpha}}\right)^{\frac{m}{2\alpha}}
 \Pi_t(x,y)
 \end{displaymath}
 with
 \begin{displaymath}
 \Pi_t(x,y) =
 \mathbb P(|Y_t| > |y|)^{\frac{1}{2}}
 \mathbb P(|X_t - x_0| > |x - x_0|)^{\frac{1}{2}}.
 \end{displaymath}
\end{enumerate}
\end{proposition}
%


%
\subsection{A class of processes $Y$ fulfilling the conditions of Proposition \ref{general_sufficient_condition_AC}}\label{example_explanatory_processes_I}
This section deals with a class of processes $Y$ fulfilling the conditions of Proposition \ref{general_sufficient_condition_AC}, but not defined by a stochastic differential equation. Explicit processes are provided in conclusion.
\\
\\
Consider $Y = g\circ H$, where $g\in C^{\infty}(\mathbb R)$, $g^{(k)}$ has polynomial growth for every $k\in\mathbb N$, and
\begin{displaymath}
H :=\int_{0}^{.}h(u)dW_2(u)
\quad\textrm{with}\quad
h\in C^0([0,T]).
\end{displaymath}
First, assume that
\begin{displaymath}
{\rm (A)}\quad h\mathbf 1_{[0,t]}\not\equiv 0
\textrm{ $;$ }\forall t\in (0,T]
\quad\textrm{and}\quad
{\rm (B)}\quad
|g'(u)| > 0\textrm{ $;$ }\forall u\in\mathbb R.
\end{displaymath}
Under Assumption \ref{assumption_drift_volatility}, the distribution of $(X_t,Y_t)$ ($t\in (0,T]$) has a density with respect to the Lebesgue measure on $\mathbb R^2$.
\\
\\
Indeed, recursively, for every $t\in [0,T]$, $k\in\mathbb N^*$, $\ell_1,\dots,\ell_k\in\{1,2\}$ and $s_1,\dots,s_k\in [0,T]$,
\begin{displaymath}
\mathbf D_{s_1,\dots,s_k}^{(k),\ell_1,\dots,\ell_k}Y_t =
\left\{
\begin{array}{ll}
 h(s_1)\cdots h(s_k)g^{(k)}(H_t)
 & \textrm{when $t\geqslant s_1\vee\dots\vee s_k$ and $\ell_1,\dots,\ell_k = 2$}\\
 0 & \textrm{otherwise}
\end{array}\right..
\end{displaymath}
Then, since $h\in C^0([0,T])$, since $g^{(k)}$ has polynomial growth for every $k\in\mathbb N$, and by the Burkholder-Davis-Gundy inequality, $Y\in\mathbb H^{\infty}$. Moreover, the conditions ${\rm (A)}$ and ${\rm (B)}$ lead to $\mathbf D^2Y_t\not\equiv 0$ for any $t\in (0,T]$. So, the distribution of $(X_t,Y_t)$ has a density with respect to the Lebesgue measure on $\mathbb R^2$ by Proposition \ref{general_sufficient_condition_AC}.(1).
\\
\\
Now, assume that $\sigma$ is bounded, $a,b,\sigma\in C^{\infty}(\mathbb R)$ with all derivatives bounded, $h$ satisfies ${\rm (A)}$, and that - instead of the condition ${\rm (B)}$ - $g$ satisfies
\begin{displaymath}
{\rm (\overline B)}\quad
\mathfrak g :=\inf_{u\in\mathbb R}|g'(u)| > 0.
\end{displaymath}
Under Assumption \ref{assumption_drift_volatility}, $Y$ and $(X,Y)$ fulfill Assumptions \ref{assumption_Y} and \ref{assumption_XY}.
\\
\\
Indeed, for every $t\in (0,T]$ and $p\geqslant 1$, by the conditions ${\rm (A)}$ and ${\rm (\overline B)}$,
\begin{displaymath}
\mathbb E\left(\frac{1}{\|\mathbf D^2Y_t\|_{T}^{p}}\right)
\leqslant\left(\frac{1}{\mathfrak g\|h\mathbf 1_{[0,t]}\|_T}\right)^p <\infty.
\end{displaymath}
So, by Proposition \ref{general_sufficient_condition_AC}.(2), the distribution of $(X_t,Y_t)$ has a smooth density $f_t$ with respect to the Lebesgue measure on $\mathbb R^2$, and for any $(x,y)\in\mathbb R^2$,
\begin{equation}\label{example_explanatory_processes_I_2}
f_s(x,y)\leqslant
\frac{\mathfrak c_{\ref{general_sufficient_condition_AC}}}{
(\mathfrak g\|h\mathbf 1_{[0,s]}\|_T)^{2m}s^m}
\textrm{ $;$ }\forall s\in (0,T].
\end{equation}
Note that by Inequality (\ref{example_explanatory_processes_I_2}), $s\mapsto f_s(x,y)$ belongs to $\mathbb L^1([t_0,T])$ for any $t_0\in (0,T)$. Finally, the condition ${\rm (\overline B)}$ also says that $g$ is one-to-one from $\mathbb R$ into itself, and since $H$ is a non-degenerate Gaussian process, $f_{Y,t}(y) > 0$ for every $t\in (0,T]$ and $y\in\mathbb R$. Therefore, $Y$ and $(X,Y)$ fulfill Assumptions \ref{assumption_Y} and \ref{assumption_XY}.
\\
\\
Let us conclude with two explicit examples of explanatory processes satisfying both conditions ${\rm (A)}$ and ${\rm (\overline B)}$:
\begin{enumerate}
 \item If $Y = W_2(1 + W_{2}^{2})$, then $h = 1$ ($\Rightarrow {\rm (A)}$), and for every $u\in\mathbb R$, $g(u) = u(1 + u^2)$, leading to $g'(u) = 1 + 3u^2\geqslant 1$ ($\Rightarrow {\rm (\overline B)}$).
 \item If $Y = H +\arctan\circ H$ with $dH_t = tdW_2(t)$, then $h = {\rm Id}_{[0,T]}$ ($\Rightarrow {\rm (A)}$), and for every $u\in\mathbb R$, $g(u) = u +\arctan(u)$, leading to $g'(u) = 1 + (1 + u^2)^{-1}\geqslant 1$ ($\Rightarrow {\rm (\overline B)}$).
\end{enumerate}
%


%
\subsection{Integrability results on the map $(s,x,y)\mapsto f_s(x,y)$}\label{section_integrability_density_XY}
Under additional conditions on $(a,b)$ and $Y$, the following proposition provides integrability results on the map $(s,x,y)\mapsto f_s(x,y)$.
%


%
\begin{proposition}\label{integrability_Pi}
Assume that $a$ and $b$ are bounded. Under the assumptions of Proposition \ref{general_sufficient_condition_AC}.(2),
\begin{enumerate}
 \item There exist $\mathfrak c_{\ref{integrability_Pi}},\mathfrak m_{\ref{integrability_Pi}} > 0$ such that, for every $t\in (0,T]$ and $(x,y)\in\mathbb R^2$,
 \begin{displaymath}
 \Pi_t(x,y)\leqslant
 \mathfrak c_{\ref{integrability_Pi}}
 \mathbb P(|Y_t| > |y|)^{\frac{1}{2}}
 \exp\left(-\mathfrak m_{\ref{integrability_Pi}}\frac{(x - x_0)^2}{t}\right).
 \end{displaymath}
 \item With the notations of Proposition \ref{general_sufficient_condition_AC}.(2), assume that there exist $t_0\in (0,T)$ and $\mu_{t_0,T}\in\mathbb L^1(\mathbb R;\mathbb R_+)$ such that $\sup_{\mathbb R}\mu_{t_0,T} <\infty$ and
 \begin{equation}\label{integrability_Pi_1}
 \sup_{t\in [t_0,T]}\left\{
 \mathbb E\left(\frac{1}{\|\mathbf D^2Y_t\|_{T}^{4\alpha}}\right)^{\frac{m}{2\alpha}}
 \mathbb P(|Y_t| > |y|)^{\frac{1}{2}}\right\}
 \leqslant\mu_{t_0,T}(y)
 \textrm{ $;$ }\forall y\in\mathbb R.
 \end{equation}
 Then,
 \begin{equation}\label{integrability_Pi_2}
 \sup_{x\in\mathbb R}\int_{[t_0,T]\times\mathbb R}f_s(x,y)dsdy <\infty
 \quad\textrm{and}\quad
 \sup_{y\in\mathbb R}\int_{[t_0,T]\times\mathbb R}f_s(x,y)dsdx <\infty.
 \end{equation}
\end{enumerate}
\end{proposition}
%


%
\begin{example}\label{example_explanatory_processes_II}
As in Section \ref{example_explanatory_processes_I}, $Y = g\circ H$, where $g\in C^{\infty}(\mathbb R)$, $g^{(k)}$ has polynomial growth for every $k\in\mathbb N$,
\begin{displaymath}
H :=\int_{0}^{.}h(u)dW_2(u)
\quad\textrm{with}\quad
h\in C^0([0,T]),
\end{displaymath}
and $h$ (resp. $g$) satisfies the condition $({\rm A})$ (resp. $(\overline{\rm B})$). First, for any $t_0\in (0,T)$, by the condition $(\overline{\rm B})$ on $g$,
\begin{displaymath}
\sup_{t\in [t_0,T]}
\mathbb E\left(\frac{1}{\|\mathbf D^2Y_t\|_{T}^{4\alpha}}\right)^{\frac{m}{2\alpha}}
\leqslant\left(\frac{1}{\mathfrak g\|h\mathbf 1_{[0,t_0]}\|_T}\right)^{2m}
\quad\textrm{with}\quad
\mathfrak g =\inf_{u\in\mathbb R}|g'(u)|.
\end{displaymath}
Now, for the sake of simplicity, assume that
\begin{equation}\label{example_explanatory_processes_II_1}
g(-u) = -g(u)
\quad\textrm{and}\quad
g'(u) > 0
\textrm{ $;$ }\forall u\in\mathbb R.
\end{equation}
For instance, the functions $u\mapsto u(1 + u^2)$ and $u\mapsto u +\arctan(u)$, satisfying $(\overline{\rm B})$ as mentioned at the end of Section \ref{example_explanatory_processes_I}, also fulfill (\ref{example_explanatory_processes_II_1}). For every $t\in [t_0,T]$,
\begin{displaymath}
H_t\rightsquigarrow\mathcal N(0,\sigma_{t}^{2})
\quad\textrm{with}\quad
\sigma_t =\|h\mathbf 1_{[0,t]}\|_T,
\end{displaymath}
leading to
\begin{eqnarray*}
 \mathbb P(|Y_t| > |y|)
 & = &
 \mathbb P(g(H_t)\in\mathbb R\backslash [-|y|,|y|]) =
 \frac{1}{\sigma_t\sqrt{2\pi}}
 \int_{\mathbb R\backslash [-g^{-1}(|y|),g^{-1}(|y|)]}
 \exp\left(-\frac{u^2}{2\sigma_{t}^{2}}\right)du\\
 & \leqslant &
 \frac{1}{\sigma_t\sqrt{2\pi}}
 \exp\left(-\frac{g^{-1}(|y|)^2}{4\sigma_{t}^{2}}\right)
 \int_{-\infty}^{\infty}
 \exp\left(-\frac{u^2}{4\sigma_{t}^{2}}\right)du =
 \sqrt 2\exp\left(-\frac{g^{-1}(|y|)^2}{4\sigma_{t}^{2}}\right).
\end{eqnarray*}
Thus,
\begin{displaymath}
\sup_{t\in [t_0,T]}\left\{
\mathbb E\left(\frac{1}{\|\mathbf D^2Y_t\|_{T}^{4\alpha}}\right)^{\frac{m}{2\alpha}}
\mathbb P(|Y_t| > |y|)^{\frac{1}{2}}\right\}
\leqslant
2^{\frac{1}{4}}
\left(\frac{1}{\mathfrak g\sigma_{t_0}}\right)^{2m}
\exp\left(-\frac{g^{-1}(|y|)^2}{8\sigma_{T}^{2}}\right) =:\mu_{t_0,T}(y).
\end{displaymath}
Finally, since $g'$ has polynomial growth,
\begin{eqnarray*}
 \int_{-\infty}^{\infty}|\mu_{t_0,T}(y)|dy
 & = &
 \mathfrak c_1\int_{0}^{\infty}\exp\left(-\frac{g^{-1}(y)^2}{8\sigma_{T}^{2}}\right)dy
 \quad\textrm{with}\quad
 \mathfrak c_1 =
 2^{\frac{5}{4}}
 \left(\frac{1}{\mathfrak g\sigma_{t_0}}\right)^{2m}\\
 & = &
 \mathfrak c_1
 \int_{0}^{\infty}\exp\left(-\frac{z^2}{8\sigma_{T}^{2}}\right)g'(z)dz <\infty,
\end{eqnarray*}
and by Proposition \ref{integrability_Pi}.(2),
\begin{displaymath}
\sup_{x\in\mathbb R}\int_{[t_0,T]\times\mathbb R}f_s(x,y)dsdy <\infty
\quad\textrm{and}\quad
\sup_{y\in\mathbb R}\int_{[t_0,T]\times\mathbb R}f_s(x,y)dsdx <\infty.
\end{displaymath}
\end{example}
\noindent
Finally, assume that $Y$ is defined by the stochastic differential equation
\begin{equation}\label{Y_diffusion}
Y_t = y_0 +
\int_{0}^{t}\mu(Y_s)ds +
\int_{0}^{t}\kappa(Y_s)dW_2(s)
\textrm{ $;$ }t\in [0,T],
\end{equation}
where $y_0\in\mathbb R$, $\mu,\kappa\in C^1(\mathbb R)$, and $\mu'$ and $\kappa'$ are bounded. In this situation, the following proposition provides a sufficient condition on $\sigma$ and $\kappa$ for $Y$ and $(X,Y)$ to fulfill Assumptions \ref{assumption_Y} and \ref{assumption_XY}, and to satisfy (\ref{integrability_Pi_2}).
%


%
\begin{proposition}\label{sufficient_condition_AC_diffusion}
Let $Y$ be the solution of Equation (\ref{Y_diffusion}), and assume that $a$, $b$ and $\sigma$ fulfill Assumption \ref{assumption_drift_volatility}. If $\kappa$ is bounded, and if $\inf_{\mathbb R}|\kappa| > 0$, then for any $t_0\in (0,T)$, $Y$ and $(X,Y)$ fulfill Assumptions \ref{assumption_Y} and \ref{assumption_XY}, and satisfy (\ref{integrability_Pi_2}).
\end{proposition}
%


%
\section{A projection least squares estimator of $(a,b)$}\label{section_projection_LSE}
%


%
\subsection{Definition of the estimator}\label{section_estimator}
This section deals with the identifiability of the model, and then with the definition of the projection least squares estimator of $(a\mathbf 1_{A_1},b\mathbf 1_{A_2})$, where $A_1$ and $A_2$ are intervals of $\mathbb R$. All the notations used in the sequel have been introduced at the end of the introduction section.
\\
\\
Clearly, the functions $a$ and $b$ are defined up to an additive constant. So, let us assume that
\begin{equation}\label{identifiability_condition}
\int_{A_2}b(y)dy = 0,
\quad
\textrm{which is usual in additive models.}
\end{equation}
\newline
We intend to build estimators of $a$ and $b$ by looking for their coefficients in projection spaces. So, consider $m_1,m_2\in\mathbb N^*$ such that $m_1,m_2\leqslant N$, and
\begin{displaymath}
S_{m_1} := {\rm span}(\varphi_1,\dots,\varphi_{m_1})
\quad (\textrm{resp. }
\Sigma_{m_2} := {\rm span}(\psi_1,\dots,\psi_{m_2})),
\end{displaymath}
where $(\varphi_j)_{j\in\mathbb N^*}$ (resp. $(\psi_k)_{k\in\mathbb N^*}$) is an orthonormal family of $\mathbb L^2(A_1)$ (resp. $\mathbb L^2(A_2)$), which elements are assumed to be continuously differentiable functions from $A_1$ (resp. $A_2$) into $\mathbb R$. For any $(\tau,\nu)$ belonging to $\mathcal S_{\bf m} := S_{m_1}\times\mathbb S_{m_2}$ with
\begin{displaymath}
\mathbf m = (m_1,m_2)
\quad {\rm and}\quad
\mathbb S_{m_2} =
\left\{\nu\in\Sigma_{m_2} :\int_{A_2}\nu(y)dy = 0\right\},
\end{displaymath}
there exist $(t_1,\dots,t_{m_1})\in\mathbb R^{m_1}$ and $(n_1,\dots,n_{m_2})\in\mathbb R^{m_2}$ such that
\begin{displaymath}
\tau =\sum_{j = 1}^{m_1}t_j\varphi_j
\quad {\rm and}\quad
\nu =\sum_{k = 1}^{m_2}n_k\psi_k
\quad\textrm{with the constraint}\quad
\sum_{k = 1}^{m_2}n_k\int_{A_2}\psi_k(y)dy = 0.
\end{displaymath}
Our purpose is to estimate $t_1,\dots,t_{m_1},n_1,\dots,n_{m_2}$ for $(\tau,\nu)$ to be as close as possible to $(a,b)$.
\\
\\
First, let $t_0\in (0,T)$ be a fixed time, set $T_0 := T - t_0$, and consider the objective function $\gamma_N$ such that, for every $(\tau,\nu)\in\mathcal S_{\bf m}$,
\begin{equation}\label{objective_function}
\gamma_N(\tau,\nu) :=
\frac{1}{NT_0}\sum_{i = 1}^{N}\left[\int_{t_0}^{T}(\tau(X_{s}^{i}) +\nu(Y_{s}^{i}))^2ds -
2\int_{t_0}^{T}(\tau(X_{s}^{i}) +\nu(Y_{s}^{i}))dX_{s}^{i}\right].
\end{equation}
The first term of the sum in the right-hand side of Equality (\ref{objective_function}) suggests to define an empirical norm by
\begin{displaymath}
\|(\tau,\nu)\|_{N}^{2} :=
\frac{1}{NT_0}\sum_{i = 1}^{N}\int_{t_0}^{T}(\tau(X_{s}^{i}) +\nu(Y_{s}^{i}))^2ds.
\end{displaymath}
%


%
\begin{remark}\label{positive_definiteness_empirical_norm}
Under Assumptions \ref{assumption_drift_volatility} and \ref{assumption_Y}, one can provide some conditions on both the $\varphi_j$'s and $Y$ for $\|.\|_N$ to be positive definite on $\mathcal S_{\bf m}$. Consider $(\tau,\nu)\in\mathcal S_{\bf m}$ such that $\|(\tau,\nu)\|_N = 0$. In particular, for every $s\in [t_0,T]$, $\tau(X_{s}^{1}) +\nu(Y_{s}^{1}) = 0$, and then $\nu(.) = 0$ and $\tau(X_{s}^{1}) = 0$ by Proposition \ref{technical_result_empirical_norm}. Since $\tau\in S_{m_1}$, there exist $t_1,\dots,t_{m_1}\in\mathbb R$ such that $\tau = t_1\varphi_1 +\dots + t_{m_1}\varphi_{m_1}$, and since the paths of $X^1$ are continuous but not constant on $[t_0,T]$, $I := X^1([t_0,T])$ is a (random) compact interval of $\mathbb R$ such that $\lambda(I(\omega)) > 0$ for every $\omega\in\Omega$. On the one hand, if $(\varphi_1,\dots,\varphi_{m_1})$ is a $\mathbb R$-supported basis (e.g. Hermite's basis), then $(\varphi_1)_{|I(\omega)},\dots,(\varphi_{m_1})_{|I(\omega)}$ ($\omega\in\Omega$) are non-zero linearly independent vectors of $\mathbb L^2(I(\omega))$, leading to
\begin{displaymath}
\tau_{|I(\omega)} =
\sum_{j = 1}^{m_1}t_j(\varphi_j)_{|I(\omega)} = 0\Longrightarrow
t_1 =\dots = t_{m_1} = 0.
\end{displaymath}
On the other hand, at least when $Y$ is defined by Equation (\ref{Y_diffusion}), thanks to the positive lower bound on $(s,x,y)\in (0,T]\times\mathbb R^2\mapsto f_s(x,y)$ in Theorem 1.2 of Menozzi et al. \cite{MPZ21}, there exists $\omega\in\Omega$ such that $J := I(\omega)\cap A_1$ is a compact interval of $\mathbb R$ satisfying $\lambda(J) > 0$. Then, as previously, $(\varphi_1)_{|J},\dots,(\varphi_{m_1})_{|J}$ are non-zero linearly independent vectors of $\mathbb L^2(J)$, leading to $t_1 =\dots = t_{m_1} = 0$.
\end{remark}
\noindent
Moreover, consider the $(m_1 + m_2)\times (m_1 + m_2)$ random matrix
\begin{displaymath}
\widehat\Psi_{\bf m} :=
\begin{pmatrix}
 \widehat\Psi_{1,1} & \widehat\Psi_{1,2}\\
 \widehat\Psi_{1,2}^{*} & \widehat\Psi_{2,2}
\end{pmatrix},
\end{displaymath}
where
\begin{eqnarray*}
 \widehat\Psi_{1,1} & := &
 \left(\frac{1}{NT_0}\sum_{i = 1}^{N}\int_{t_0}^{T}
 \varphi_j(X_{s}^{i})\varphi_{j'}(X_{s}^{i})ds\right)_{1\leqslant j,j'\leqslant m_1},\\
 \widehat\Psi_{2,2} & := &
 \left(\frac{1}{NT_0}\sum_{i = 1}^{N}\int_{t_0}^{T}
 \psi_k(Y_{s}^{i})\psi_{k'}(Y_{s}^{i})ds\right)_{1\leqslant k,k'\leqslant m_2}\quad {\rm and}\\
 \widehat\Psi_{1,2} & := &
 \left(\frac{1}{NT_0}\sum_{i = 1}^{N}\int_{t_0}^{T}
 \varphi_j(X_{s}^{i})\psi_k(Y_{s}^{i})ds\right)_{(j,k)\in\{1,\dots,m_1\}\times\{1,\dots,m_2\}}.
\end{eqnarray*}
The matrix $\widehat\Psi_{\bf m}$ is related with the empirical norm $\|.\|_N$ in the following way: for every vector $\mathbf x = (x_1,\dots,x_{m_1 + m_2})$ of $\mathbb R^{m_1 + m_2}$,
\begin{eqnarray*}
 \mathbf x^*\widehat\Psi_{\bf m}\mathbf x & = &
 \left\|\left(\sum_{j = 1}^{m_1}x_j\varphi_j,
 \sum_{k = 1}^{m_2}x_{m_1 + k}\psi_k\right)\right\|_{N}^{2}\\
 & = &
 \|(\tau,\nu)\|_{N}^{2}\geqslant 0
 \quad {\rm with}\quad
 \tau =\sum_{j = 1}^{m_1}x_j\varphi_j
 \quad {\rm and}\quad
 \nu =\sum_{k = 1}^{m_2}x_{m_1 + k}\psi_k.
\end{eqnarray*}
Then, the symmetric matrix $\widehat\Psi_{\bf m}$ is positive semidefinite, and when $\|.\|_N$ is a norm on $\mathcal S_{\bf m}$ (see Remark \ref{positive_definiteness_empirical_norm}), $\widehat\Psi_{\bf m}$ is even positive definite in $\mathcal M_{m_1,m_2}(\mathbb X_{\bf m})$ with $\mathbb X_{\bf m}\subset\mathbb R^{m_1 + m_2}$ satisfying $\mathbb X_{\bf m}\cong\mathcal S_{\bf m}$.
\\
\\
Now, assume that $\widehat\Psi_{\bf m}$ is invertible, and let $(\widehat a_{m_1},\widehat b_{m_2})$ be the minimizer of $\gamma_N$ over $\mathcal S_{\bf m}$. Precisely,
\begin{displaymath}
\widehat a_{m_1} =\sum_{j = 1}^{m_1}\widehat\theta_j\varphi_j
\quad {\rm and}\quad
\widehat b_{m_2} =\sum_{k = 1}^{m_2}\widehat\theta_{m_1 + k}\psi_k
\end{displaymath}
with
\begin{equation}\label{projection_LS_estimator}
\widehat\theta =
\underset{\theta\in\mathbb R^{m_1 + m_2} : h(\theta) = 0}{\rm argmin}\mathcal J_N(\theta),
\end{equation}
where
\begin{eqnarray*}
 \mathcal J_N(\theta) & := &
 \frac{1}{NT_0}\sum_{i = 1}^{N}\left[\int_{t_0}^{T}\left(
 \sum_{j = 1}^{m_1}\theta_j\varphi_j(X_{s}^{i}) +
 \sum_{k = 1}^{m_2}\theta_{m_1 + k}\psi_k(Y_{s}^{i})\right)^2ds\right.\\
 & &
 \hspace{4cm} -\left.
 2\int_{t_0}^{T}\left(
 \sum_{j = 1}^{m_1}\theta_j\varphi_j(X_{s}^{i}) +
 \sum_{k = 1}^{m_2}\theta_{m_1 + k}\psi_k(Y_{s}^{i})\right)dX_{s}^{i}\right]
\end{eqnarray*}
and
\begin{displaymath}
h(\theta) :=
\sum_{k = 1}^{m_2}\theta_{m_1 + k}\int_{A_2}\psi_k(y)dy =
\langle\theta,\mathbf d_{\bf m}\rangle_{2,\mathbb R^{m_1 + m_2}}
\end{displaymath}
with
\begin{equation}\label{constraint}
\mathbf d_{\bf m} =
(0,\dots,0,\delta_{m_2})^*
\quad {\rm and}\quad
\delta_{m_2} =
\left(\int_{A_2}\psi_1(y)dy,\dots,\int_{A_2}\psi_{m_2}(y)dy\right).
\end{equation}
Consider
\begin{displaymath}
\widehat{\bf Z}_{\bf m} :=
\begin{bmatrix}
 \displaystyle{\left(\frac{1}{NT_0}\sum_{i = 1}^{N}\int_{t_0}^{T}
 \varphi_j(X_{s}^{i})dX_{s}^{i}\right)_{1\leqslant j\leqslant m_1}}\\
 \displaystyle{\left(\frac{1}{NT_0}\sum_{i = 1}^{N}\int_{t_0}^{T}
 \psi_k(Y_{s}^{i})dX_{s}^{i}\right)_{1\leqslant k\leqslant m_2}}
\end{bmatrix},
\end{displaymath}
and let $\mathcal L_N$ be the Lagrangian for Problem (\ref{projection_LS_estimator}):
\begin{displaymath}
\mathcal L_N(\theta,\lambda) :=
\mathcal J_N(\theta) -\lambda h(\theta)\textrm{ $;$ }
(\theta,\lambda)\in\mathbb R^{m_1 + m_2}\times\mathbb R.
\end{displaymath}
Necessarily,
\begin{displaymath}
\nabla\mathcal L_N(\widehat\theta,\widehat\lambda) =
\begin{pmatrix}
 2(\widehat\Psi_{\bf m}\widehat\theta -
 \widehat{\bf Z}_{\bf m}) -\widehat\lambda\mathbf d_{\bf m}\\
 -h(\widehat\theta)
\end{pmatrix} = 0.
\end{displaymath}
So, if $\mathbf d_m\neq 0$ (or equivalently if $\delta_{m_2}\neq 0$), then
\begin{displaymath}
\widehat\theta =
\widehat\Psi_{\bf m}^{-1}\left(\widehat{\bf Z}_{\bf m} +
\frac{\widehat\lambda}{2}\mathbf d_{\bf m}\right),
\quad\textrm{leading to}\quad
\widehat\lambda =
-2\cdot\frac{\langle
\widehat\Psi_{\bf m}^{-1}\widehat{\bf Z}_{\bf m},
\mathbf d_{\bf m}\rangle_{2,\mathbb R^{m_1 + m_2}}}{
\langle
\widehat\Psi_{\bf m}^{-1}\mathbf d_{\bf m},\mathbf d_{\bf m}
\rangle_{2,\mathbb R^{m_1 + m_2}}}.
\end{displaymath}
Therefore,
\begin{equation}\label{coefficients_pLS_estimator}
\widehat\theta =
\widehat\Psi_{\bf m}^{-1}\widehat{\bf Z}_{\bf m} -
\frac{\mathbf d_{\bf m}^{*}\widehat\Psi_{\bf m}^{-1}\widehat{\bf Z}_{\bf m}}{
\mathbf d_{\bf m}^{*}\widehat\Psi_{\bf m}^{-1}\mathbf d_{\bf m}}\cdot
\widehat\Psi_{\bf m}^{-1}\mathbf d_{\bf m}
\quad {\rm when}\quad
\delta_{m_2}\neq 0.
\end{equation}
Clearly, if $\delta_{m_2} = 0$, then the constraint is automatically satisfied and vanishes, leading to $\widehat\theta =\widehat\Psi_{\bf m}^{-1}\widehat{\bf Z}_{\bf m}$.
\\
\\
Finally, under Assumption \ref{assumption_XY} - satisfied when $Y$ fulfills Assumption \ref{assumption_Y} (resp. $Y$ is the solution of Equation (\ref{Y_diffusion})) by Proposition \ref{general_sufficient_condition_AC} (resp. Proposition \ref{sufficient_condition_AC_diffusion}) - the distribution of $(X_t,Y_t)$ ($t\in [t_0,T]$) has a density $f_t$ with respect to the Lebesgue measure on $\mathbb R^2$ and, for every $(x,y)\in\mathbb R^2$, the map $s\mapsto f_s(x,y)$ belongs to $\mathbb L^1([t_0,T])$. This legitimates to consider the density function $f$ defined by
\begin{displaymath}
f(x,y) :=\frac{1}{T_0}\int_{t_0}^{T}f_s(x,y)ds
\textrm{ $;$ }\forall (x,y)\in\mathbb R^2.
\end{displaymath}
Then, we can explain why minimizing the objective fonction $\gamma_N$ is meaningful in order to estimate the $\mathbb R^2$-valued function $(a,b)$. Indeed, for every $(\tau,\nu)\in\mathcal S_{\bf m}$,
\begin{eqnarray}
 \mathbb E(\gamma_N(\tau,\nu)) & = &
 \mathbb E\left[\int_{t_0}^{T}(\tau(X_s) +\nu(Y_s) - (a(X_s) + b(Y_s)))^2ds\right] -
 \mathbb E\left[\int_{t_0}^{T}(a(X_s) + b(Y_s))^2ds\right]
 \nonumber\\
 \label{criterion_justification}
 & = &
 \int_{\mathbb R^2}(\tau(x) +\nu(y) - (a(x) + b(y)))^2f(x,y)dxdy -
 \int_{\mathbb R^2}(a(x) + b(y))^2f(x,y)dxdy,
\end{eqnarray}
which is minimal for $(\tau,\nu) = (a,b)$, and then justifies our estimation procedure of $(a,b)$. Note that in Equality (\ref{criterion_justification}), the theoretical counterpart to the empirical norm appears and leads us to set
\begin{eqnarray*}
 \|(\tau,\nu)\|_{f}^{2} & := &
 \int_{\mathbb R^2}(\tau(x) +\nu(y))^2f(x,y)dxdy\\
 & = &
 \mathbb E\left(\int_{t_0}^{T}(\tau(X_s) +\nu(Y_s))^2ds\right) =\mathbb E(\|(\tau,\nu)\|_{N}^{2}).
\end{eqnarray*}
Consider $\Psi_{\bf m} :=\mathbb E(\widehat\Psi_{\bf m})$. As for the empirical norm, for every $\mathbf x = (x_1,\dots,x_{m_1 + m_2})\in\mathbb R^{m_1 + m_2}$,
\begin{eqnarray*}
 \mathbf x^*\Psi_{\bf m}\mathbf x & = & \|(\tau,\nu)\|_{f}^{2}
 \quad {\rm with}\quad
 \tau =\sum_{j = 1}^{m_1}x_j\varphi_j
 \quad {\rm and}\quad
 \nu =\sum_{k = 1}^{m_2}x_{m_1 + k}\psi_k.
\end{eqnarray*}
%


%
\subsection{Risk bound for the fixed model $\mathcal S_{\bf m}$}\label{section_non_adaptive_risk_bound}
Let us define
\begin{displaymath}
\mathfrak L_{\varphi}(m_1) :=\sup_{x\in A_1}\sum_{j = 1}^{m_1}\varphi_j(x)^2
\quad {\rm and}\quad
\mathfrak L_{\psi}(m_2) :=\sup_{y\in A_2}\sum_{k = 1}^{m_2}\psi_k(y)^2.
\end{displaymath}
First, in the sequel, $\mathbf m = (m_1,m_2)$ needs to fulfill the following stability condition.
%


%
\begin{assumption}\label{stability_condition}
There exists $r > 0$ such that
\begin{displaymath}
(\mathfrak L_{\varphi}(m_1) +\mathfrak L_{\psi}(m_2))
(\|\Psi_{\bf m}^{-1}\|_{\rm op}\vee 1)\leqslant
\frac{\mathfrak c_r}{2}\cdot\frac{N}{\log(N)}
\quad {\rm with}\quad
\mathfrak c_r =\frac{1 -\log(2)}{1 + r}.
\end{displaymath}
\end{assumption}
\noindent
Assumption \ref{stability_condition} ensures that $\mathfrak L_{\varphi}(m_1)$, $\mathfrak L_{\psi}(m_2)$ and $\|\Psi_{\bf m}^{-1}\|_{\rm op}$ are finite, and thanks to the matrix Chernov's inequality established in Tropp \cite{TROPP12}, Theorem 1.1, one can establish the following key lemma in the spirit of Cohen et al. \cite{CDL13}.
%


%
\begin{lemma}\label{bound_Omega}
Under Assumptions \ref{assumption_drift_volatility}, \ref{assumption_XY} and \ref{stability_condition}, there exists a constant $\mathfrak c_{\ref{bound_Omega}} > 0$, not depending on $\mathbf m$ and $N$, such that
\begin{displaymath}
\mathbb P(\Omega_{\bf m}^{c})
\leqslant
\frac{\mathfrak c_{\ref{bound_Omega}}}{N^r},
\quad\textrm{where}\quad
\Omega_{\bf m} :=
\left\{\|\Psi_{\bf m}^{-\frac{1}{2}}\widehat\Psi_{\bf m}\Psi_{\bf m}^{-\frac{1}{2}} -
I_{m_1 + m_2}\|_{\rm op}\leqslant\frac{1}{2}\right\},
\end{displaymath}
and $\Psi_{\bf m}^{-1/2}$ is the square root of the positive definite symmetric matrix $\Psi_{\bf m}^{-1}$.
\end{lemma}
\noindent
Lemma \ref{bound_Omega} says that there exists an event, with probability near of one, on which the empirical and theoretical norms are equivalent on $\mathcal S_{\bf m}$. Precisely, on $\Omega_{\bf m}$, for every $(\tau,\nu)\in\mathcal S_{\bf m}$,
\begin{displaymath}
\frac{1}{2}\|(\tau,\nu)\|_{f}^{2}
\leqslant
\|(\tau,\nu)\|_{N}^{2}
\leqslant
\frac{3}{2}\|(\tau,\nu)\|_{f}^{2}.
\end{displaymath}
Now, let us consider the empirical counterpart of Assumption \ref{stability_condition}, that is the event
\begin{displaymath}
\Lambda_{\bf m} :=
\left\{
(\mathfrak L_{\varphi}(m_1) +\mathfrak L_{\psi}(m_2))
(\|\widehat\Psi_{\bf m}^{-1}\|_{\rm op}\vee 1)\leqslant
\mathfrak c_r\frac{N}{\log(N)}
\right\}.
\end{displaymath}
Note that on the event $\Lambda_{\bf m}$, the lowest eigenvalue of $\widehat\Psi_{\bf m}$ is positive, and then $\widehat\Psi_{\bf m}$ is invertible. The following proposition provides a risk bound - with respect to the empirical norm - on the truncated estimator
\begin{displaymath}
(\widetilde a_{m_1},\widetilde b_{m_2}) :=
(\widehat a_{m_1},\widehat b_{m_2})\mathbf 1_{\Lambda_{\bf m}}.
\end{displaymath}
%


%
\begin{proposition}\label{nonadaptive_risk_bound}
Consider $A = A_1\times A_2$. Under Assumptions \ref{assumption_drift_volatility}, \ref{assumption_XY} and \ref{stability_condition} with $r\geqslant 5$, if $a + b\in\mathbb L^4(A,f(x,y)dxdy)$, then there exists a constant $\mathfrak c_{\ref{nonadaptive_risk_bound}} > 0$, not depending on $\mathbf m$ and $N$, such that
\begin{equation}\label{nonadaptive_risk_bound_1}
\mathbb E(\|(\widetilde a_{m_1},\widetilde b_{m_2}) - (a,b)\mathbf 1_A\|_{N}^{2})
\leqslant
\min_{(\tau,\nu)\in\mathcal S_{\bf m}}\|(\tau,\nu) - (a,b)\mathbf 1_A\|_{f}^{2} +
\frac{2\|\sigma\|_{\infty}^{2}}{T_0}\cdot\frac{m_1 + m_2}{N} +
\frac{\mathfrak c_{\ref{nonadaptive_risk_bound}}}{N}.
\end{equation}
\end{proposition}
\noindent
Moreover, one can control the $f$-weighted risk of $(\widehat a_{m_1},\widehat b_{m_2})$.
%


%
\begin{proposition}\label{nonadaptive_risk_bound_f_weighted}
Under Assumptions \ref{assumption_drift_volatility}, \ref{assumption_XY} and \ref{stability_condition} with $r\geqslant 7$, if $a + b\in\mathbb L^4(A,f(x,y)dxdy)$, then there exists a constant $\mathfrak c_{\ref{nonadaptive_risk_bound_f_weighted}} > 0$, not depending on $\mathbf m$ and $N$, such that
\begin{equation}\label{nonadaptive_risk_bound_f_weighted_1}
\mathbb E(\|(\widetilde a_{m_1},\widetilde b_{m_2}) - (a,b)\mathbf 1_A\|_{f}^{2})
\leqslant
13\min_{(\tau,\nu)\in\mathcal S_{\bf m}}\|(\tau,\nu) - (a,b)\mathbf 1_A\|_{f}^{2} +
\frac{8\|\sigma\|_{\infty}^{2}}{T_0}\cdot\frac{m_1 + m_2}{N} +
\frac{\mathfrak c_{\ref{nonadaptive_risk_bound_f_weighted}}}{N}.
\end{equation}
\end{proposition}
\noindent
First, at least when $(a,b)\mathbf 1_A$ is bounded, or by Theorem 1.2 in Menozzi et al. \cite{MPZ21} when $(a,b)\mathbf 1_A$ is not bounded but $Y$ is the solution of Equation (\ref{Y_diffusion}), $a + b\in\mathbb L^4(A,f(x,y)dxdy)$.
\\
\\
Now, let us say few words about the three terms in the right-hand sides of Inequalities (\ref{nonadaptive_risk_bound_1}) and (\ref{nonadaptive_risk_bound_f_weighted_1}):
\begin{enumerate}
 \item The first term is the squared bias of our estimator of $(a,b)$. By Proposition \ref{integrability_Pi} when $(a,b)$ is bounded, or by Theorem 1.2 in Menozzi et al. \cite{MPZ21} when $Y$ is the solution of Equation (\ref{Y_diffusion}),
 \begin{displaymath}
 \mathfrak c_{f,1} :=\sup_{x\in A_1}\int_{A_2}f(x,y)dy <\infty
 \quad {\rm and}\quad
 \mathfrak c_{f,2} :=\sup_{y\in A_2}\int_{A_1}f(x,y)dx <\infty,
 \end{displaymath}
 leading to
 \begin{equation}\label{bias_term_bound}
 \min_{(\tau,\nu)\in\mathcal S_{\bf m}}\|(\tau,\nu) - (a,b)\mathbf 1_A\|_{f}^{2}
 \leqslant
 2\left(\mathfrak c_{f,1}\|a_{m_1} - a\mathbf 1_{A_1}\|^2 +
 \mathfrak c_{f,2}\min_{\nu\in\mathbb S_{m_2}}\|\nu - b\mathbf 1_{A_2}\|^2\right),
 \end{equation}
 where $a_{m_1}$ is the orthogonal projection of $a$ on $S_{m_1}$ for the usual inner product in $\mathbb L^2(A_1)$. Inequality (\ref{bias_term_bound}) is refined in Section \ref{section_bias_term_degenerate_constraint} (resp. Section \ref{section_bias_term_non_degenerate_constraint}) when $\delta_{m_2} = 0$ (resp. $\delta_{m_2}\neq 0$).
 \item The second term is a control of order $(m_1 + m_2)/N$, which is standard in the nonparametric regression framework, of the variance of our estimator of $(a,b)$.
 \item The last term is a negligible remainder of order $1/N$.
\end{enumerate}
%


%
\subsection{Refined bound on the bias when $\delta_{m_2} = 0$ and estimation rate}\label{section_bias_term_degenerate_constraint}
When $\delta_{m_2} = 0$, $\mathbb S_{m_2} =\Sigma_{m_2}$, and then Inequality (\ref{bias_term_bound}) is equivalent to
\begin{equation}\label{bias_term_bound_delta_0}
\min_{(\tau,\nu)\in\mathcal S_{\bf m}}\|(\tau,\nu) - (a,b)\mathbf 1_A\|_{f}^{2}
\leqslant
2(\mathfrak c_{f,1}\|a_{m_1} - a\mathbf 1_{A_1}\|^2 +
\mathfrak c_{f,2}\|b_{m_2} - b\mathbf 1_{A_2}\|^2),
\end{equation}
where $b_{m_2}$ is the orthogonal projection of $b$ on $\Sigma_{m_2}$ for the usual inner product in $\mathbb L^2(A_2)$.
%


%
\begin{example}\label{rate_trigonometric_basis}
(Trigonometric basis) First, assume that $A_1 = [0,1]$, and that $(\varphi_1,\dots,\varphi_{m_1})$ is the $[0,1]$-supported $m_1$-dimensional trigonometric basis: for every $x\in [0,1]$ and $j\in\mathbb N^*$ satisfying $2j + 1\leqslant m_1$,
\begin{displaymath}
\varphi_1(x) = 1,\quad
\varphi_{2j + 1}(x) =\sqrt 2\sin(2\pi jx)\quad\textrm{and}\quad
\varphi_{2j}(x) =\sqrt 2\cos(2\pi jx).
\end{displaymath}
Assume also that $A_2 = [0,1]$, and that $(\psi_1,\dots,\psi_{m_2})$ is the $[0,1]$-supported $m_2$-dimensional trigonometric basis with no constant function: for every $x\in [0,1]$ and $j\in\mathbb N^*$ satisfying $2j\leqslant m_2$,
\begin{displaymath}
\psi_{2j - 1}(x) =\sqrt 2\cos(2\pi jx)\quad\textrm{and}\quad
\psi_{2j}(x) =\sqrt 2\sin(2\pi jx).
\end{displaymath}
The function $\psi_0\equiv 1$ may be omitted, because in our setting
\begin{displaymath}
\int_{A_2}b(y)dy = 0,
\quad\textrm{leading to}\quad
\langle b,\psi_0\rangle = 0.
\end{displaymath}
Clearly, for this version of the trigonometric basis, $\delta_{m_2} = 0$. Now, let us define the Fourier-Sobolev spaces:
\begin{displaymath}
\mathbb W_{2}^{\gamma}([0,1]) :=
\left\{\varphi : [0,1]\rightarrow\mathbb R
\textrm{ $\gamma$ times differentiable} :
\int_{0}^{1}\varphi^{(\gamma)}(x)^2dx <\infty\right\}
\textrm{ $;$ }\gamma > 0.
\end{displaymath}
Consider $\alpha,\beta > 0$, and assume that $a$ (resp. $b$) belongs to $\mathbb W_{2}^{\alpha}([0,1])$ (resp. $\mathbb W_{2}^{\beta}([0,1])$). So, by DeVore and Lorentz \cite{DL93}, Corollary 2.4 p. 205, there exist two constants $\mathfrak c_{\alpha},\mathfrak c_{\beta} > 0$, not depending on $m_1$ and $m_2$ respectively, such that
\begin{displaymath}
\|a_{m_1} - a\mathbf 1_{A_1}\|^2
\leqslant\mathfrak c_{\alpha}m_{1}^{-2\alpha}
\quad\textrm{and}\quad
\|b_{m_2} - b\mathbf 1_{A_2}\|^2
\leqslant\mathfrak c_{\beta}m_{2}^{-2\beta},
\end{displaymath}
leading to
\begin{displaymath}
\min_{(\tau,\nu)\in\mathcal S_{\bf m}}\|(\tau,\nu) - (a,b)\mathbf 1_A\|_{f}^{2}
\lesssim
m_{1}^{-2\alpha} + m_{2}^{-2\beta}.
\end{displaymath}
Therefore, under the assumptions of Proposition \ref{nonadaptive_risk_bound},
\begin{displaymath}
\mathbb E(\|(\widetilde a_{m_1},\widetilde b_{m_2}) - (a,b)\mathbf 1_A\|_{N}^{2})
\lesssim
m_{1}^{-2\alpha} + m_{2}^{-2\beta} +
\frac{m_1 + m_2}{N}.
\end{displaymath}
By choosing $m_{1}^{\star}\asymp N^{1/(2\alpha + 1)}$ and $m_{2}^{\star}\asymp N^{1/(2\beta + 1)}$,
\begin{displaymath}
\mathbb E(\|(\widetilde a_{m_{1}^{\star}},\widetilde b_{m_{2}^{\star}}) - (a,b)\mathbf 1_A\|_{N}^{2})
\lesssim
N^{-\frac{2\alpha}{2\alpha + 1}} + N^{-\frac{2\beta}{2\beta + 1}},
\end{displaymath}
and then the rate coincides with the worst of the two $1$-dimensional optimal rates for the projection least squares estimation of the regression function in compact setting. However, the rate doesn't suffer from the curse of dimensionality. Finally, if there exists $f_0 > 0$ such that
\begin{equation}\label{lower_bound_time_average_density}
f(x,y)\geqslant f_0\textrm{ $;$ }\forall (x,y)\in A,
\end{equation}
then for every $(\tau,\nu)\in\mathcal S_{\bf m}$,
\begin{eqnarray*}
 \|(\tau,\nu)\|_{f}^{2}
 & \geqslant &
 f_0\int_A(\tau(x) +\nu(y))^2dxdy\\
 & = &
 f_0\cdot (\|\tau\mathbf 1_{A_1}\|^2 +\|\nu\mathbf 1_{A_2}\|^2)
 \quad\textrm{because}\quad
 \delta_{m_2} = 0.
\end{eqnarray*}
So,
\begin{displaymath}
\|\Psi_{\bf m}^{-1}\|_{\rm op}\leqslant\frac{1}{f_0}.
\end{displaymath}
In that case, Assumption \ref{stability_condition} says that $m_1$ and $m_2$ need to be of order $N/\log(N)$, which is a mild condition making the optimal choices $m_{1}^{\star}$ and $m_{2}^{\star}$ possible.
\end{example}
%


%
\subsection{Refined bound on the bias when $\delta_{m_2}\neq 0$ and estimation rate}\label{section_bias_term_non_degenerate_constraint}
First, the following Lemma provides a refinement of Inequality (\ref{bias_term_bound}) when $\delta_{m_2}\neq 0$.
%


%
\begin{lemma}\label{control_bias_term_non_degenerate_constraint}
Assume that $\delta_{m_2}\neq 0$. Then,
\begin{eqnarray}
 \min_{(\tau,\nu)\in\mathcal S_{\bf m}}\|(\tau,\nu) - (a,b)\mathbf 1_A\|_{f}^{2}
 & \leqslant &
 2\mathfrak c_{f,1}\|a_{m_1} - a\mathbf 1_{A_1}\|^2
 \nonumber\\
 \label{control_bias_term_non_degenerate_constraint_1}
 & &
 \hspace{0.5cm} +
 2\mathfrak c_{f,2}\left[\|b_{m_2} - b\mathbf 1_{A_2}\|^2 +
\frac{1}{\|\delta_{m_2}\|_{2,\mathbb R^{m_2}}^{2}}\left(
\int_{A_2}(b_{m_2}(y) - b(y))dy\right)^2\right].
\end{eqnarray}
Moreover, if $\lambda(A_2) <\infty$, then
\begin{displaymath}
\min_{(\tau,\nu)\in\mathcal S_{\bf m}}\|(\tau,\nu) - (a,b)\mathbf 1_A\|_{f}^{2}
\leqslant
2\mathfrak c_{f,1}\|a_{m_1} - a\mathbf 1_{A_1}\|^2 +
2\mathfrak c_{f,2}\left(1 +
\frac{\lambda(A_2)}{\|\delta_{m_2}\|_{2,\mathbb R^{m_2}}^{2}}\right)
\|b_{m_2} - b\mathbf 1_{A_2}\|^2.
\end{displaymath}
\end{lemma}
\noindent
Now, consider an interval $I\subset\mathbb R$ and an (arbitrary) orthonormal family $(\theta_k)_{k\in\mathbb N}$ of $\mathbb L^2(I)$. Moreover, let us define the associated general Sobolev spaces:
\begin{displaymath}
\mathbb W_{\theta}^{\gamma}(I,L) :=
\left\{s\in\mathbb L^2(I) :\sum_{k = 0}^{\infty}k^{\gamma}\langle s,\theta_k\rangle^2\leqslant L\right\}
\textrm{ $;$ }
\gamma,L > 0.
\end{displaymath}
For any $m\in\mathbb N$ and $s\in\mathbb L^2(I)$, let $s_m$ be the orthogonal projection of $s$ on ${\rm span}(\theta_0,\dots,\theta_{m - 1})$, and note that if $s\in\mathbb W_{\theta}^{\gamma}(I,L)$, then
\begin{displaymath}
\|s - s_m\|^2 =
\sum_{k = m}^{\infty}\langle s,\theta_k\rangle^2 =
\sum_{k = m}^{\infty}\langle s,\theta_k\rangle^2k^{\gamma}k^{-\gamma}\leqslant Lm^{-\gamma}.
\end{displaymath}
Thanks to Proposition \ref{nonadaptive_risk_bound} together with Lemma \ref{control_bias_term_non_degenerate_constraint}, the following proposition provides a rate for our projection least squares estimator of $(a,b)$ when $\delta_{m_2}$ is lower bounded and $a$ (resp. $b$) belongs to a general Sobolev space associated to the $\varphi_j$'s (resp. the $\psi_k$'s).
%


%
\begin{proposition}\label{rate_Sobolev}
Assume that there exist $\omega\geqslant -1$, and two positive constants $\mathfrak c_{\psi,1}$ and $\mathfrak c_{\psi,2}$ such that, for every $m_2\in\mathbb N^*$,
\begin{equation}\label{rate_Sobolev_1}
\|\delta_{m_2}\|_{2,\mathbb R^{m_2}}^{2}\geqslant\mathfrak c_{\psi,1}m_{2}^{\omega + 1}
\end{equation}
and
\begin{equation}\label{rate_Sobolev_2}
\left|\int_{A_2}\psi_k(y)dy\right|^2\leqslant\mathfrak c_{\psi,2}k^{\omega}
\textrm{ $;$ }\forall k\in\{1,\dots,m_2\}.
\end{equation}
Moreover, consider $\alpha,L_1,L_2 > 0$ and $\beta >\omega + 1$. Under the assumptions of Proposition \ref{nonadaptive_risk_bound}, if $a\in\mathbb W_{\varphi}^{\alpha}(A_1,L_1)$ and $b\in\mathbb W_{\psi}^{\beta}(A_2,L_2)$ then, for $m_{1}^{\star}\asymp N^{1/(\alpha + 1)}$ and $m_{2}^{\star}\asymp N^{1/(\beta + 1)}$,
\begin{displaymath}
\mathbb E(\|(\widetilde a_{m_{1}^{\star}},\widetilde b_{m_{2}^{\star}}) - (a,b)\mathbf 1_A\|_{N}^{2})
\lesssim
N^{-\frac{\alpha}{\alpha + 1}} + N^{-\frac{\beta}{\beta + 1}}.
\end{displaymath}
\end{proposition}
%


%
\begin{example}\label{rate_Laguerre_basis}
(Laguerre basis) Consider $I =\mathbb R_+$, and let $(\ell_k)_{k\in\mathbb N}$ be the Laguerre basis: for every $x\in I$, $\ell_0(x) =\sqrt 2e^{-x}\mathbf 1_I(x)$, and for every $k\in\mathbb N^*$,
\begin{displaymath}
\ell_k(x) =\sqrt 2L_k(2x)e^{-x}\mathbf 1_I(x)
\quad\textrm{with}\quad
L_k(x) =\sum_{j = 0}^{k}(-1)^j
\begin{pmatrix}
 k\\
 j
\end{pmatrix}\frac{x^j}{j!}.
\end{displaymath}
The $L_k$'s are the Laguerre polynomials, and $(\ell_k)_{k\in\mathbb N}$ is a Hilbert basis of $\mathbb L^2(I)$ for the usual inner product $\langle .,.\rangle$. Moreover, $\ell_0(0) =\sqrt 2$ and, for every $k\in\mathbb N^*$,
\begin{displaymath}
\|\ell_k\|_{\infty}\leqslant\sqrt 2
\quad\textrm{and}\quad
\int_{0}^{\infty}\ell_k(y)dy =\sqrt 2(-1)^k,
\end{displaymath}
leading to $\mathfrak L_{\ell}(m) = 2m$ ($m\in\mathbb N$) and to the conditions (\ref{rate_Sobolev_1}) and (\ref{rate_Sobolev_2}) with $\omega = 0$. Therefore, if $\varphi_k =\psi_k =\ell_{k - 1}$ for every $k\in\mathbb N^*$, and if $a\in\mathbb W_{\varphi}^{\alpha}(I,L_1)$ and $b\in\mathbb W_{\psi}^{\beta}(I,L_2)$ ($\alpha,L_1,L_2 > 0$ and $\beta > 1$) then, by Proposition \ref{rate_Sobolev},
\begin{displaymath}
\mathbb E(\|(\widetilde a_{m_{1}^{\star}},\widetilde b_{m_{2}^{\star}}) - (a,b)\mathbf 1_A\|_{N}^{2})
\lesssim
N^{-\frac{\alpha}{\alpha + 1}} + N^{-\frac{\beta}{\beta + 1}}.
\end{displaymath}
Finally, about the relationship between the natural Laguerre-Sobolev spaces introduced in Bongioanni and Torrea \cite{BT09} and the coefficients-based Sobolev spaces $\mathbb W_{\ell}^{\gamma}(I,L)$ ($\gamma,L > 0$), and for additional details as regularity properties, the reader can refer to Comte and Genon-Catalot \cite{CGC15} and Belomestny et al. \cite{BCG16}.
\end{example}
%


%
\begin{example}\label{rate_Hermite_basis}
(Hermite basis) Consider $I =\mathbb R$, and let $(h_k)_{k\in\mathbb N}$ be the Hermite basis: for every $x\in I$ and $k\in\mathbb N$,
\begin{displaymath}
h_k(x) = (2^kk!\sqrt\pi)^{-\frac{1}{2}}H_k(x)e^{-\frac{x^2}{2}}\mathbf 1_I(x)
\quad\textrm{with}\quad
H_k(x) = (-1)^ke^{x^2}\frac{d^k}{dx^k}e^{-x^2}.
\end{displaymath}
The $H_k$'s are the Hermite polynomials, $(h_k)_{k\in\mathbb N}$ is a Hilbert basis of $\mathbb L^2(I)$ for the usual inner product $\langle .,.\rangle$, and $\mathfrak L_{\ell}(m)\lesssim\sqrt m$ ($m\in\mathbb N$) by Lemma 1 in Comte and Lacour \cite{CL23}. Moreover, for any $k\in\mathbb N$,
\begin{itemize}
 \item Since $h_{2k + 1}$ is an odd function, $\displaystyle{\int_{-\infty}^{\infty}h_{2k + 1}(y)dy = 0}$.
 \item Since $\sqrt 2h_j' =\sqrt jh_{j - 1} -\sqrt{j + 1}h_{j + 1}$ ($j\in\mathbb N^*$),
 \begin{displaymath}
 \int_{-\infty}^{\infty}h_{2j}(y)dy =
 \sqrt{\frac{2j - 1}{2j}}\int_{-\infty}^{\infty}h_{2j - 2}(y)dy,
 \end{displaymath}
 leading to
 \begin{displaymath}
 \int_{-\infty}^{\infty}h_{2k}(y)dy =
 \sqrt 2\pi^{\frac{1}{4}}\frac{\sqrt{(2k)!}}{2^kk!}
 \underset{k\rightarrow\infty}{\sim}
 \sqrt 2k^{-\frac{1}{4}}
 \quad\textrm{thanks to the Stirling formula.}
 \end{displaymath}
\end{itemize}
So, the Hermite basis satisfies the conditions (\ref{rate_Sobolev_1}) and (\ref{rate_Sobolev_2}) with $\omega = -1/2$. Therefore, if $\varphi_k =\psi_k = h_{k - 1}$ for every $k\in\mathbb N^*$, and if $a\in\mathbb W_{\varphi}^{\alpha}(I,L_1)$ and $b\in\mathbb W_{\psi}^{\beta}(I,L_2)$ ($\alpha,L_1,L_2 > 0$ and $\beta > 1/2$) then, by Proposition \ref{rate_Sobolev},
\begin{displaymath}
\mathbb E(\|(\widetilde a_{m_{1}^{\star}},\widetilde b_{m_{2}^{\star}}) - (a,b)\mathbf 1_A\|_{N}^{2})
\lesssim
N^{-\frac{\alpha}{\alpha + 1}} + N^{-\frac{\beta}{\beta + 1}}.
\end{displaymath}
Finally, about the relationship between the natural Hermite-Sobolev spaces introduced in Bongioanni and Torrea \cite{BT06} and the coefficients-based Sobolev spaces $\mathbb W_{h}^{\gamma}(I,L)$ ($\gamma,L > 0$), and for additional details as regularity properties, the reader can refer to Belomestny et al. \cite{BCG18}.
\end{example}
\noindent
Note that one can mix the bases choice. For instance, $(\varphi_j)_{j\in\mathbb N}$ could be the trigonometric basis, and $(\psi_k)_{k\in\mathbb N}$ the Laguerre basis.
%


%
\subsection{Model selection}\label{section_model_selection}
In this section, Assumption \ref{stability_condition} is set through the collection of models
\begin{displaymath}
\mathcal M_N :=
\left\{
\mathbf m = (m_1,m_2)\in\{1,\dots,N\}^2 :
(\mathfrak L_{\varphi}(m_1) +\mathfrak L_{\psi}(m_2))
(\|\Psi_{\bf m}^{-1}\|_{\rm op}\vee 1)\leqslant
\frac{\mathfrak c_r}{2}\cdot\frac{N}{\log(N)}\right\}.
\end{displaymath}
For instance, in Example \ref{rate_trigonometric_basis}, since $\mathfrak L_{\varphi}(m_1)\leqslant 2m_1$ and $\mathfrak L_{\psi}(m_2)\leqslant 2m_2$, and under the condition (\ref{lower_bound_time_average_density}) on $f$ which leads to $\|\Psi_{\bf m}^{-1}\|_{\rm op}\leqslant f_{0}^{-1}$, if $(m_1,m_2)\in\{1,\dots,N\}^2$ satisfies
\begin{displaymath}
m_1 + m_2\leqslant
\frac{\mathfrak c_r}{4(f_{0}^{-1}\vee 1)}\cdot\frac{N}{\log(N)},
\quad {\rm then}\quad
(m_1,m_2)\in\mathcal M_N.
\end{displaymath}
So, at least for compactly supported bases, $\mathcal M_N$ is a large collection of models.
\\
\\
Now, let $(a_{m_1}^{N},b_{m_2}^{N})$ be the minimizer of $(u,v)\mapsto \|(u,v) - (a,b)\mathbf 1_A\|_{N}^{2}$ over $\mathcal S_{\bf m}$, and note that:
\begin{itemize}
 \item[(A)] By using the definition (\ref{nonadaptive_risk_bound_4}) of $\langle .,.\rangle_N$,
 \begin{displaymath}
 \|(a_{m_1}^{N},b_{m_2}^{N}) - (a,b)\mathbf 1_A\|_{N}^{2} =
 \|(a_{m_1}^{N},b_{m_2}^{N})\|_{N}^{2} -\|(a,b)\mathbf 1_A\|_{N}^{2}.
 \end{displaymath}
 \item[(B)] By the definition (\ref{objective_function}) of $\gamma_N$, for every $\tau =\theta_1\varphi_1 +\dots +\theta_{m_1}\varphi_{m_1}\in S_{m_1}$ and $\nu =\theta_{m_1 + 1}\psi_1 +\dots +\theta_{m_1 + m_2}\psi_{m_2}\in\Sigma_{m_2}$,
 \begin{displaymath}
 \gamma_N(\tau,\nu) =
 -\|(\tau,\nu)\|_{N}^{2} + 2R_{\bf m}(\theta),
 \end{displaymath}
 where
 \begin{displaymath}
 R_{\bf m}(\theta) :=\|(\tau,\nu)\|_{N}^{2} -
 \langle\theta,\widehat{\bf Z}_{\bf m}\rangle_{2,\mathbb R^{m_1 + m_2}} =
 \theta^*(\widehat\Psi_{\bf m}\theta -\widehat{\bf Z}_{\bf m}).
 \end{displaymath}
 Then, by (\ref{coefficients_pLS_estimator}),
 \begin{eqnarray*}
  R_{\bf m}(\widehat\theta) & = &
  -\left(\widehat\Psi_{\bf m}^{-1}\widehat{\bf Z}_{\bf m} -
  \frac{\mathbf d_{\bf m}^{*}\widehat\Psi_{\bf m}^{-1}\widehat{\bf Z}_{\bf m}}{
  \mathbf d_{\bf m}^{*}\widehat\Psi_{\bf m}^{-1}\mathbf d_{\bf m}}\cdot
  \widehat\Psi_{\bf m}^{-1}\mathbf d_{\bf m}\right)^*
  \frac{\mathbf d_{\bf m}^{*}\widehat\Psi_{\bf m}^{-1}\widehat{\bf Z}_{\bf m}}{
  \mathbf d_{\bf m}^{*}\widehat\Psi_{\bf m}^{-1}\mathbf d_{\bf m}}\cdot\mathbf d_{\bf m}\\
  & = &
  -\frac{\mathbf d_{\bf m}^{*}\widehat\Psi_{\bf m}^{-1}\widehat{\bf Z}_{\bf m}}{
  \mathbf d_{\bf m}^{*}\widehat\Psi_{\bf m}^{-1}\mathbf d_{\bf m}}\cdot
  \widehat{\bf Z}_{\bf m}^{*}\widehat\Psi_{\bf m}^{-1}\mathbf d_{\bf m} +
  \frac{(\mathbf d_{\bf m}^{*}\widehat\Psi_{\bf m}^{-1}\widehat{\bf Z}_{\bf m})^2}{
  \mathbf d_{\bf m}^{*}\widehat\Psi_{\bf m}^{-1}\mathbf d_{\bf m}} = 0,
 \end{eqnarray*}
 leading to
 \begin{displaymath}
 \gamma_N(\widehat a_{m_1},\widehat b_{m_2}) =
 -\|(\widehat a_{m_1},\widehat b_{m_2})\|_{N}^{2}.
 \end{displaymath}
\end{itemize}
Thanks to (A) and (B), $\gamma_N(\widehat a_{m_1},\widehat b_{m_2})$ may be interpreted as an empirical version of the bias of the projection least squares estimator of $(a,b)$ up to the additive constant $\|(a,b)\|_{N}^{2}$. So, in order to approach the bias-variance compromise, it makes sense to consider the model
\begin{displaymath}
(\widehat m_1,\widehat m_2) :=
\underset{(m_1,m_2)\in\widehat{\mathcal M}_N}{\rm argmin}
\{\gamma_N(\widehat a_{m_1},\widehat b_{m_2}) + {\rm pen}(m_1,m_2)\}
\quad {\rm with}\quad
{\rm pen}(m_1,m_2) =\frac{\kappa\|\sigma\|_{\infty}^{2}}{T_0}\cdot\frac{m_1 + m_2}{N},
\end{displaymath}
where $\kappa > 0$ is a constant to calibrate, and
\begin{displaymath}
\widehat{\mathcal M}_N :=
\left\{\mathbf m = (m_1,m_2)\in\{1,\dots,N\}^2 :
(\mathfrak L_{\varphi}(m_1) +\mathfrak L_{\psi}(m_2))
(\|\widehat\Psi_{\bf m}^{-1}\|_{\rm op}\vee 1)\leqslant
\mathfrak c_r\frac{N}{\log(N)}\right\}.
\end{displaymath}
In the sequel, the $S_{m_1}$'s and the $\Sigma_{m_2}$'s are assumed to be nested. Then, for every $\mathbf m = (m_1,m_2)$ and $\mathbf m' = (m_1',m_2')$ belonging to $\{1,\dots,N\}^2$,
\begin{equation}\label{consequence_nested}
\mathcal S_{\bf m} +\mathcal S_{\mathbf m'}\subset\mathcal S_{(m_1\vee m_1',m_2\vee m_2')}.
\end{equation}
%


%
\begin{theorem}\label{adaptive_risk_bound}
Assume that $r\geqslant 7$ in both the definitions of $\mathcal M_N$ and $\widehat{\mathcal M}_N$. Under Assumptions \ref{assumption_drift_volatility} and \ref{assumption_XY}, if the $S_{m_1}$'s and the $\Sigma_{m_2}$'s are nested, then there exist two positive constants $\kappa_0$ and $\mathfrak c_{\ref{adaptive_risk_bound}}$, not depending on $N$, such that for every $\kappa\geqslant\kappa_0$,
\begin{displaymath}
\mathbb E(\|(\widehat a_{\widehat m_1},\widehat b_{\widehat m_2}) - (a,b)\mathbf 1_A\|_{N}^{2})
\leqslant
\mathfrak c_{\ref{adaptive_risk_bound}}\left(
\min_{\mathbf m\in\mathcal M_N}\left\{
\mathbb E(\|(\widehat a_{m_1},\widehat b_{m_2}) - (a,b)\mathbf 1_A\|_{N}^{2}\mathbf 1_{\Omega_N}) +
\kappa\frac{m_1 + m_2}{N}\right\} +\frac{1}{N}\right).
\end{displaymath}
\end{theorem}
\noindent
We underline that model selection is a finite-sample procedure and the result stated in Theorem \ref{adaptive_risk_bound} is nonasymptotic, contrary to rate results of sections \ref{section_bias_term_degenerate_constraint} and \ref{section_bias_term_non_degenerate_constraint}.  This result holds for the trigonometric basis as well as for Hermite's and Laguerre's bases, when $\delta_{m_2} = 0$ or $\delta_{m_2}\neq 0$. Theorem \ref{adaptive_risk_bound} shows that the adaptive projection least squares estimator of $(a,b)$ reaches a data-driven bias-variance compromise for a large enough constant $\kappa$. This one is calibrated once for all along preliminary simulation experiments. The reader can refer to Baudry et al. \cite{BMM12} to understand how works the calibration procedure.
%


%
\section{Numerical experiments}\label{section_numerical_experiments}
\begin{figure}[h!]
\includegraphics[width=12cm,height=6cm]{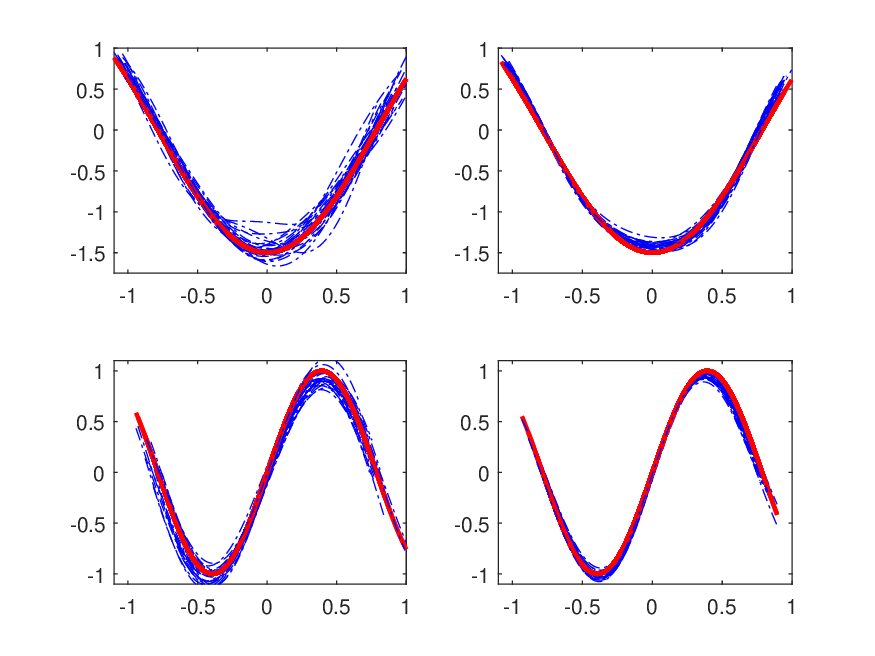}
\caption{\footnotesize The true functions in bold red, and 25 adaptive cpLS estimated (blue) $a_1$ on the first row and $b_1$ on the second one, for $N = 400$ (left) and $N = 1000$ (right), with $Y$ of type (A). 100*MSE: 2.86 (left) and 1.05 (right) for $a_1$, and 1.27 (left) and 0.54 (right) for $b_1$. Mean of selected dimensions: 16.7 (left) and 19.0 (right) for $a_1$, and 12.7 (left) and 20.0 (right) for $b_1$.}\label{fig1}
\end{figure}

\begin{figure}[h!]
\includegraphics[width=12cm,height=6cm]{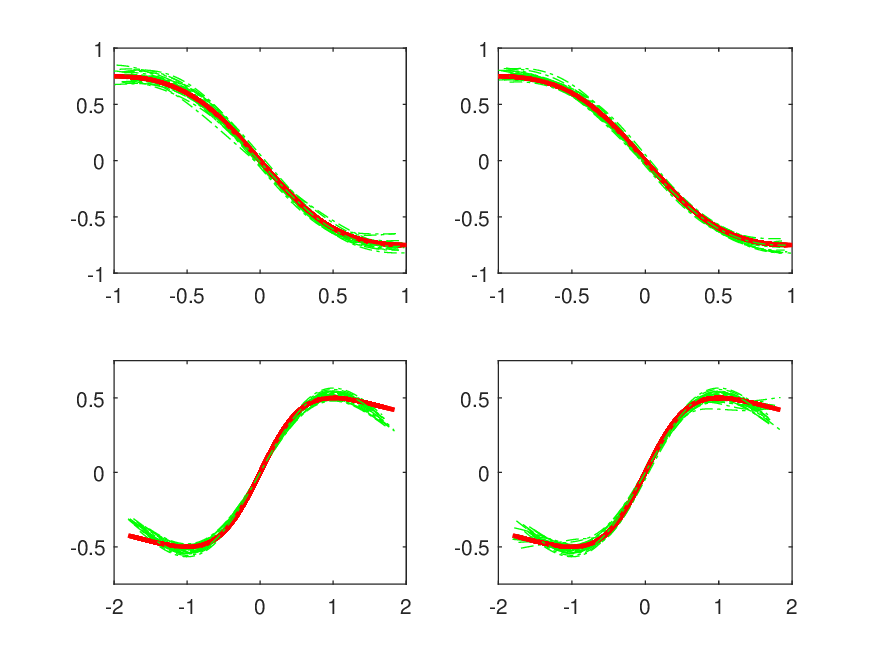}
\caption{\footnotesize The true functions in bold red. For $N = 1000$ and with $Y$ of type (A), in green and on the left (resp. right), 25 adaptive cpLS (resp. oracle-based) estimated $a_2$ on the first row and $b_2$ on the second one. 100*MSE: 0.23 (left) and 0.14 (right) for $a_2$, and 0.36 (left) and 0.29 (right) for $b_2$. Mean of selected dimensions: 7.28 (left) and 9.00 for $a_2$, and 2.00 (left) and 3.80 (right) for $b_2$.}
\label{fig2}
\end{figure}

\begin{figure}[h!]
\includegraphics[width=12cm,height=6cm]{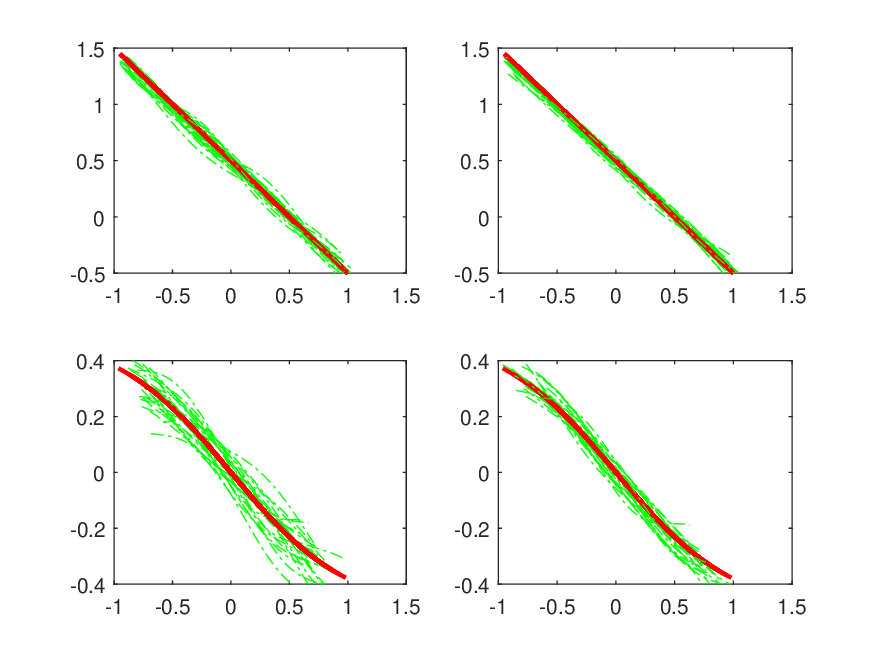}
\caption{\footnotesize The true functions in bold red. For $N = 1000$ and with $Y$ of type (B), in green and on the left (resp. right), 25 adaptive cpLS (resp. oracle-based) estimated $a_3$ on the first row and $b_3$ on the second one. 100*MSE: 0.52 (left) and 0.30 (right) for $a_3$, and 0.40 (left) and 0.14 (right) for $b_3$. Mean of selected dimensions: 12.4 (left) and 12.6 for $a_3$, and 10.6 (left) and 9.84 (right) for $b_3$.}
\label{fig3}
\end{figure}

Throughout this section, the paths of $X$ and $Y$ are generated from discrete-time approximations along the following dissection of $[0,T]$ ($T = 10$): $\{\ell\Delta\textrm{ $;$ }\ell = 0,\dots,n\}$ with $n = 500$ and $\Delta = 0.02$. Two values of $N$ are considered: $N = 400$ and $N = 1000$. 
\\
\\
{\bf SDE models:} The numerical experiments are carried out with the following explanatory processes:
\begin{itemize}
 \item[(A)] $Y =\sigma_YW_2(1 + W_{2}^{2})$, which is not defined by a stochastic differential equation, and
 \item[(B)] $Y =\sigma_YU$, where $U$ is the Ornstein-Uhlenbeck process defined by the following Langevin equation:
 \begin{displaymath}
 dU_t = -\frac{r}{2}U_tdt +\frac{\gamma}{2}dW_2(t)
 \quad {\rm with}\quad
 r = 2,\quad
 \gamma = 1\quad {\rm and}\quad
 U_0\rightsquigarrow\mathcal N\left(0,\frac{\gamma^2}{4r}\right).
 \end{displaymath}
\end{itemize}
For both of these processes, $\sigma_Y = 2$, and the first 20 observations are dropped out. Note also that the paths of the Ornstein-Uhlenbeck process $U$ are simulated thanks to an exact discretization scheme. Moreover, for the three following couples of functions $(a,b)$, the paths of the process $X$ are simulated thanks to the Euler scheme derived from Equation (\ref{main_equation}):
\begin{enumerate}
 \item $a_1(x) = -1.5\cos(2x)$ and $b_1(y) =\sin(4y)$,
 \item $a_2(x) = -1.5x/(1 + x^2)$ and $b_2(y) = y/(1 + y^2)$, and
 \item $a_3(x) = -x + 0.5$ and $b(y) = -0.5\tanh(y)$.
\end{enumerate}
The parameters involved in the definition of the $(a_{\ell},b_{\ell})$'s are carefully chosen to ensure that they take values of same order. Moreover, in all these models, the function $\sigma$ is assumed to be constantly equal to $1.5$.
\\
\\
{\bf Statistical implementation:} The constrained projection least squares (cpLS) estimator $(\widehat a_{m_1},\widehat b_{m_2})$ of $(a,b)$ is computed in the Hermite basis, which has been defined in Example \ref{rate_Hermite_basis}. Under the cutoff condition
\begin{displaymath}
(m_1 + m_2)\|\widehat\Psi_{\bf m}^{-1}\|_{\rm op}\leqslant
e^{14\log(10)}\cdot\frac{N}{\log(N)},
\end{displaymath}
the couple of dimensions $(\widehat m_1,\widehat m_2)$ is selected by minimizing the map
\begin{displaymath}
(m_1,m_2)\longmapsto
-\|(\widehat a_{m_1},\widehat b_{m_2})\|_{N}^{2} +
\kappa_H\frac{\sigma^2(m_1 + m_2)}{NT_0},
\end{displaymath}
where $\kappa_H = 8$, and $T_0 = T$ for the sake of simplicity. We assume that $\sigma$ is known, but it may be estimated. Moreover, let $(\texttt a_X,\texttt b_X)$ (resp. $(\texttt a_Y,\texttt b_Y)$) be the 98\% and 2\% (resp. 99\% and 1\%) quantiles of an arbitrary path of $X$ (resp. $Y$). The couple of dimensions $(m_{1}^{\star},m_{2}^{\star})$, minimizing the map
\begin{displaymath}
\Delta :
(m_1,m_2)\longmapsto
\int_{\texttt a_X}^{\texttt b_X}(\widehat a_{m_1}(x) - a(x))^2dx +
\int_{\texttt a_Y}^{\texttt b_Y}(\widehat b_{m_2}(y) - b(y))^2dy,
\end{displaymath}
is also computed and called "oracles" because the unknown true function $(a,b)$ is involved in the definition of $\Delta$. Note that the MSE of the cpLS estimator of $(a,b)$ is computed thanks to a formula depending on $(\texttt a_X,\texttt b_X)$ and $(\texttt a_Y,\texttt b_Y)$ in the same way as $\Delta$.
\\
\\
{\bf Results and comments:} First, all numerical results are gathered in Table \ref{table_MSE}. The MSE of the adaptive cpLS estimator decreases when $N$ increases, while the selected dimensions both increase. Similar results are observed for the two explanatory processes under consideration ((A) and (B)), except when $(a,b) = (a_1,b_1)$ for which the Ornstein-Uhlenbeck process leads to a significantly higher MSE for both $N = 400$ and $N = 1000$. The errors of the adaptive cpLS estimator are generally twice those of the oracle-based estimator, which seems satisfactory. However, note that several selected dimensions are too small compared to the oracles. This suggests that the penalty constant could be improved, but let us mention that this choice is rather "sensitive".
\begin{table}[h!]
\begin{tabular}{cc|cccc|cccc}
 & & \multicolumn{4}{c|}{$N = 400$} & \multicolumn{4}{c}{$N = 1000$}\\
 $(a,b)$ & $Y$ & MSE$_{\rm (std)}$ & MSE-O$_{\rm (std)}$ & Dim & Dim-O & MSE$_{\rm (std)}$ & MSE-O$_{\rm (std)}$ & Dim & Dim-O\\
 \hline
 & & & & & & & & & \\
 $a_1$ & (A) & 2.56$_{(1.88)}$ & 1.19$_{(0.33)}$ & 16.6 & 14.4 & 1.47$_{(0.88)}$ & 0.81$_{(1.59)}$ & 18.9 & 15.5\\
 $b_1$ & & 1.52$_{(1.18)}$ & 0.67$_{(0.54)}$ & 15.2 & 13.8 & 0.80$_{(0.53)}$ & 0.44$_{(0.34)}$ & 19.5 & 14.5\\
 $a_1$ & (B) & 20.9$_{(12.5)}$ & 3.17$_{1.73)}$ & 15.3 & 13.3 & 13.4$_{(7.97)}$ & 2.69$_{(1.23)}$ & 18.3 & 13.5\\
 $b_1$ & & 34.5$_{(18.7)}$ & 4.42$_{(2.73)}$ & 8.13 & 8.94 & 22.0$_{(12.7)}$ & 3.61$_{(2.14)}$ &8.2 & 9.34\\
 \hline
 $a_2$ & (A) & 0.60$_{(0.51)}$ & 0.32$_{(0.25)}$ & 5.16 & 6.66 & 0.27$_{(0.21)}$ & 0.15$_{(0.11)}$ & 7.06 & 8.46\\
 $b_2$ & & 0.23$_{(0.33)}$ & 0.17$_{(0.16)}$ & 2.03 & 3.64 & 0.15$_{(0.19)}$ & 0.08$_{(0.08)}$ & 2.78 & 4.22\\
 $a_2$ & (B) & 0.67$_{(0.58)}$ & 0.36$_{(0.27)}$ & 5.40 & 7.08 & 0.27$_{(1.98)}$ & 0.17$_{(0.12)}$ & 7.68 & 9.02\\
 $b_2$ & & 0.56$_{(0.43)}$ & 0.48$_{(0.34)}$ & 2.00 & 1.92 & 0.37$_{(0.15)}$ & 0.28$_{(0.13)}$ & 2.00 & 4.11\\
 \hline
 $a_3$ & (A) & 1.08$_{(0.74)}$ & 0.57$_{(0.43)}$ & 10.7 & 10.8 & 0.45$_{(0.33)}$ & 0.28$_{(0.21)}$ & 12.6 & 12.7\\
 $b_3$ & & 0.73$_{(0.60)}$ & 0.24$_{(0.23)}$ & 4.25 & 6.12 & 0.35$_{(0.26)}$ & 0.12$_{(0.08)}$ & 9.46 & 9.09\\
 $a_3$ & (B) & 0.86$_{(0.63)}$ & 0.54$_{(0.35)}$ & 10.3 & 11.7 & 0.35$_{(0.21)}$ & 0.23$_{(0.16)}$ & 12.2 & 13.2\\
 $b_3$ & & 2.53$_{(1.11)}$ & 0.96$_{0.72)}$ & 2.08 & 4.77 & 1.65$_{(0.97)}$ & 0.45$_{(0.34)}$ & 2.73 & 4.97\\
 & & & & & & & & & \\
\end{tabular}
\caption{\footnotesize Numerical results over 200 repetitions, for the three couples $(a_{\ell},b_{\ell})$ ($\ell = 1,2,3$), the two processes $Y$ defined by (A) and (B), and two sample sizes $N = 400$ and $N = 1000$. The columns MSE$_{\rm (std)}$ and MSE-O$_{\rm (std)}$ provide the 100*MSE with 100*standard deviation for the adaptive estimators and the oracles, while Dim and Dim-O provide the mean selected dimensions for both the estimators and oracles.}
\label{table_MSE}
\end{table}
\newline
Now, in order to illustrate what the cpLS estimator errors in Table \ref{table_MSE} concretely mean, we provide illustrations for each $(a_{\ell},b_{\ell})$ ($\ell = 1,2,3$). For $(a,b) = (a_1,b_1)$ and $Y$ of type (A), Figure \ref{fig1} allows to compare the adaptive cpLS estimator of $(a,b)$ when $N = 400$ and $N = 1000$. For $(a,b) = (a_2,b_2)$ and $Y$ of type (A), Figure \ref{fig2} presents 25 adaptive cpLS and oracle-based estimations for $N = 1000$. Lastly, Figure \ref{fig3} allows the same comparison as Figure \ref{fig2}, but with $(a,b) = (a_3,b_3)$ and $Y$ of type (B). Obviously, the oracle-based estimator performs better than the adaptive cpLS one. In particular, the beam of adaptive cpLS estimations of $b_3$ is significantly dispersed, which is not obvious from the MSE, probably due to the range of the function. However, globally, the adaptive cpLS estimator captures in a stable way the shape of the functions under consideration, and the method works very convincingly.
\appendix
%


%
\section{Proofs of probabilistic results (Section \ref{section_preliminaries})}\label{section_proofs_probabilistic_results}
%


%
\subsection{Proof of Proposition \ref{technical_result_empirical_norm}}\label{section_proof_technical_result_empirical_norm}
The proof of Proposition \ref{technical_result_empirical_norm}, but also that of Proposition \ref{general_sufficient_condition_AC}, rely on the following technical lemma.
%


%
\begin{lemma}\label{weak_differentiability_X}
Let $\xi$ be the solution of the stochastic differential equation
\begin{equation}\label{weak_differentiability_X_1}
\xi_t = x_0 +
\int_{0}^{t}\varphi(\xi_s,\zeta_s)ds +
\int_{0}^{t}\psi(\xi_s,\zeta_s)dW_1(s)
\textrm{ $;$ }t\in [0,T],
\end{equation}
where $\varphi$ and $\psi$ are continuously differentiable functions from $\mathbb R^2$ into $\mathbb R$, $\psi$ is bounded, both $\varphi$ and $\psi$ have bounded partial derivatives, and $\zeta$ is a $\mathbb F_2$-adapted process such that
\begin{equation}\label{weak_differentiability_X_2}
\mathbb E\left(\sup_{t\in [0,T]}\zeta_{t}^{2}\right) +
\sup_{s\in [0,T]}\mathbb E\left(\sup_{t\in [s,T]}
(\mathbf D_{s}^{2}\zeta_t)^2\right) <\infty.
\end{equation}
For any $t\in [0,T]$, $\xi_t\in\mathbb D^{1,2}$ and
\begin{displaymath}
\mathbf D_s\xi_t =
\begin{pmatrix}
 \alpha_1(s,t)\\
 \alpha_2(s,t)
\end{pmatrix}
e^{\beta(s,t)}\mathbf 1_{t\geqslant s}
\textrm{ $;$ }
\forall s\in [0,T]
\end{displaymath}
where, for every $s\in [0,t]$,
\begin{displaymath}
\beta(s,t) :=
\int_{s}^{t}\partial_1\psi(\xi_u,\zeta_u)dW_1(u) +
\int_{s}^{t}\left(\partial_1\varphi(\xi_u,\zeta_u) -
\frac{1}{2}\partial_1\psi(\xi_u,\zeta_u)^2\right)du
\end{displaymath}
and
\begin{eqnarray*}
 \alpha_{\ell}(s,t) & := &
 \psi(\xi_s,\zeta_s)\mathbf 1_{\ell = 1} +
 \int_{s}^{t}e^{-\beta(s,u)}\partial_2\varphi(\xi_u,\zeta_u)\mathbf D_{s}^{\ell}\zeta_udu\\
 & &
 \hspace{1cm} +
 \int_{s}^{t}e^{-\beta(s,u)}\partial_2\psi(\xi_u,\zeta_u)\mathbf D_{s}^{\ell}\zeta_u
 (dW_1(u) -\partial_1\psi(\xi_u,\zeta_u)du)
 \textrm{ $;$ }\ell = 1,2.
\end{eqnarray*}
\end{lemma}
\noindent
The proof of Lemma \ref{weak_differentiability_X} is postponed to Section \ref{section_proof_weak_differentiability_X}.
\\
\\
Let $X$ be the solution of Equation (\ref{main_equation}), which coincides with Equation (\ref{weak_differentiability_X_1}) by taking $\zeta := Y$ and, for every $(x,y)\in\mathbb R^2$, $\varphi(x,y) := a(x) + b(y)$ and $\psi(x,y) :=\sigma(x)$. Consider also $\tau\in\mathcal S_1$, $\nu\in\mathcal S_2$ fulfilling the condition (\ref{technical_result_empirical_norm_1}), and $t\in (0,T]$ satisfying
\begin{equation}\label{technical_result_empirical_norm_2}
\tau(X_t) +\nu(Y_t) = 0.
\end{equation}
First, $\mathbf D^1Y_t = 0$ because $Y$ is $\mathbb F_2$-adapted while $W_1$ and $W_2$ are independent. Then, by Equality (\ref{technical_result_empirical_norm_2}) and the chain rule for the Malliavin derivative (see Nualart \cite{NUALART06}, Proposition 1.2.3),
\begin{equation}\label{technical_result_empirical_norm_3}
\tau'(X_t)\mathbf D^1X_t = 0.
\end{equation}
Moreover, by Lemma \ref{weak_differentiability_X}, and since $\partial_2\psi = 0$,
\begin{equation}\label{technical_result_empirical_norm_4}
\mathbf D_{s}^{1}X_t =
\sigma(X_s)e^{\beta(s,t)}
\textrm{ $;$ }\forall s\in [0,t].
\end{equation}
Since $\inf_{\mathbb R}|\sigma| > 0$ (see Assumption \ref{assumption_drift_volatility}), $\mathbf D^1X_t\not\equiv 0$ by Equality (\ref{technical_result_empirical_norm_4}), and thus Equality (\ref{technical_result_empirical_norm_3}) leads to $\tau'(X_t) = 0$. Now, by Equality (\ref{technical_result_empirical_norm_2}), by the chain rule for the Malliavin derivative, and since $\tau'(X_t) = 0$,
\begin{displaymath}
\nu'(Y_t)\mathbf D^2Y_t = 0.
\end{displaymath}
Then, $\nu'(Y_t) = 0$ because $\mathbf D^2Y_t\not\equiv 0$ (see Assumption \ref{assumption_Y}.(1)). In particular, $\mathbb E(|\nu'(Y_t)|) = 0$, and since the distribution of $Y_t$ has a positive and continuously differentiable density with respect to the Lebesgue measure on $\mathbb R$ (see Assumption \ref{assumption_Y}.(2)), $\nu' = 0$. In conclusion, $\nu(.) = 0$ because $\nu$ fulfills (\ref{technical_result_empirical_norm_1}), and $\tau(X_t) = 0$ by Equality (\ref{technical_result_empirical_norm_2}).
%


%
\subsubsection{Proof of Lemma \ref{weak_differentiability_X}}\label{section_proof_weak_differentiability_X}
Let $(H_n)_{n\in\mathbb N}$ be the sequence defined by $H_0 = x_0$ and, for every $n\in\mathbb N$, $H_{n + 1} =\Phi(H_n)$, where
\begin{displaymath}
\Phi_t(H) :=
x_0 +
\int_{0}^{t}\varphi(H_s,\zeta_s)ds +
\int_{0}^{t}\psi(H_s,\zeta_s)dW_1(s)
\textrm{ $;$ }
t\in [0,T]\textrm{, }
H\in\mathcal C.
\end{displaymath}
Consider $n\in\mathbb N$ - for instance $n = 0$ - such that for any $\tau\in [0,T]$ and $r\in [0,\tau]$,
\begin{equation}\label{weak_differentiability_X_3}
\mu_n(r,\tau) :=
\mathbb E\left(
\sup_{t\in [r,\tau]}\|\mathbf D_rH_n(t)\|_{2,\mathbb R^2}^{2}\right)\leqslant
\mathfrak c_2\sum_{k = 0}^{n}\frac{\mathfrak c_{1}^{k}(\tau - r)^k}{k!},
\end{equation}
where
\begin{displaymath}
 \mathfrak c_1 :=
 3\max\{\|\psi\|_{\infty}^{2}\textrm{ $;$ }
 2(T\|\partial_1\varphi\|_{\infty}^{2} + 4\|\partial_1\psi\|_{\infty}^{2})\textrm{ $;$ }
 2(T\|\partial_2\varphi\|_{\infty}^{2} + 4\|\partial_2\psi\|_{\infty}^{2})\}
\end{displaymath}
and $\mathfrak c_2 :=\mathfrak c_{2,1} +\mathfrak c_{2,2}$ with
\begin{displaymath}
\mathfrak c_{2,\ell} =
\mathfrak c_1\left(1 + T\sup_{s\in [0,T]}\mathbb E\left(\sup_{t\in [s,T]}
(\mathbf D_{s}^{\ell}\zeta_t)^2\right)\right)
\textrm{ $;$ }\ell = 1,2.
\end{displaymath}
By Inequalities (\ref{weak_differentiability_X_2}) and (\ref{weak_differentiability_X_3}), $\zeta_u,H_n(u)\in\mathbb D^{1,2}$ for every $u\in [0,T]$, and by Nualart \cite{NUALART06}, Propositions 1.2.3 and 1.3.8, for every $t\in [0,T]$, $s\in [0,t]$ and $\ell = 1,2$,
\begin{eqnarray*}
 \mathbf D_{s}^{\ell}H_{n + 1}(t) & = &
 \int_{s}^{t}\mathbf D_{s}^{\ell}[\varphi(H_n(u),\zeta_u)]du +
 \psi(H_n(s),\zeta_s)\mathbf 1_{\ell = 1} +
 \int_{s}^{t}\mathbf D_{s}^{\ell}[\psi(H_n(u),\zeta_u)]dW_1(u)\\
 & = &
 \psi(H_n(s),\zeta_s)\mathbf 1_{\ell = 1} +
 \int_{s}^{t}[\partial_1\varphi(H_n(u),\zeta_u)\mathbf D_{s}^{\ell}H_n(u) +
 \partial_2\varphi(H_n(u),\zeta_u)\mathbf D_{s}^{\ell}\zeta_u]du\\
 & &
 \hspace{3cm} +
 \int_{s}^{t}[\partial_1\psi(H_n(u),\zeta_u)\mathbf D_{s}^{\ell}H_n(u) +
 \partial_2\psi(H_n(u),\zeta_u)\mathbf D_{s}^{\ell}\zeta_u]dW_1(u).
\end{eqnarray*}
Then, by the Doob inequality, by the isometry property of It\^o's integral, and since $\psi$ and the derivatives of both $\varphi$ and $\psi$ are bounded, for $\ell = 1,2$,
\begin{eqnarray*}
 \mathbb E\left(
 \sup_{t\in [r,\tau]}(\mathbf D_{r}^{\ell}H_{n + 1}(t))^2\right)
 & \leqslant &
 3\mathbb E(\psi(H_n(r),\zeta_r)^2)\mathbf 1_{\ell = 1}\\
 & &
 + 3(\tau - r)\int_{r}^{\tau}
 \mathbb E[(\partial_1\varphi(H_n(u),\zeta_u)\mathbf D_{r}^{\ell}H_n(u) +
 \partial_2\varphi(H_n(u),\zeta_u)\mathbf D_{r}^{\ell}\zeta_u)^2]du\\
 & &
 + 12\int_{r}^{\tau}\mathbb E[(\partial_1\psi(H_n(u),\zeta_u)\mathbf D_{r}^{\ell}H_n(u) +
 \partial_2\psi(H_n(u),\zeta_u)\mathbf D_{r}^{\ell}\zeta_u)^2]du\\
 & \leqslant &
 \mathfrak c_1\left(1 +
 \int_{r}^{\tau}\mathbb E((\mathbf D_{r}^{\ell}H_n(u))^2)du +
 \int_{r}^{\tau}\mathbb E((\mathbf D_{r}^{\ell}\zeta_u)^2)du\right)\\
 & \leqslant &
 \mathfrak c_{2,\ell} +
 \mathfrak c_1\int_{r}^{\tau}\mathbb E((\mathbf D_{r}^{\ell}H_n(u))^2)du,
\end{eqnarray*}
leading to
\begin{eqnarray*}
 \mu_{n + 1}(r,\tau)
 & \leqslant &
 \mathfrak c_2 +
 \mathfrak c_1\int_{r}^{\tau}
 \mathbb E\left(\sup_{v\in [r,u]}\|\mathbf D_rH_n(v)\|_{2,\mathbb R^2}^{2}\right)du\\
 & \leqslant &
 \mathfrak c_2 +
 \mathfrak c_1\int_{r}^{\tau}\left(
 \mathfrak c_2\sum_{k = 0}^{n}\frac{\mathfrak c_{1}^{k}(u - r)^k}{k!}\right)du
 \quad\textrm{by Inequality (\ref{weak_differentiability_X_3})}\\
 & = &
 \mathfrak c_2\left(1 +
 \sum_{k = 0}^{n}\frac{\mathfrak c_{1}^{k + 1}}{k!}\cdot
 \frac{(\tau - r)^{k + 1}}{k + 1}\right) =
 \mathfrak c_2\sum_{k = 0}^{n + 1}\frac{\mathfrak c_{1}^{k}(\tau - r)^k}{k!}.
\end{eqnarray*}
So, Inequality (\ref{weak_differentiability_X_3}) remains true for $\mu_{n + 1}$, and by induction
\begin{displaymath}
\sup_{n\in\mathbb N}\left\{
\sup_{s\in [0,T]}\mathbb E\left(
\sup_{t\in [s,T]}\|\mathbf D_sH_n(t)\|_{2,\mathbb R^2}^{2}\right)\right\}
\leqslant
\mathfrak c_2\sum_{k = 0}^{\infty}\frac{\mathfrak c_{1}^{k}T^k}{k!} =
\mathfrak c_2e^{\mathfrak c_1T} <\infty.
\end{displaymath}
Since, in addition, the Picard scheme $(H_n)_{n\in\mathbb N}$ converges to $\xi$ in $(\mathcal C,\|.\|_{\mathcal C})$ (see the proof of Nualart \cite{NUALART06}, Lemma 2.2.1), by Nualart \cite{NUALART06}, Lemma 1.2.3, $\xi_t\in\mathbb D^{1,2}$ for every $t\in [0,T]$. Therefore, by Nualart \cite{NUALART06}, Propositions 1.2.3 and 1.3.8, for every $t\in [0,T]$, $s\in [0,t]$ and $\ell = 1,2$,
\begin{eqnarray}
 \mathbf D_{s}^{\ell}\xi_t & = &
 \psi(\xi_s,\zeta_s)\mathbf 1_{\ell = 1} +
 \int_{s}^{t}[\partial_1\varphi(\xi_u,\zeta_u)\mathbf D_{s}^{\ell}\xi_u +
 \partial_2\varphi(\xi_u,\zeta_u)\mathbf D_{s}^{\ell}\zeta_u]du
 \nonumber\\
 & &
 \hspace{4cm} +
 \int_{s}^{t}[\partial_1\psi(\xi_u,\zeta_u)\mathbf D_{s}^{\ell}\xi_u +
 \partial_2\psi(\xi_u,\zeta_u)\mathbf D_{s}^{\ell}\zeta_u]dW_1(u)
 \nonumber\\
 \label{weak_differentiability_X_4}
 & = &
 \psi(\xi_s,\zeta_s)\mathbf 1_{\ell = 1} +
 \int_{s}^{t}\mathbf D_{s}^{\ell}\xi_udZ(s,u) + R_{\ell}(s,t),
\end{eqnarray}
where
\begin{eqnarray*}
 Z(s,t) & := &
 \int_{s}^{t}\partial_1\varphi(\xi_u,\zeta_u)du +
 \int_{s}^{t}\partial_1\psi(\xi_u,\zeta_u)dW_1(u)\\
 & &
 \hspace{1.5cm}{\rm and}\quad
 R_{\ell}(s,t) :=
 \int_{s}^{t}\mathbf D_{s}^{\ell}\zeta_u(
 \partial_2\varphi(\xi_u,\zeta_u)du +\partial_2\psi(\xi_u,\zeta_u)dW_1(u)).
\end{eqnarray*}
In conclusion, for any $s\in [0,T]$ and $\ell = 1,2$, since $(\mathbf D_{s}^{\ell}\xi_t)_{t\in [s,T]}$ is the solution of the linear stochastic differential equation (\ref{weak_differentiability_X_4}), $\mathbf D_{s}^{\ell}\xi_t =\alpha_{\ell}(s,t)e^{\beta(s,t)}$ for any $t\in [s,T]$, where
\begin{eqnarray*}
 \beta(s,t)
 & := &
 Z(s,t) -\frac{1}{2}\langle Z(s,.)\rangle_t\\
 & = &
 \int_{s}^{t}\partial_1\psi(\xi_u,\zeta_u)dW_1(u) +
 \int_{s}^{t}\left(\partial_1\varphi(\xi_u,\zeta_u) -
 \frac{1}{2}\partial_1\psi(\xi_u,\zeta_u)^2\right)du
\end{eqnarray*}
and
\begin{eqnarray*}
 \alpha_{\ell}(s,t)
 & := &
 \psi(\xi_s,\zeta_s)\mathbf 1_{\ell = 1} +
 \int_{s}^{t}e^{-\beta(s,u)}dR_{\ell}(s,u) -
 \int_{s}^{t}e^{-\beta(s,u)}d\langle Z(s,.),R_{\ell}(s,.)\rangle_u\\
 & = &
 \psi(\xi_s,\zeta_s)\mathbf 1_{\ell = 1} +
 \int_{s}^{t}e^{-\beta(s,u)}\partial_2\varphi(\xi_u,\zeta_u)\mathbf D_{s}^{\ell}\zeta_udu\\
 & &
 \hspace{3cm} +
 \int_{s}^{t}e^{-\beta(s,u)}\partial_2\psi(\xi_u,\zeta_u)\mathbf D_{s}^{\ell}\zeta_u
 (dW_1(u) -\partial_1\psi(\xi_u,\zeta_u)du).
\end{eqnarray*}
%


%
\subsection{Proof of Proposition \ref{general_sufficient_condition_AC}}\label{section_proof_general_sufficient_condition_AC}
Let $X$ be the solution of Equation (\ref{main_equation}), which coincides with Equation (\ref{weak_differentiability_X_1}) by taking $\zeta := Y$ and, for every $(x,y)\in\mathbb R^2$, $\varphi(x,y) := a(x) + b(y)$ and $\psi(x,y) :=\sigma(x)$. By Lemma \ref{weak_differentiability_X}, since $\partial_2\psi = 0$, and since $Y$ is $\mathbb F_2$-adapted while $W_1$ and $W_2$ are independent, for every $t\in [0,T]$ and $s\in [0,t]$,
\begin{equation}\label{general_sufficient_condition_AC_1}
\mathbf D_{s}^{1}X_t =
\sigma(X_s)e^{\beta(s,t)}
\quad {\rm and}\quad
\mathbf D_{s}^{2}X_t =
\int_{s}^{t}e^{\beta(s,u)}b'(Y_u)\mathbf D_sY_udu.
\end{equation}
The proof is dissected in two steps: Proposition \ref{general_sufficient_condition_AC}.(1) is established in Step 1, and Proposition \ref{general_sufficient_condition_AC}.(2) is established in Step 2.
\\
\\
{\bf Step 1.} For any $t\in (0,T]$, let $\Gamma_t$ be the Malliavin matrix of $(X_t,Y_t)$, which is defined by
\begin{displaymath}
\Gamma_t :=
\begin{pmatrix}
 \|\mathbf DX_t\|_{T}^{2} & \langle\mathbf DX_t,\mathbf DY_t\rangle_T\\
 \langle\mathbf DY_t,\mathbf DX_t\rangle_T & \|\mathbf DY_t\|_{T}^{2}
\end{pmatrix}.
\end{displaymath}
Since $\mathbf D^1Y_t = 0$, and by the Cauchy-Schwarz inequality,
\begin{eqnarray}
 \det(\Gamma_t) & = &
 \|\mathbf DX_t\|_{T}^{2}
 \|\mathbf DY_t\|_{T}^{2} -
 \langle\mathbf DX_t,\mathbf DY_t\rangle_{T}^{2}
 \nonumber\\
 & = &
 (\|\mathbf D^1X_t\|_{T}^{2} +\|\mathbf D^2X_t\|_{T}^{2})\|\mathbf D^2Y_t\|_{T}^{2} -
 \langle\mathbf D^2X_t,\mathbf D^2Y_t\rangle_{T}^{2}
 \nonumber\\
 \label{general_sufficient_condition_AC_2}
 & \geqslant &
 \|\mathbf D^1X_t\|_{T}^{2}\|\mathbf D^2Y_t\|_{T}^{2}.
\end{eqnarray}
Moreover, $\mathbf D^2Y_t\not\equiv 0$ (see Assumption \ref{assumption_Y}.(1)), and since $\inf_{\mathbb R}|\sigma| > 0$ (see Assumption \ref{assumption_drift_volatility}), (\ref{general_sufficient_condition_AC_1}) leads to $\mathbf D^1X_t\not\equiv 0$. Thus, $\det(\Gamma_t) > 0$, and by the Bouleau-Hirsch criterion (see Nualart \cite{NUALART06}, Theorem 2.1.2), the distribution of $(X_t,Y_t)$ has a density with respect to the Lebesgue measure on $\mathbb R^2$.
\\
\\
{\bf Step 2.} Assume that $\sigma$ is bounded, and that $a,b,\sigma\in C^{\infty}(\mathbb R)$ with all derivatives bounded. Assume also that $Y\in\mathbb H^{\infty}$, and that $1/\|\mathbf D^2Y_t\|_T$ belongs to $\mathbb L^p(\Omega)$ for every $p\geqslant 1$. First of all, note that since $Y\in\mathbb H^{\infty}$, by (\ref{general_sufficient_condition_AC_1}) together with Nualart \cite{NUALART06}, Lemma 2.2.1,
\begin{equation}\label{general_sufficient_condition_AC_3}
\mathbb E\left(\sup_{t\in [0,T]}|X_t|^p\right) +
\sup_{s\in [0,T]}\mathbb E\left(\sup_{t\in [s,T]}\|\mathbf D_sX_t\|_{p,\mathbb R^2}^{p}\right) <\infty
\textrm{ $;$ }\forall p\geqslant 1.
\end{equation}
{\it Notations:}
\begin{itemize}
 \item For every $k\in\mathbb N^*$, $\Pi_k$ is the set of all partitions of $\{1,\dots,k\}$, and $\Pi_{k}^{*} :=\Pi_k\backslash\{\{\{1,\dots,k\}\}\}$.
 \item For every $k\in\mathbb N^*$ and $j\in\{1,\dots,k\}$, $\Pi_{k\backslash j}$ is the set of all partitions of $\{1,\dots,k\}\backslash\{j\}$.
 \item Consider a set $E$, $k\in\mathbb N^*$ and $J\subset\{1,\dots,k\}$. For every $\mathbf x = (x_1,\dots,x_k)\in E^k$, $\mathbf x_J := (x_j)_{j\in J}$.
\end{itemize}
{\it Step 2.1.} Consider $t\in (0,T]$. By (\ref{general_sufficient_condition_AC_1}), since $\rho : x\in (0,\infty)\mapsto 1/x$ is a convex function, by the Jensen inequality, and since $\inf_{\mathbb R}|\sigma| > 0$,
\begin{eqnarray*}
 \frac{1}{\|\mathbf D^1X_t\|_{T}^{2}}
 & = &
 \frac{1}{t}\rho\left(\int_{0}^{t}\sigma(X_s)^2e^{2\beta(s,t)}\frac{ds}{t}\right)\\
 & \leqslant &
 \frac{1}{t^2}\int_{0}^{t}\rho(\sigma(X_s)^2e^{2\beta(s,t)})ds\\
 & \leqslant &
 \frac{\mathfrak c_1}{t^2}\int_{0}^{t}e^{-2\beta(s,t)}ds
 \quad {\rm with}\quad
 \mathfrak c_1 =\left(\inf_{x\in\mathbb R}|\sigma(x)|\right)^{-2}.
\end{eqnarray*}
Then, for any $\alpha > 1$,
\begin{eqnarray*}
 \mathbb E\left(\frac{1}{\|\mathbf D^1X_t\|_{T}^{4\alpha}}\right)
 & \leqslant &
 \frac{\mathfrak c_{1}^{2\alpha}}{t^{2\alpha + 1}}
 \int_{0}^{t}\mathbb E(e^{-4\alpha\beta(s,t)})ds\\
 & \leqslant &
 \frac{\mathfrak c_{1}^{2\alpha}}{t^{2\alpha + 1}}
 e^{2\alpha T\|a'\|_{\infty} + (8\alpha^2 + 2\alpha)T\|\sigma'\|_{\infty}^{2}}\\
 & &
 \hspace{0.75cm}\times
 \int_{0}^{t}\mathbb E\left[\exp\left(
 \int_{s}^{t}(-4\alpha)\sigma'(X_u)dW_1(u) -
 \frac{1}{2}\int_{s}^{t}(-4\alpha)^2\sigma'(X_u)^2du\right)\right]ds\\
 & \leqslant &
 \frac{\mathfrak c_2(\alpha)^2}{t^{2\alpha}}
 \quad {\rm with}\quad
 \mathfrak c_2(\alpha) =
 \mathfrak c_{1}^{\alpha}e^{6\alpha^2T\max\{\|a'\|_{\infty},\|\sigma'\|_{\infty}^{2}\}}.
\end{eqnarray*}
Therefore, by Inequality (\ref{general_sufficient_condition_AC_2}),
\begin{displaymath}
\mathbb E\left(\left|\frac{1}{\det(\Gamma_t)}\right|^{\alpha}\right)
\leqslant
\frac{\mathfrak c_2(\alpha)}{t^{\alpha}}
\mathbb E\left(\frac{1}{\|\mathbf DY_t\|_{T}^{4\alpha}}\right)^{\frac{1}{2}} <\infty.
\end{displaymath}
{\it Step 2.2.} The purpose of this step is to recursively establish that for every $k\in\mathbb N^*$,
\begin{itemize}
 \item $\mathbf P(k)$: For every $\mathbf l\in\{1,2\}^k$, $\mathbf s\in\mathbb [0,T]^k$ and $t\in [\max(\mathbf s),T]$,
 \begin{eqnarray*}
  \mathbf D_{\bf s}^{(k),\mathbf l}X_t
  & = &
  \alpha_{\bf l}(\mathbf s,t) +\beta_{\bf l}(\mathbf s,t) +
  \int_{\max(\mathbf s)}^{t}a'(X_u)
  \mathbf D_{\bf s}^{(k),\mathbf l}X_udu\\
  & &
  \hspace{4cm} +
  \int_{\max(\mathbf s)}^{t}\sigma'(X_u)
  \mathbf D_{\bf s}^{(k),\mathbf l}X_udW_1(u),
 \end{eqnarray*}
 where
 \begin{eqnarray*}
  \alpha_{\bf l}(\mathbf s,t) & := &
  \int_{\max(\mathbf s)}^{t}\sum_{\pi\in\Pi_k}
  b^{(|\pi|)}(Y_u)\prod_{J\in\pi}\mathbf D_{\mathbf s_J}^{(|J|),\mathbf l_J}Y_udu\\
  & &
  \hspace{1.5cm} +
  \mathbf 1_{k > 1}\int_{\max(\mathbf s)}^{t}\sum_{\pi\in\Pi_{k}^{*}}
  a^{(|\pi|)}(X_u)\prod_{J\in\pi}\mathbf D_{\mathbf s_J}^{(|J|),\mathbf l_J}X_udu\\
  & &
  \hspace{3cm} +
  \mathbf 1_{k > 1}\int_{\max(\mathbf s)}^{t}\sum_{\pi\in\Pi_{k}^{*}}
  \sigma^{(|\pi|)}(X_u)\prod_{J\in\pi}\mathbf D_{\mathbf s_J}^{(|J|),\mathbf l_J}X_udW_1(u)
 \end{eqnarray*}
 and
 \begin{eqnarray*}
  \beta_{\bf l}(\mathbf s,t)
  & := &
  \sigma(X_{s_1})\mathbf 1_{k = 1,\ell_1 = 1} +
  \mathbf 1_{k > 1}\left[\mathbf 1_{\ell_1 = 1}\sum_{\pi\in\Pi_{k\backslash 1}}\sigma^{(|\pi|)}(X_{s_1})
  \prod_{J\in\pi}\mathbf D_{\mathbf s_J}^{(|J|),\mathbf l_J}X_{s_1}\right.\\
  & &
  \hspace{4.5cm}\left. +\mathbf 1_{\ell_k = 1}
  \sum_{\pi\in\Pi_{k - 1}}\sigma^{(|\pi|)}(X_{\max(\mathbf s)})
  \prod_{J\in\pi}\mathbf D_{\mathbf s_J}^{(|J|),\mathbf l_J}X_{\max(\mathbf s)}
  \right]\\
  & &
  \hspace{0.75cm} +
  \mathbf 1_{k > 2}\sum_{j = 2}^{k - 1}\mathbf 1_{\ell_j = 1}\sum_{\pi\in\Pi_{k\backslash j}}
  \sigma^{(|\pi|)}(X_{\max\{s_1,\dots,s_j\}})
  \prod_{J\in\pi}\mathbf D_{\mathbf s_J}^{(|J|),\mathbf l_J}X_{\max\{s_1,\dots,s_j\}}.
 \end{eqnarray*}
 \item $\mathbf Q(k)$: For every $p\geqslant 1$ and $\mathbf l\in\{1,2\}^{k + 1}$,
 \begin{displaymath}
 \sup_{\mathbf r\in [0,T]^{k + 1}}
 \mathbb E\left(\sup_{u\in [\max(\mathbf r),T]}
 |\mathbf D_{\mathbf r}^{(k + 1),\mathbf l}X_u|^p\right) <\infty.
 \end{displaymath}
\end{itemize}
First, for $k = 1$ and every $t\in [0,T]$, as established in the proof of Lemma \ref{weak_differentiability_X}, for any $\ell\in\{1,2\}$ and every $s\in [t,T]$,
\begin{eqnarray*}
 \mathbf D_{s}^{\ell}X_t & = &
 \sigma(X_s)\mathbf 1_{\ell = 1} +
 \int_{s}^{t}b'(Y_u)\mathbf D_{s}^{\ell}Y_udu\\
 & &
 \hspace{3cm} +
 \int_{s}^{t}a'(X_u)\mathbf D_{s}^{\ell}X_udu +
 \int_{s}^{t}\sigma'(X_u)\mathbf D_{s}^{\ell}X_udW_1(u)\\
 & = &
 \alpha_{\ell}(s,t) +\beta_{\ell}(s,t) +
 \int_{s}^{t}a'(X_u)\mathbf D_{s}^{\ell}X_udu +
 \int_{s}^{t}\sigma'(X_u)\mathbf D_{s}^{\ell}X_udW_1(u).
\end{eqnarray*}
Then, since $\sigma$ is bounded, $b,\sigma\in C^{\infty}(\mathbb R)$ with all derivatives bounded, and since $Y\in\mathbb H^{\infty}$,
\begin{displaymath}
\sup_{r\in [0,T]}\mathbb E\left(
\sup_{u\in [r,T]}|\alpha_{\ell}(r,u) +\beta_{\ell}(r,u)|^p\right) <\infty
\textrm{ $;$ }\forall p\geqslant 1,
\end{displaymath}
and for every $\ell_1,\ell_2\in\{1,2\}$, $s_1,s_2\in [0,T]$ and $t\in [s_1\vee s_2,T]$,
\begin{eqnarray*}
 \mathbf D_{s_2}^{\ell_2}[\alpha_{\ell_1}(s_1,t) +\beta_{\ell_1}(s_1,t)]
 & = &
 \sigma'(X_{s_1})\mathbf D_{s_2}^{\ell_2}X_{s_1}\mathbf 1_{\ell_1 = 1} +
 \int_{s_1\vee s_2}^{t}b''(Y_u)\mathbf D_{s_1}^{\ell_1}Y_u\mathbf D_{s_2}^{\ell_2}Y_udu\\
 & &
 \hspace{5cm} +
 \int_{s_1\vee s_2}^{t}b'(Y_u)\mathbf D_{s_1,s_2}^{(2),\ell_1,\ell_2}Y_udu,
\end{eqnarray*}
leading - together with (\ref{general_sufficient_condition_AC_3}) - to
\begin{displaymath}
\sup_{r_1,r_2\in [0,T]}\mathbb E\left(\sup_{u\in [r_1\vee r_2,T]}
|\mathbf D_{r_2}^{\ell_2}[\alpha_{\ell_1}(r_1,u) +\beta_{\ell_1}(r_1,u)]|^p\right) <\infty
\textrm{ $;$ }\forall p\geqslant 1.
\end{displaymath}
So, since $a,\sigma\in C^{\infty}(\mathbb R)$ with all derivatives bounded, by Nualart \cite{NUALART06}, Lemma 2.2.2, for every $p\geqslant 1$ and $\ell_1,\ell_2\in\{1,2\}$,
\begin{displaymath}
\sup_{r_1,r_2\in [0,T]}\mathbb E\left(\sup_{u\in [r_1\vee r_2,T]}
|\mathbf D_{r_1,r_2}^{(2),\ell_1,\ell_2}X_u|^p\right) <\infty.
\end{displaymath}
Now, consider $k\in\mathbb N^*$ such that $\mathbf P(j)$ and $\mathbf Q(j)$ are true for every $j\in\{1,\dots,k\}$. Then, for any $\mathbf l = (\ell_1,\dots,\ell_k)\in\{1,2\}^k$, $\ell\in\{1,2\}$, $\mathbf s = (s_1,\dots,s_k)\in\mathbb [0,T]^k$, $s\in [0,T]$ and $t\in [\max(\mathbf s)\vee s,T]$,
\begin{eqnarray*}
 \mathbf D_{\mathbf s,s}^{(k + 1),\mathbf l,\ell}X_t
 & = &
 \mathbf D_{s}^{\ell}\alpha_{\bf l}(\mathbf s,t) +
 \mathbf D_{s}^{\ell}\beta_{\bf l}(\mathbf s,t) +
 \sigma'(X_{\max(\mathbf s)\vee s})
 \mathbf D_{\bf s}^{(k),\mathbf l}X_{\max(\mathbf s)\vee s}\\
 & &
 \hspace{0.75cm} +
 \int_{\max(\mathbf s)\vee s}^{t}
 \mathbf D_{s}^{\ell}[a'(X_u)
 \mathbf D_{\bf s}^{(k),\mathbf l}X_u]du +
 \int_{\max(\mathbf s)\vee s}^{t}
 \mathbf D_{s}^{\ell}[\sigma'(X_u)
 \mathbf D_{\bf s}^{(k),\mathbf l}X_u]dW_1(u)\\
 & = &
 A_{\mathbf s,s}^{\mathbf l,\ell}(t) + B_{\mathbf s,s}^{\mathbf l,\ell}(t) +
 \int_{\max(\mathbf s)\vee s}^{t}
 a'(X_u)\mathbf D_{\mathbf s,s}^{(k + 1),\mathbf l,\ell}X_udu\\
 & &
 \hspace{5cm} +
 \int_{\max(\mathbf s)\vee s}^{t}
 \sigma'(X_u)\mathbf D_{\bf s}^{(k + 1),\mathbf l,\ell}X_udW_1(u),
\end{eqnarray*}
where
\begin{eqnarray*}
 A_{\mathbf s,s}^{\mathbf l,\ell}(t) & := &
 \mathbf D_{s}^{\ell}\alpha_{\bf l}(\mathbf s,t) - R_{\mathbf s,s}^{\mathbf l,\ell}(t) +
 \int_{\max(\mathbf s)\vee s}^{t}
 a''(X_u)\mathbf D_{\bf s}^{(k),\mathbf l}X_u\mathbf D_{s}^{\ell}X_udu\\
 & &
 \hspace{6cm} +
 \int_{\max(\mathbf s)\vee s}^{t}
 \sigma''(X_u)\mathbf D_{\bf s}^{(k),\mathbf l}X_u\mathbf D_{s}^{\ell}X_udW_1(u)\\
 & &
 \hspace{0.5cm}{\rm with}\quad
 R_{\mathbf s,s}^{\mathbf l,\ell}(t) =\mathbf 1_{k > 1,\ell = 1}\sum_{\pi\in\Pi_{k}^{*}}
 \sigma^{(|\pi|)}(X_{\max(\mathbf s)\vee s})\prod_{J\in\pi}
 \mathbf D_{\mathbf s_J}^{(|J|),\mathbf l_J}X_{\max(\mathbf s)\vee s},
\end{eqnarray*}
and
\begin{displaymath}
B_{\mathbf s,s}^{\mathbf l,\ell}(t) :=
\mathbf D_{s}^{\ell}\beta_{\bf l}(\mathbf s,t) + R_{\mathbf s,s}^{\mathbf l,\ell}(t) +
\sigma'(X_{\max(\mathbf s)\vee s})
\mathbf D_{\bf s}^{(k),\mathbf l}X_{\max(\mathbf s)\vee s}\mathbf 1_{\ell = 1}.
\end{displaymath}
On the one hand, note that if
\begin{displaymath}
{\bf (A)}\quad
A_{\mathbf s,s}^{\mathbf l,\ell}(t) =\alpha_{\mathbf l,\ell}(\mathbf s,s,t)
\quad {\rm and}\quad
{\bf (B)}\quad
B_{\mathbf s,s}^{\mathbf l,\ell}(t) =\beta_{\mathbf l,\ell}(\mathbf s,s,t),
\end{displaymath}
then $\mathbf P(k + 1)$ is true. So, let us establish \textbf{(A)} and \textbf{(B)}.
\begin{itemize}
 \item[\textbf{(A)}] For $k > 1$ and every $u\in [\max(\mathbf s)\vee s,t]$,
 \begin{eqnarray*}
  & &
  \mathbf D_{s}^{\ell}\left[\sum_{\pi\in\Pi_{k}^{*}}
  a^{(|\pi|)}(X_u)\prod_{J\in\pi}\mathbf D_{\mathbf s_J}^{(|J|),\mathbf l_J}X_u\right]\\
  & &
  \hspace{3cm} =
  \sum_{\pi\in\Pi_{k}^{*}}
  a^{(|\pi| + 1)}(X_u)\mathbf D_{s}^{\ell}X_u
  \prod_{J\in\pi}\mathbf D_{\mathbf s_J}^{(|J|),\mathbf l_J}X_u\\
  & &
  \hspace{6cm} +
  \sum_{\pi\in\Pi_{k}^{*}}
  a^{(|\pi|)}(X_u)\sum_{I\in\pi}\mathbf D_{\mathbf s_I,s}^{(|I| + 1),\mathbf l_I,\ell}X_u
  \prod_{J\in\pi\backslash\{I\}}\mathbf D_{\mathbf s_J}^{(|J|),\mathbf l_J}X_u\\
  & &
  \hspace{3cm} =
  \sum_{\pi\in\Pi_{k}^{*}}
  a^{(|\pi| + 1)}(X_u)
  \prod_{J\in\pi\cup\{\{k + 1\}\}}\mathbf D_{(\mathbf s,s)_J}^{(|J|),(\mathbf l,\ell)_J}X_u\\
  & &
  \hspace{6cm} +
  \sum_{\pi\in\Pi_{k}^{*}}
  a^{(|\pi|)}(X_u)\sum_{I\in\pi}\prod_{J\in (\pi\backslash\{I\})\cup\{I\cup\{k + 1\}\}}
  \mathbf D_{(\mathbf s,s)_J}^{(|J|),(\mathbf l,\ell)_J}X_u\\
  & &
  \hspace{3cm} =
  \sum_{\overline\pi\in\mathcal U_{k}^{*}\cup\mathcal V_{k}^{*}}
  a^{(|\overline\pi|)}(X_u)
  \prod_{J\in\overline\pi}\mathbf D_{(\mathbf s,s)_J}^{(|J|),(\mathbf l,\ell)_J}X_u,
 \end{eqnarray*}
 where
 \begin{eqnarray*}
  \mathcal U_{k}^{*} & := &
  \{\pi\cup\{\{k + 1\}\}\textrm{ $;$ }\pi\in\Pi_{k}^{*}\}\\
  & &
  \hspace{2cm}{\rm and}\quad
  \mathcal V_{k}^{*} :=
  \{(\pi\backslash\{I\})\cup\{I\cup\{k + 1\}\}
  \textrm{ $;$ }\pi\in\Pi_{k}^{*}\textrm{, }I\in\pi\}.
 \end{eqnarray*}
 Moreover, consider
 \begin{eqnarray*}
  \mathcal U_k & := &
  \{\pi\cup\{\{k + 1\}\}\textrm{ $;$ }\pi\in\Pi_k\}\\
  & = &
  \mathcal U_{k}^{*}\cup\{\pi_{k,k + 1}\}
  \quad {\rm with}\quad
  \pi_{k,k + 1} =\{\{1,\dots,k\},\{k + 1\}\}.
 \end{eqnarray*}
 Then,
 \begin{eqnarray*}
  A_{\mathbf s,s}^{\mathbf l,\ell}(t;a) & := &
  \int_{\max(\mathbf s)\vee s}^{t}\mathbf D_{s}^{\ell}\left[\mathbf 1_{k > 1}\sum_{\pi\in\Pi_{k}^{*}}
  a^{(|\pi|)}(X_u)\prod_{J\in\pi}\mathbf D_{\mathbf s_J}^{(|J|),\mathbf l_J}X_u\right]du\\
  & &
  \hspace{4.5cm} +
  \int_{\max(\mathbf s)\vee s}^{t}
  a''(X_u)\mathbf D_{\bf s}^{(k),\mathbf l}X_u\mathbf D_{s}^{\ell}X_udu\\
  & = &
  \int_{\max(\mathbf s)\vee s}^{t}\left[\mathbf 1_{k > 1}
  \sum_{\overline\pi\in\mathcal U_{k}^{*}\cup\mathcal V_{k}^{*}}
  a^{(|\overline\pi|)}(X_u)
  \prod_{J\in\overline\pi}\mathbf D_{(\mathbf s,s)_J}^{(|J|),(\mathbf l,\ell)_J}X_u\right]du\\
  & &
  \hspace{4.5cm} +
  \int_{\max(\mathbf s)\vee s}^{t}
  a^{(|\pi_{k,k + 1}|)}(X_u)\mathbf D_{\bf s}^{(k),\mathbf l}X_u\mathbf D_{s}^{\ell}X_udu\\
  & = &
  \int_{\max(\mathbf s)\vee s}^{t}\left[
  \sum_{\overline\pi\in\mathcal U_k\cup\mathcal V_{k}^{*}}
  a^{(|\overline\pi|)}(X_u)
  \prod_{J\in\overline\pi}\mathbf D_{(\mathbf s,s)_J}^{(|J|),(\mathbf l,\ell)_J}X_u\right]du.
 \end{eqnarray*}
 In the same way,
 \begin{eqnarray*}
  A_{\mathbf s,s}^{\mathbf l,\ell}(t;\sigma) & := &
  \int_{\max(\mathbf s)\vee s}^{t}\mathbf D_{s}^{\ell}\left[\mathbf 1_{k > 1}\sum_{\pi\in\Pi_{k}^{*}}
  \sigma^{(|\pi|)}(X_u)\prod_{J\in\pi}\mathbf D_{\mathbf s_J}^{(|J|),\mathbf l_J}X_u\right]dW_1(u)\\
  & &
  \hspace{4.5cm} +
  \int_{\max(\mathbf s)\vee s}^{t}
  \sigma''(X_u)\mathbf D_{\bf s}^{(k),\mathbf l}X_u\mathbf D_{s}^{\ell}X_udW_1(u)\\
  & = &
  \int_{\max(\mathbf s)\vee s}^{t}\left[
  \sum_{\overline\pi\in\mathcal U_k\cup\mathcal V_{k}^{*}}
  \sigma^{(|\overline\pi|)}(X_u)
  \prod_{J\in\overline\pi}\mathbf D_{(\mathbf s,s)_J}^{(|J|),(\mathbf l,\ell)_J}X_u\right]dW_1(u),
 \end{eqnarray*}
 and
 \begin{eqnarray*}
  A_{\mathbf s,s}^{\mathbf l,\ell}(t;b) & := &
  \int_{\max(\mathbf s)\vee s}^{t}\mathbf D_{s}^{\ell}\left(\sum_{\pi\in\Pi_k}
  b^{(|\pi|)}(Y_u)\prod_{J\in\pi}\mathbf D_{\mathbf s_J}^{(|J|),\mathbf l_J}Y_u\right)du\\
  & = &
  \int_{\max(\mathbf s)\vee s}^{t}\left(
  \sum_{\overline\pi\in\mathcal U_k\cup\mathcal V_k}
  b^{(|\overline\pi|)}(Y_u)
  \prod_{J\in\overline\pi}\mathbf D_{(\mathbf s,s)_J}^{(|J|),(\mathbf l,\ell)_J}Y_u\right)du
 \end{eqnarray*}
 with
 \begin{displaymath}
 \mathcal V_k =
 \{(\pi\backslash\{I\})\cup\{I\cup\{k + 1\}\}
 \textrm{ $;$ }\pi\in\Pi_k\textrm{, }I\in\pi\}.
 \end{displaymath}
 Therefore, since $\Pi_{k + 1} =\mathcal U_k\cup\mathcal V_k$ and $\Pi_{k + 1}^{*} =\mathcal U_k\cup\mathcal V_{k}^{*}$,
 \begin{eqnarray*}
  A_{\mathbf s,s}^{\mathbf l,\ell}(t) & = &
  A_{\mathbf s,s}^{\mathbf l,\ell}(t;b) +
  A_{\mathbf s,s}^{\mathbf l,\ell}(t;a) +
  A_{\mathbf s,s}^{\mathbf l,\ell}(t;\sigma)\\
  & = &
  \int_{\max(\mathbf s)\vee s}^{t}
  \sum_{\overline\pi\in\Pi_{k + 1}}
  b^{(|\overline\pi|)}(Y_u)
  \prod_{J\in\overline\pi}\mathbf D_{(\mathbf s,s)_J}^{(|J|),(\mathbf l,\ell)_J}Y_udu\\
  & &
  \hspace{1.5cm} +
  \int_{\max(\mathbf s)\vee s}^{t}
  \sum_{\overline\pi\in\Pi_{k + 1}^{*}}
  a^{(|\overline\pi|)}(X_u)
  \prod_{J\in\overline\pi}\mathbf D_{(\mathbf s,s)_J}^{(|J|),(\mathbf l,\ell)_J}X_udu\\
  & &
  \hspace{3cm} +
  \int_{\max(\mathbf s)\vee s}^{t}
  \sum_{\overline\pi\in\Pi_{k + 1}^{*}}
  \sigma^{(|\overline\pi|)}(X_u)
  \prod_{J\in\overline\pi}\mathbf D_{(\mathbf s,s)_J}^{(|J|),(\mathbf l,\ell)_J}X_udW_1(u)\\
  & = &
  \alpha_{\mathbf l,\ell}(\mathbf s,s,t).
 \end{eqnarray*}
 \item[\textbf{(B)}] For the sake of simplicity, assume that $k > 2$. Then,
 \begin{displaymath}
 \beta_{\bf l}(\mathbf s,t) =
 \sum_{j = 1}^{k}\beta_{\mathbf l,j}(\mathbf s,t)\mathbf 1_{\ell_j = 1}
 \end{displaymath}
 where, for every $j\in\{1,\dots,k\}$,
 \begin{displaymath}
 \beta_{\mathbf l,j}(\mathbf s,t) :=
 \sum_{\pi\in\Pi_{k\backslash j}}
 \sigma^{(|\pi|)}(X_{\overline{\bf s}_j})
 \prod_{J\in\pi}\mathbf D_{\mathbf s_J}^{(|J|),\mathbf l_J}X_{\overline{\bf s}_j}
 \quad {\rm and}\quad
 \overline{\bf s}_j :=\max\{s_1,\dots,s_j\}.
 \end{displaymath}
 By following the same line as in the proof of \textbf{(A)}, for every $j\in\{1,\dots,k\}$,
 \begin{displaymath}
 \mathbf D_{s}^{\ell}\beta_{\mathbf l,j}(\mathbf s,t) =
 \sum_{\overline\pi\in\mathcal U_{k\backslash j}\cup\mathcal V_{k\backslash j}}
 \sigma^{(|\overline\pi|)}(X_{\overline{\bf s}_j})
 \prod_{J\in\overline\pi}\mathbf D_{(\mathbf s,s)_J}^{(|J|),(\mathbf l,\ell)_J}X_{\overline{\bf s}_j}
 \end{displaymath}
 with
 \begin{eqnarray*}
  \mathcal U_{k\backslash j} & = &
  \{\pi\cup\{\{k + 1\}\}\textrm{ $;$ }\pi\in\Pi_{k\backslash j}\}\\
  & &
  \hspace{2cm}{\rm and}\quad
  \mathcal V_{k\backslash j} =
  \{(\pi\backslash\{I\})\cup\{I\cup\{k + 1\}\}
  \textrm{ $;$ }\pi\in\Pi_{k\backslash j}\textrm{, }I\in\pi\}.
 \end{eqnarray*}
 Moreover,
 \begin{displaymath}
 \Pi_{(k + 1)\backslash (k + 1)} =\Pi_k =\Pi_{k}^{*}\cup\{\pi_k\}
 \quad {\rm with}\quad\pi_k =\{\{1,\dots,k\}\},
 \end{displaymath}
 and
 \begin{displaymath}
 \Pi_{(k + 1)\backslash j} =
 \mathcal U_{k\backslash j}\cup\mathcal V_{k\backslash j}
 \textrm{ $;$ }\forall j\in\{1,\dots,k\}.
 \end{displaymath}
 Therefore,
 \begin{eqnarray*}
  B_{\mathbf s,s}^{\mathbf l,\ell}(t)
  & = &
  \sum_{j = 1}^{k}\mathbf 1_{\ell_j = 1}
  \sum_{\overline\pi\in\mathcal U_{k\backslash j}\cup\mathcal V_{k\backslash j}}
  \sigma^{(|\overline\pi|)}(X_{\overline{\bf s}_j})
  \prod_{J\in\overline\pi}\mathbf D_{(\mathbf s,s)_J}^{(|J|),(\mathbf l,\ell)_J}X_{\overline{\bf s}_j}\\
  & &
  \hspace{1.5cm} +
  \mathbf 1_{\ell = 1}\sum_{\pi\in\Pi_{k}^{*}}
  \sigma^{(|\pi|)}(X_{\max(\mathbf s)\vee s})\prod_{J\in\pi}
  \mathbf D_{\mathbf s_J}^{(|J|),\mathbf l_J}X_{\max(\mathbf s)\vee s}\\
  & &
  \hspace{3cm} +
  \mathbf 1_{\ell = 1}
  \sigma^{(|\pi_k|)}(X_{\max(\mathbf s)\vee s})
  \mathbf D_{\bf s}^{(k),\mathbf l}X_{\max(\mathbf s)\vee s}\\
  & = &
  \sum_{j = 1}^{k}\mathbf 1_{\ell_j = 1}
  \sum_{\overline\pi\in\Pi_{(k + 1)\backslash j}}
  \sigma^{(|\overline\pi|)}(X_{\overline{\bf s}_j})
  \prod_{J\in\overline\pi}\mathbf D_{(\mathbf s,s)_J}^{(|J|),(\mathbf l,\ell)_J}X_{\overline{\bf s}_j}\\
  & &
  \hspace{1.5cm} +
  \mathbf 1_{\ell = 1}\sum_{\pi\in\Pi_k}
  \sigma^{(|\pi|)}(X_{\max(\mathbf s)\vee s})\prod_{J\in\pi}
  \mathbf D_{\mathbf s_J}^{(|J|),\mathbf l_J}X_{\max(\mathbf s)\vee s} =
  \beta_{\mathbf l,\ell}(\mathbf s,s,t).
 \end{eqnarray*}
\end{itemize}
On the other hand, let us show that $\mathbf Q(k + 1)$ is true. By $\mathbf Q(1),\dots,\mathbf Q(k - 1)$, since $\sigma$ is bounded, $a,b,\sigma\in C^{\infty}(\mathbb R)$ with all derivatives bounded, and since $Y\in\mathbb H^{\infty}$,
\begin{displaymath}
\sup_{\mathbf r\in [0,T]^{k + 1}}\mathbb E\left(
\sup_{u\in[\max(\mathbf r),T]}|\alpha_{\mathbf l,\ell}(\mathbf r,u) +
\beta_{\mathbf l,\ell}(\mathbf r,u)|^p\right) <\infty
\textrm{ $;$ }\forall p\geqslant 1.
\end{displaymath}
Consider $\mathbf l_1\in\{1,2\}^{k + 1}$ and $\ell_2\in\{1,2\}$. For every $\mathbf s_1\in [0,T]^{k + 1}$, $s_2\in [0,T]$ and $t\in [\max(\mathbf s_1)\vee s_2,T]$, by following the same line as in the proofs of \textbf{(A)} and \textbf{(B)}, $\mathbf D_{s_2}^{\ell_2}\alpha_{\mathbf l_1}(\mathbf s_1,t)$ and $\mathbf D_{s_2}^{\ell_2}\beta_{\mathbf l_1}(\mathbf s_1,t)$ can be written as integrals of multilinear maps in
\begin{displaymath}
\mathbf D^{(1)}X,\dots,\mathbf D^{(k + 1)}X,
\mathbf D^{(1)}Y,\dots,\mathbf D^{(k + 2)}Y,\quad
\textrm{but not depending on $\mathbf D^{(k + 2)}X$.}
\end{displaymath}
Then, by $\mathbf Q(1),\dots,\mathbf Q(k)$, since $\sigma$ is bounded, $a,b,\sigma\in C^{\infty}(\mathbb R)$ with all derivatives bounded, and since $Y\in\mathbb H^{\infty}$,
\begin{displaymath}
\sup_{(\mathbf r_1,r_2)\in [0,T]^{k + 1}\times [0,T]}\mathbb E\left(
\sup_{u\in[\max(\mathbf r_1)\vee r_2,T]}
|\mathbf D_{r_2}^{\ell_2}[\alpha_{\mathbf l_1}(\mathbf r_1,u) +
\beta_{\mathbf l_1}(\mathbf r_1,u)]|^p\right) <\infty
\textrm{ $;$ }\forall p\geqslant 1.
\end{displaymath}
So, by Nualart \cite{NUALART06}, Lemma 2.2.2,
\begin{displaymath}
\sup_{\mathbf r\in [0,T]^{k + 2}}
\mathbb E\left(\sup_{u\in [\max(\mathbf r),T]}
|\mathbf D_{\mathbf r}^{(k + 2),(\mathbf l_1,\ell_2)}X_u|^p\right) <\infty
\textrm{ $;$ }\forall p\geqslant 1.
\end{displaymath}
{\it Step 2.3 (conclusion).} Consider $t\in (0,T]$. By the two previous steps, and by Nualart \cite{NUALART06}, Inequality (2.32) and Proposition 2.1.5, the distribution of $(X_t,Y_t)$ has a smooth density $f_t$ with respect to the Lebesgue measure on $\mathbb R^2$, and there exist $\mathfrak c_3 > 0$ and $(m,\alpha),(N,p)\in\mathbb N^*\times (2,\infty)$, depending on $T$ but not on $t$, such that for every $(x,y)\in\mathbb R^2$,
\begin{eqnarray*}
 f_t(x,y)
 & \leqslant &
 \mathfrak c_3
 \mathbb E\left(\left|\frac{1}{\det(\Gamma_t)}\right|^{\alpha}\right)^{\frac{m}{\alpha}}
 \|\mathbf D(X_t,Y_t)\|_{\mathbb D^{2,p}}^{N}
 \mathbb P(|Y_t| > |y|)^{\frac{1}{2}}
 \mathbb P(|X_t - x_0| > |x - x_0|)^{\frac{1}{2}}\\
 & \leqslant &
 \frac{\mathfrak c_4\mathfrak c_2(\alpha)^{m/\alpha}}{t^m}
 \mathbb E\left(\frac{1}{\|\mathbf DY_t\|_{T}^{4\alpha}}\right)^{\frac{m}{2\alpha}}
 \Pi_t(x,y)
\end{eqnarray*}
with
\begin{eqnarray*}
 \mathfrak c_{4}^{\frac{p}{N}} & = &
 \mathfrak c_3\sum_{k = 0}^{3}\sum_{j = 1}^{p}\left[
 \sup_{r_1,\dots,r_k\in [0,T]}\mathbb E\left(
 \sup_{u\in [r_1\vee\dots\vee r_k,T]}
 \|\mathbf D_{r_1,\dots,r_k}^{(k)}X_u\|_{j,\mathbb R^{2^k}}^{j}\right)\right.\\
 & &
 \hspace{4cm}
 \left. +
 \sup_{r_1,\dots,r_k\in [0,T]}\mathbb E\left(
 \sup_{u\in [r_1\vee\dots\vee r_k,T]}
 \|\mathbf D_{r_1,\dots,r_k}^{(k)}Y_u\|_{j,\mathbb R^{2^k}}^{j}\right)\right]
 <\infty.
\end{eqnarray*}
%


%
\subsection{Proof of Proposition \ref{integrability_Pi}}\label{section_proof_integrability_Pi}
The proof is dissected in two steps.
\\
\\
{\bf Step 1.} Let $M = (M_t)_{t\in [0,T]}$ be the $\mathbb F$-martingale defined by
\begin{displaymath}
M_t :=\int_{0}^{t}\sigma(X_s)dW_1(s)
\textrm{ $;$ }\forall t\in [0,T].
\end{displaymath}
Consider $t\in (0,T]$. First, since $\sigma$ is bounded,
\begin{displaymath}
\langle M\rangle_t =
\int_{0}^{t}\sigma(X_s)^2ds\leqslant
\|\sigma\|_{\infty}^{2}t,
\end{displaymath}
and since $(a,b)$ is bounded, one may consider $\gamma :=\|a\|_{\infty} +\|b\|_{\infty}$. Now, by the Bernstein inequality for local martingales (see Revuz and Yor \cite{RY99}, p. 153), for every $x\in\mathbb R$ such that $|x - x_0|\geqslant\gamma t$,
\begin{eqnarray*}
 \mathbb P(|X_t - x_0| > |x - x_0|) & = &
 \mathbb P\left(\langle M\rangle_t\leqslant\|\sigma\|_{\infty}^{2}t,
 \left|M_t +\int_{0}^{t}(a(X_s) + b(Y_s))ds\right| > |x - x_0|\right)\\
 & \leqslant &
 \mathbb P\left[\langle M\rangle_t\leqslant\|\sigma\|_{\infty}^{2}t,
 \sup_{s\in [0,t]}|M_s| > |x - x_0| -\gamma t\right]\\
 & \leqslant &
 2\exp\left(-\frac{(|x - x_0| -\gamma t)^2}{2\|\sigma\|_{\infty}^{2}t}\right)
 \leqslant
 2\exp\left(-\frac{(x - x_0)^2}{2\|\sigma\|_{\infty}^{2}t}
 \left(1 -\frac{2\gamma t}{|x - x_0|}\right)\right).
\end{eqnarray*}
So, for every $x\in\mathbb R$ such that $|x - x_0|\geqslant (2\gamma + 1)t >\gamma t$,
\begin{displaymath}
\mathbb P(|X_t - x_0| > |x - x_0|)\leqslant
2\exp\left(-\mathfrak c_1\frac{(x - x_0)^2}{t}\right)
\quad {\rm with}\quad
\mathfrak c_1 =\frac{1}{2(2\gamma + 1)\|\sigma\|_{\infty}^{2}},
\end{displaymath}
Therefore, there exist $\mathfrak c_2,\mathfrak c_3 > 0$, depending on $T$ but not on $t$, such that for every $(x,y)\in\mathbb R^2$,
\begin{displaymath}
\Pi_t(x,y)\leqslant
\mathfrak c_2\mathbb P(|Y_t| > |y|)^{\frac{1}{2}}
\exp\left(-\mathfrak c_3\frac{(x - x_0)^2}{t}\right).
\end{displaymath}
{\bf Step 2.} Under the additional condition (\ref{integrability_Pi_1}),
\begin{eqnarray*}
 \sup_{x\in\mathbb R}\int_{[t_0,T]\times\mathbb R}f_s(x,y)dsdy
 & \leqslant &
 \frac{\mathfrak c_{\ref{general_sufficient_condition_AC}}}{t_{0}^{m}}
 \int_{[t_0,T]\times\mathbb R}
 \mathbb E\left(\frac{1}{\|\mathbf D^2Y_s\|_{T}^{4\alpha}}\right)^{\frac{m}{2\alpha}}
 \mathbb P(|Y_s| > |y|)^{\frac{1}{2}}dsdy\\
 & \leqslant &
 \mathfrak c_{\ref{general_sufficient_condition_AC}}
 \frac{T - t_0}{t_{0}^{m}}
 \int_{-\infty}^{\infty}
 \mu_{t_0,T}(y)dy <\infty
\end{eqnarray*}
and, by Step 1,
\begin{eqnarray*}
 \sup_{y\in\mathbb R}\int_{[t_0,T]\times\mathbb R}f_s(x,y)dsdx
 & \leqslant &
 \frac{\mathfrak c_2
 \mathfrak c_{\ref{general_sufficient_condition_AC}}}{t_{0}^{m}}
 \int_{[t_0,T]\times\mathbb R}
 \mathbb E\left(\frac{1}{\|\mathbf D^2Y_s\|_{T}^{4\alpha}}\right)^{\frac{m}{2\alpha}}
 \mathbb P(|Y_s| > |y|)^{\frac{1}{2}}\\
 & &
 \hspace{5.5cm}\times
 \exp\left(-\mathfrak c_3\frac{(x - x_0)^2}{s}\right)dsdx\\
 & \leqslant &
 \mathfrak c_2
 \mathfrak c_{\ref{general_sufficient_condition_AC}}
 \frac{T - t_0}{t_{0}^{m}}\|\mu_{t_0,T}\|_{\infty}
 \int_{-\infty}^{\infty}
 \exp\left(-\mathfrak c_3\frac{x^2}{T}\right)dx <\infty.
\end{eqnarray*}
%


%
\subsection{Proof of Proposition \ref{sufficient_condition_AC_diffusion}}\label{section_proof_sufficient_condition_AC_diffusion}
The process $\mathbf X = (X,Y)$ is the solution of the two-dimensional stochastic differential equation
\begin{equation}\label{sufficient_condition_AC_diffusion_1}
\mathbf X_t =\mathbf X_0 +
\int_{0}^{t}\Phi(\mathbf X_s)ds +
\int_{0}^{t}\Psi(\mathbf X_s)dW_s
\textrm{ $;$ }t\in [0,T]
\end{equation}
where, for every $(x,y)\in\mathbb R^2$,
\begin{displaymath}
\Phi(x,y) :=
\begin{pmatrix}
 a(x) + b(y)\\
 \mu(y)
\end{pmatrix}
\quad {\rm and}\quad
\Psi(x,y) :=
\begin{pmatrix}
 \sigma(x) & 0\\
 0 & \kappa(y)
\end{pmatrix}.
\end{displaymath}
For every $\mathbf x,\mathbf v\in\mathbb R^2$,
\begin{eqnarray*}
 \langle\Psi(\mathbf x)\Psi(\mathbf x)^*\mathbf v,\mathbf v\rangle_{2,\mathbb R^2}
 & = &
 \mathbf v^*\Psi(\mathbf x)^2\mathbf v\\
 & = &
 \|\Psi(\mathbf x)\mathbf v\|_{2,\mathbb R^2}^{2} =
 (\sigma(\mathbf x_1)\mathbf v_1)^2 +
 (\kappa(\mathbf x_2)\mathbf v_2)^2,
\end{eqnarray*}
leading to
\begin{displaymath}
m^2\|\mathbf v\|_{2,\mathbb R^2}^{2}\leqslant
\langle\Psi(\mathbf x)\Psi(\mathbf x)^*\mathbf v,\mathbf v\rangle_{2,\mathbb R^2}
\leqslant
(\|\sigma\|_{\infty}\vee\|\kappa\|_{\infty})^2\|\mathbf v\|_{2,\mathbb R^2}^{2}
\end{displaymath}
with
\begin{displaymath}
m =\left(\inf_{x\in\mathbb R}|\sigma(x)|\right)\wedge
\left(\inf_{y\in\mathbb R}|\kappa(y)|\right) > 0.
\end{displaymath}
Then, $\Psi$ satisfies the non-degeneracy condition (1.5) in Menozzi et al. \cite{MPZ21}. Since, in addition, $\Phi :\mathbb R^2\rightarrow\mathbb R^2$ and $\Psi :\mathbb R^2\rightarrow\mathcal M_2(\mathbb R)$ are Lipschitz continuous maps, $Y$ and $(X,Y)$ fulfill Assumptions \ref{assumption_Y} and \ref{assumption_XY}, and satisfy (\ref{integrability_Pi_2}), by Menozzi et al. \cite{MPZ21}, Theorem 1.2.
%


%
\section{Proofs of statistical results (Section \ref{section_projection_LSE})}\label{section_proofs_statistical_results}
%


%
\subsection{Proof of Lemma \ref{bound_Omega}}\label{section_proof_bound_Omega}
The proof of Lemma \ref{bound_Omega} relies on the following matrix Chernov's inequality.
%


%
\begin{proposition}\label{matrix_Chernov_inequality}
Let $\mathbb X$ be a $d\times d$ positive semidefinite symmetric random matrix such that $\lambda_{\max}(\mathbb X)\leqslant R$ a.s, where $d\in\mathbb N^*$ and $R > 0$ is a deterministic constant. Consider $\mathbb G :=\mathbb X_1 +\dots +\mathbb X_n$, $\mu_{\min} :=\lambda_{\min}(\mathbb E(\mathbb G))$ and $\mu_{\max} :=\lambda_{\max}(\mathbb E(\mathbb G))$, where $n\in\mathbb N^*$ and $\mathbb X_1,\dots,\mathbb X_n$ are independent copies of $\mathbb X$. Then,
\begin{displaymath}
\mathbb P(
\lambda_{\min}(\mathbb G)
\leqslant (1 -\delta)\mu_{\min})
\leqslant d\left(\frac{e^{-\delta}}{(1 -\delta)^{1 -\delta}}
\right)^{\frac{\mu_{\min}}{R}}
\textrm{$;$ }\forall\delta\in [0,1],
\end{displaymath}
and
\begin{displaymath}
\mathbb P(\lambda_{\max}(\mathbb G)\geqslant (1 +\delta)\mu_{\max})
\leqslant d\left(\frac{e^{\delta}}{(1 +\delta)^{1 +\delta}}
\right)^{\frac{\mu_{\max}}{R}}
\textrm{$;$ }\forall\delta\geqslant 0.
\end{displaymath}
\end{proposition}
\noindent
See Tropp \cite{TROPP12}, Theorem 1.1 for a proof.
\\
\\
Note that
\begin{displaymath}
\Omega_{\bf m} =
\left\{\|\mathbb G - I_{m_1 + m_2}\|_{\rm op}\leqslant\frac{1}{2}\right\},
\end{displaymath}
where $\mathbb G :=\mathbb X_1 +\dots +\mathbb X_N$ with, for every $i\in\{1,\dots,N\}$,
\begin{displaymath}
\mathbb X_i =
\frac{1}{N}\Psi_{\bf m}^{-\frac{1}{2}}\mathbb X_{i}^{0}\Psi_{\bf m}^{-\frac{1}{2}}
\end{displaymath}
and
\begin{displaymath}
\mathbb X_{i}^{0} =
\begin{bmatrix}
 \displaystyle{\left(\frac{1}{T_0}\int_{t_0}^{T}
 \varphi_j(X_{s}^{i})\varphi_{j'}(X_{s}^{i})ds\right)_{j,j'}} &
 \displaystyle{\left(\frac{1}{T_0}\int_{t_0}^{T}
 \varphi_j(X_{s}^{i})\psi_k(Y_{s}^{i})ds\right)_{j,k}}\\
 \displaystyle{\left(\frac{1}{T_0}\int_{t_0}^{T}
 \varphi_j(X_{s}^{i})\psi_k(Y_{s}^{i})ds\right)_{k,j}} &
 \displaystyle{\left(\frac{1}{T_0}\int_{t_0}^{T}
 \psi_k(Y_{s}^{i})\psi_{k'}(Y_{s}^{i})ds\right)_{k,k'}}
\end{bmatrix}.
\end{displaymath}
In order to prove Lemma \ref{bound_Omega} by applying Proposition \ref{matrix_Chernov_inequality} to $\mathbb G$, let us determine $\mu_{\min}$, $\mu_{\max}$ and $R$. First, since the $\mathbb X_i$'s are i.i.d. positive semidefinite symmetric matrices such that $\mathbb E(\mathbb X_i) = N^{-1}I_{m_1 + m_2}$ for every $i\in\{1,\dots,N\}$, $\mu_{\min} =\mu_{\max} = 1$. Now, for every $i\in\{1,\dots,N\}$ and $\mathbf x\in\mathbb R^{m_1 + m_2}$, by setting $\mathbf y =\Psi_{\bf m}^{-1/2}\mathbf x$,
\begin{eqnarray*}
 \mathbf x^*\mathbb X_i\mathbf x
 & = &
 \frac{1}{NT_0}\int_{t_0}^{T}\left(
 \sum_{j = 1}^{m_1}y_j\varphi_j(X_{s}^{i}) +
 \sum_{k = 1}^{m_2}y_{m_1 + k}\psi_k(Y_{s}^{i})\right)^2ds\\
 & \leqslant &
 \frac{1}{NT_0}\underbrace{\left(\sum_{j = 1}^{m_1 + m_2}y_{j}^{2}\right)}_{=\|\mathbf y\|^2}
 \int_{t_0}^{T}\underbrace{\left(\sum_{j = 1}^{m_1}\varphi_j(X_{s}^{i})^{2} +
 \sum_{k = 1}^{m_2}
 \psi_k(Y_{s}^{i})^{2}\right)}_{\leqslant\mathfrak L_{\varphi}(m_1) +\mathfrak L_{\psi}(m_2)}ds
 \leqslant
 R\|\mathbf x\|^2
\end{eqnarray*}
with
\begin{displaymath}
R =
\frac{1}{N}\|\Psi_{\bf m}^{-1}\|_{\rm op}
(\mathfrak L_{\varphi}(m_1) +\mathfrak L_{\psi}(m_2)).
\end{displaymath}
Then, for every $i\in\{1,\dots,N\}$,
\begin{displaymath}
\lambda_{\max}(\mathbb X_i) =
\sup_{\|\mathbf x\|_{2,\mathbb R^{m_1 + m_2}} = 1}\mathbf x^*\mathbb X_i\mathbf x\leqslant R.
\end{displaymath}
So, by Proposition \ref{matrix_Chernov_inequality}, for every $\delta\in (0,1)$,
\begin{eqnarray}
 \mathbb P(\|\mathbb G - I_{m_1 + m_2}\|_{\rm op} >\delta)
 & \leqslant &
 \mathbb P(\lambda_{\min}(\mathbb G)\leqslant 1 -\delta) +
 \mathbb P(\lambda_{\max}(\mathbb G)\geqslant 1 +\delta)
 \nonumber\\
 \label{bound_Omega_1}
 & \leqslant &
 2(m_1 + m_2)\left(
 \frac{e^{-\delta}}{(1 -\delta)^{1 -\delta}}\right)^{\frac{1}{R}}.
\end{eqnarray}
In conclusion, by Inequality (\ref{bound_Omega_1}) with $\delta = 1/2$, and by Assumption \ref{stability_condition},
\begin{eqnarray*}
 \mathbb P(\Omega_{\bf m}^{c}) & = &
 \mathbb P\left(\|\mathbb G - I_{m_1 + m_2}\|_{\rm op} >\frac{1}{2}\right)\\
 & \leqslant &
 2(m_1 + m_2)(2^{\frac{1}{2}}e^{-\frac{1}{2}})^{\frac{1}{R}} =
 2(m_1 + m_2)\exp\left(-\frac{1}{2R}(1 -\log(2))\right)\\
 & \leqslant &
 2(m_1 + m_2)\exp\left(-\frac{\log(N)}{\mathfrak c_r}(1 -\log(2))\right) =
 \frac{2(m_1 + m_2)}{N^{r + 1}}.
\end{eqnarray*}
%


%
\subsection{Proof of Proposition \ref{nonadaptive_risk_bound}}\label{section_proof_nonadaptive_risk_bound}
First of all, one can show that
\begin{displaymath}
\Omega_{\bf m} =
\left\{\left|\frac{\|(\tau,\nu)\|_{N}^{2}}{\|(\tau,\nu)\|_{f}^{2}} - 1\right|\leqslant\frac{1}{2}\textrm{ $;$ }
\forall (\tau,\nu)\in S_{m_1}\times\Sigma_{m_2}\right\}
\end{displaymath}
by following the same line as Comte and Genon-Catalot in the beginning of the proof of \cite{CGC20}, Proposition 2.1, and that there exists a constant $\mathfrak c_1 > 0$, not depending on $\mathbf m$ and $N$, such that
\begin{equation}\label{nonadaptive_risk_bound_2}
\mathbb P(\Lambda_{\bf m}^{c})\leqslant
\mathbb P(\Omega_{\bf m}^{c})\leqslant
\frac{\mathfrak c_1}{N^r}
\end{equation}
with the same arguments as in the proof of \cite{CGC20}, Lemma 6.1. Note also that
\begin{equation}\label{nonadaptive_risk_bound_3}
\widehat{\bf Z}_{\bf m} =
\begin{pmatrix}
 (\langle (\varphi_j,0),(a,b)\rangle_N)_j\\
 (\langle (0,\psi_k),(a,b)\rangle_N)_k
\end{pmatrix}
+\mathbf E_{\bf m}
\quad {\rm with}\quad
\mathbf E_{\bf m} =
\begin{pmatrix}
 (\zeta_N(\varphi_j,0))_j\\
 (\zeta_N(0,\psi_k))_k
\end{pmatrix},
\end{equation}
where $\zeta_N$ is the centered empirical process defined by
\begin{displaymath}
\zeta_N(\tau,\nu) :=
\frac{1}{NT_0}\sum_{i = 1}^{N}\int_{t_0}^{T}
(\tau(X_{s}^{i}) +\nu(Y_{s}^{i}))\sigma(X_{s}^{i})dW_{1}^{i}(s),
\end{displaymath}
and $\langle .,.\rangle_N$ is the empirical inner product defined by
\begin{equation}\label{nonadaptive_risk_bound_4}
\langle (\tau_1,\nu_1),(\tau_2,\nu_2)\rangle_N :=
\frac{1}{NT_0}\sum_{i = 1}^{N}\int_{t_0}^{T}
(\tau_1(X_{s}^{i}) +\nu_1(Y_{s}^{i}))
(\tau_2(X_{s}^{i}) +\nu_2(Y_{s}^{i}))ds.
\end{equation}
To conclude these preliminaries, let us provide a suitable control of the second order moment of $\mathbf E_{\bf m}\mathbf 1_{\Omega_{\bf m}^{c}}$.
%


%
\begin{lemma}\label{2nd_order_moment_E}
Under the assumptions of Proposition \ref{nonadaptive_risk_bound}, there exists a constant $\mathfrak c_{\ref{2nd_order_moment_E}} > 0$, not depending on $\mathbf m$ and $N$, such that
\begin{displaymath}
\mathbb E(\|\mathbf E_{\bf m}\|_{2,\mathbb R^{m_1 + m_2}}^{2}\mathbf 1_{\Omega_{\bf m}^{c}})
\leqslant
\frac{\mathfrak c_{\ref{2nd_order_moment_E}}}{N^{(r - 1)/2}}.
\end{displaymath}
\end{lemma}
\noindent
The proof of Lemma \ref{2nd_order_moment_E} is postponed to Section \ref{section_proof_2nd_order_moment_E}.
%


%
\subsubsection{Steps of the proof}
The proof of Proposition \ref{nonadaptive_risk_bound} is dissected in two steps.
\\
\\
{\bf Step 1.} Let $(a_{m_1}^{N},b_{m_2}^{N})$ be the minimizer of $(u,v)\mapsto \|(u,v) - (a,b)\mathbf 1_A\|_{N}^{2}$ over $\mathcal S_{\bf m}$. Precisely,
\begin{displaymath}
a_{m_1}^{N} =\sum_{j = 1}^{m_1}\theta_{j}^{N}\varphi_j
\quad {\rm and}\quad
b_{m_2}^{N} =\sum_{k = 1}^{m_2}\theta_{m_1 + k}^{N}\psi_k
\end{displaymath}
with
\begin{equation}\label{nonadaptive_risk_bound_5}
\theta^N =
\underset{\theta\in\mathbb R^{m_1 + m_2} : h(\theta) = 0}{\rm argmin}\mathbb J_N(\theta),
\end{equation}
where
\begin{displaymath}
\mathbb J_N(\theta) :=
\frac{1}{NT_0}\sum_{i = 1}^{N}\int_{t_0}^{T}\left(
\sum_{j = 1}^{m_1}\theta_j\varphi_j(X_{s}^{i}) - (a\mathbf 1_{A_1})(X_{s}^{i}) +
\sum_{k = 1}^{m_2}\theta_{m_1 + k}\psi_k(Y_{s}^{i}) - (b\mathbf 1_{A_2})(Y_{s}^{i})\right)^2ds.
\end{displaymath}
Consider
\begin{displaymath}
\mathbf Z_{\bf m} =
\begin{pmatrix}
 (\langle (\varphi_j,0),(a,b)\rangle_N)_j\\
 (\langle (0,\psi_k),(a,b)\rangle_N)_k
\end{pmatrix},
\end{displaymath}
and let $\mathbb L_N$ be the Lagrangian for Problem (\ref{nonadaptive_risk_bound_5}):
\begin{displaymath}
\mathbb L_N(\theta,\lambda) :=
\mathbb J_N(\theta) -\lambda h(\theta)\textrm{ $;$ }
(\theta,\lambda)\in\mathbb R^{m_1 + m_2}\times\mathbb R.
\end{displaymath}
Necessarily,
\begin{displaymath}
\nabla\mathbb L_N(\theta^N,\lambda^N) =
\begin{pmatrix}
 2(\widehat\Psi_{\bf m}\theta^N -
 \mathbf Z_{\bf m}) -\lambda^N\mathbf d_{\bf m}\\
 -h(\theta^N)
\end{pmatrix} = 0,
\end{displaymath}
leading to
\begin{displaymath}
\theta^N =
\widehat\Psi_{\bf m}^{-1}\left(\mathbf Z_{\bf m} +
\frac{\lambda^N}{2}\mathbf d_{\bf m}\right),
\quad\textrm{and then}\quad
\lambda^N =
-2\cdot\frac{\langle
\widehat\Psi_{\bf m}^{-1}\mathbf Z_{\bf m},
\mathbf d_{\bf m}\rangle_{2,\mathbb R^{m_1 + m_2}}}{
\langle
\widehat\Psi_{\bf m}^{-1}\mathbf d_{\bf m},\mathbf d_{\bf m}
\rangle_{2,\mathbb R^{m_1 + m_2}}}.
\end{displaymath}
So,
\begin{displaymath}
\theta^N =
\widehat\Psi_{\bf m}^{-1}\mathbf Z_{\bf m} -
\frac{\mathbf d_{\bf m}^{*}\widehat\Psi_{\bf m}^{-1}\mathbf Z_{\bf m}}{
\mathbf d_{\bf m}^{*}\widehat\Psi_{\bf m}^{-1}\mathbf d_{\bf m}}\cdot
\widehat\Psi_{\bf m}^{-1}\mathbf d_{\bf m}.
\end{displaymath}
{\bf Step 2.} First, for any $(\tau,\nu)\in\mathcal S_{\bf m}$ and $t\in (0,1)$, since $t(\tau,\nu) + (1 - t)(a_{m_1}^{N},b_{m_2}^{N})$ belongs to $\mathcal S_{\bf m}$,
\begin{eqnarray*}
 \|(a_{m_1}^{N},b_{m_2}^{N}) - (a,b)\mathbf 1_A\|_{N}^{2}
 & \leqslant &
 \|t((\tau,\nu) - (a_{m_1}^{N},b_{m_2}^{N})) +
 (a_{m_1}^{N},b_{m_2}^{N}) - (a,b)\mathbf 1_A\|_{N}^{2}\\
 & = &
 t^2\|(\tau,\nu) - (a_{m_1}^{N},b_{m_2}^{N})\|_{N}^{2}\\
 & &
 \hspace{1cm} -
 2t\langle (\tau,\nu) - (a_{m_1}^{N},b_{m_2}^{N}),
 (a,b)\mathbf 1_A - (a_{m_1}^{N},b_{m_2}^{N})\rangle_N\\
 & &
 \hspace{1cm} +
 \|(a_{m_1}^{N},b_{m_2}^{N}) - (a,b)\mathbf 1_A\|_{N}^{2}.
\end{eqnarray*}
Then,
\begin{equation}\label{nonadaptive_risk_bound_6}
\langle (\tau,\nu) - (a_{m_1}^{N},b_{m_2}^{N}),
(a,b)\mathbf 1_A - (a_{m_1}^{N},b_{m_2}^{N})\rangle_N\leqslant
\frac{t}{2}\|(\tau,\nu) - (a_{m_1}^{N},b_{m_2}^{N})\|_{N}^{2}
\xrightarrow[t\rightarrow 0]{} 0.
\end{equation}
By applying Inequality (\ref{nonadaptive_risk_bound_6}) to $(\tau,\nu) = (\widehat a_{m_1},\widehat b_{m_2})$ and to $(\tau,\nu) = 2(a_{m_1}^{N},b_{m_2}^{N}) - (\widehat a_{m_1},\widehat b_{m_2})$,
\begin{displaymath}
\langle
(\widehat a_{m_1},\widehat b_{m_2}) - (a_{m_1}^{N},b_{m_2}^{N}),
(a,b)\mathbf 1_A - (a_{m_1}^{N},b_{m_2}^{N})\rangle_N = 0.
\end{displaymath}
So,
\begin{displaymath}
\|(\widehat a_{m_1},\widehat b_{m_2}) - (a,b)\mathbf 1_A\|_{N}^{2} =
\min_{(\tau,\nu)\in\mathcal S_{\bf m}}\|(\tau,\nu) - (a,b)\mathbf 1_A\|_{N}^{2} +
\|(\widehat a_{m_1},\widehat b_{m_2}) - (a_{m_1}^{N},b_{m_2}^{N})\|_{N}^{2},
\end{displaymath}
leading to
\begin{eqnarray*}
 \mathbb E(\|(\widetilde a_{m_1},\widetilde b_{m_2}) - (a,b)\mathbf 1_A\|_{N}^{2})
 & \leqslant &
 \mathbb E\left(\min_{(\tau,\nu)\in\mathcal S_{\bf m}}
 \|(\tau,\nu) - (a,b)\mathbf 1_A\|_{N}^{2}\right)\\
 & &
 \hspace{0.5cm} +
 \mathbb E(\|(\widehat a_{m_1},\widehat b_{m_2}) - (a_{m_1}^{N},b_{m_2}^{N})\|_{N}^{2}
 \mathbf 1_{\Lambda_{\bf m}\cap\Omega_{\bf m}})\\
 & &
 \hspace{0.5cm} +
 \mathbb E(\|(\widehat a_{m_1},\widehat b_{m_2}) - (a_{m_1}^{N},b_{m_2}^{N})\|_{N}^{2}
 \mathbf 1_{\Lambda_{\bf m}\cap\Omega_{\bf m}^{c}}) +
 \mathbb E(\|(a,b)\mathbf 1_A\|_{N}^{2}\mathbf 1_{\Lambda_{\bf m}^{c}})\\
 & =: &
 \min_{(\tau,\nu)\in\mathcal S_{\bf m}}
 \|(\tau,\nu) - (a,b)\mathbf 1_A\|_{f}^{2} + V + R_1 + R_2.
\end{eqnarray*}
Now, let us control the variance term $V$ and the remainder terms $R_1$ and $R_2$.
\begin{itemize}
 \item {\bf Control of $V$.} By Equality (\ref{nonadaptive_risk_bound_3}), $\widehat{\bf Z}_{\bf m} -\mathbf Z_{\bf m} =\mathbf E_{\bf m}$, leading to
 \begin{displaymath}
 (\widehat a_{m_1},\widehat b_{m_2}) - (a_{m_1}^{N},b_{m_2}^{N}) =
 \left(\sum_{j = 1}^{m_1}\delta_j\varphi_j,\sum_{k = 1}^{m_2}\delta_{m_1 + k}\psi_k\right)
 \quad {\rm with}\quad
 \delta =\widehat\Psi_{\bf m}^{-1}\mathbf E_{\bf m} -
 \frac{\mathbf d_{\bf m}^{*}\widehat\Psi_{\bf m}^{-1}\mathbf E_{\bf m}}{
 \mathbf d_{\bf m}^{*}\widehat\Psi_{\bf m}^{-1}\mathbf d_{\bf m}}
 \widehat\Psi_{\bf m}^{-1}\mathbf d_{\bf m},
 \end{displaymath}
 and then
 \begin{eqnarray}
  \|(\widehat a_{m_1},\widehat b_{m_2}) - (a_{m_1}^{N},b_{m_2}^{N})\|_{N}^{2}
  & = &
  \frac{1}{NT_0}\sum_{i = 1}^{N}\int_{t_0}^{T}
  \left(\sum_{j = 1}^{m_1}\delta_j\varphi_j(X_{s}^{i}) +
  \sum_{k = 1}^{m_2}\delta_{m_1 + k}\psi_k(Y_{s}^{i})\right)^2ds
  \nonumber\\
  & = &
  \delta^*\widehat\Psi_{\bf m}\delta
  \nonumber\\
  & = &
  \mathbf E_{\bf m}^{*}\widehat\Psi_{\bf m}^{-1}\mathbf E_{\bf m} -
  \frac{\mathbf d_{\bf m}^{*}\widehat\Psi_{\bf m}^{-1}\mathbf E_{\bf m}}{
  \mathbf d_{\bf m}^{*}\widehat\Psi_{\bf m}^{-1}\mathbf d_{\bf m}}
  (\mathbf d_{\bf m}^{*}\widehat\Psi_{\bf m}^{-1}\mathbf E_{\bf m} +
  \mathbf E_{\bf m}^{*}\widehat\Psi_{\bf m}^{-1}\mathbf d_{\bf m})
  \nonumber\\
  & &
  \hspace{4.25cm} +
  \left(\frac{\mathbf d_{\bf m}^{*}\widehat\Psi_{\bf m}^{-1}\mathbf E_{\bf m}}{
  \mathbf d_{\bf m}^{*}\widehat\Psi_{\bf m}^{-1}\mathbf d_{\bf m}}\right)^2
  \mathbf d_{\bf m}^{*}\widehat\Psi_{\bf m}^{-1}\mathbf d_{\bf m}
  \nonumber\\
  \label{nonadaptive_risk_bound_6}
  & = &
  \mathbf E_{\bf m}^{*}\widehat\Psi_{\bf m}^{-1}\mathbf E_{\bf m} -
  \frac{(\mathbf d_{\bf m}^{*}\widehat\Psi_{\bf m}^{-1}\mathbf E_{\bf m})^2}{
  \mathbf d_{\bf m}^{*}\widehat\Psi_{\bf m}^{-1}\mathbf d_{\bf m}}
  \leqslant
  \mathbf E_{\bf m}^{*}\widehat\Psi_{\bf m}^{-1}\mathbf E_{\bf m}.
 \end{eqnarray}
 Moreover, on $\Omega_{\bf m}$,
 \begin{displaymath}
 \mathbf E_{\bf m}^{*}\widehat\Psi_{\bf m}^{-1}\mathbf E_{\bf m}\leqslant
 2\mathbf E_{\bf m}^{*}\Psi_{\bf m}^{-1}\mathbf E_{\bf m}.
 \end{displaymath}
 Then,
 \begin{displaymath}
 V\leqslant 2\mathbb E(
 \mathbf E_{\bf m}^{*}\Psi_{\bf m}^{-1}\mathbf E_{\bf m}) =
 \frac{2}{NT_0}{\rm Tr}(\Psi_{\bf m}^{-1}\Psi_{\bf m}^{\sigma}),
 \end{displaymath}
 where
 \begin{displaymath}
 \Psi_{\bf m}^{\sigma} :=
 \begin{pmatrix}
  \Psi_{1,1}^{\sigma} & \Psi_{1,2}^{\sigma}\\
  \Psi_{1,2}^{\sigma,*} & \Psi_{2,2}^{\sigma}
 \end{pmatrix}
 \end{displaymath}
 with
 \begin{eqnarray*}
  \Psi_{1,1}^{\sigma} & = &
  \left(\int_{\mathbb R}
  \varphi_j(x)\varphi_{j'}(x)\sigma(x)^2f_X(x)dx\right)_{1\leqslant j,j'\leqslant m_1},\\
  \Psi_{2,2}^{\sigma} & = &
  \left(\int_{\mathbb R^2}
  \psi_k(y)\psi_{k'}(y)\sigma(x)^2f(x,y)dxdy\right)_{1\leqslant k,k'\leqslant m_2}\quad {\rm and}\\
  \Psi_{1,2}^{\sigma} & = &
  \left(\int_{\mathbb R^2}
  \varphi_j(x)\psi_k(y)\sigma(x)^2f(x,y)dxdy\right)_{(j,k)\in\{1,\dots,m_1\}\times\{1,\dots,m_2\}}.
 \end{eqnarray*}
 Therefore, by following the same line as in the proof of Comte and Genon-Catalot \cite{CGC20}, Proposition 2.2,
 \begin{displaymath}
 V\leqslant
 \frac{2\|\sigma\|_{\infty}^{2}}{T_0}\cdot\frac{m_1 + m_2}{N}.
 \end{displaymath}
 \item {\bf Control of $R_1$.} By Inequality (\ref{nonadaptive_risk_bound_6}), on $\Lambda_{\bf m}$,
 \begin{displaymath}
 \|(\widehat a_{m_1},\widehat b_{m_2}) - (a_{m_1}^{N},b_{m_2}^{N})\|_{N}^{2}
 \leqslant
 \|\widehat\Psi_{\bf m}^{-1}\|_{\rm op}\|\mathbf E_{\bf m}\|_{2,\mathbb R^{m_1 + m_2}}^{2}
 \leqslant
 \mathfrak c_rN\|\mathbf E_{\bf m}\|_{2,\mathbb R^{m_1 + m_2}}^{2}.
 \end{displaymath}
 Then, by Lemma \ref{2nd_order_moment_E},
 \begin{displaymath}
 R_1\leqslant
 \mathfrak c_rN\mathbb E(
 \|\mathbf E_{\bf m}\|_{2,\mathbb R^{m_1 + m_2}}^{2}\mathbf 1_{\Omega_{\bf m}^{c}})
 \leqslant
 \frac{\mathfrak c_r\mathfrak c_{\ref{2nd_order_moment_E}}}{N^{(r - 3)/2}}.
 \end{displaymath}
 \item {\bf Control of $R_2$.} Since $a + b\in\mathbb L^4(A,f(x,y)dxdy)$,
 \begin{displaymath}
 \mathbb E(\|(a,b)\mathbf 1_A\|_{N}^{4})
 \leqslant
 \int_A(a(x) + b(y))^4f(x,y)dxdy <\infty
 \end{displaymath}
 and then, by Inequality (\ref{nonadaptive_risk_bound_2}),
 \begin{displaymath}
 R_2\leqslant
 \mathbb E(\|(a,b)\mathbf 1_A\|_{N}^{4})^{\frac{1}{2}}\mathbb P(\Lambda_{\bf m}^{c})
 \leqslant
 \mathfrak c_{1}^{\frac{1}{2}}\|(a + b)^2\mathbf 1_A\|_f\frac{1}{N^{r/2}}.
 \end{displaymath}
\end{itemize}
Since $r\geqslant 5$, gathering these controls of $V$, $R_1$ and $R_2$ gives Inequality (\ref{nonadaptive_risk_bound_1}).
%


%
\subsubsection{Proof of Lemma \ref{2nd_order_moment_E}}\label{section_proof_2nd_order_moment_E}
First,
\begin{equation}\label{2nd_order_moment_E_1}
\mathbb E(\|\mathbf E_{\bf m}\|_{2,\mathbb R^{m_1 + m_2}}^{2}\mathbf 1_{\Omega_{\bf m}^{c}})
\leqslant
\mathbb E(\|\mathbf E_{\bf m}\|_{2,\mathbb R^{m_1 + m_2}}^{4})^{\frac{1}{2}}
\mathbb P(\Omega_{\bf m}^{c})^{\frac{1}{2}},
\end{equation}
and
\begin{eqnarray*}
 \mathbb E(\|\mathbf E_{\bf m}\|_{2,\mathbb R^{m_1 + m_2}}^{4})
 & = &
 \mathbb E\left[\left(
 \sum_{j = 1}^{m_1}\zeta_N(\varphi_j,0)^2 +
 \sum_{k = 1}^{m_2}\zeta_N(0,\psi_k)^2
 \right)^2\right]\\
 & = &
 \frac{1}{(NT_0)^4}\mathbb E\left[\left(
 \sum_{j = 1}^{m_1}M_j(T)^2 +
 \sum_{k = 1}^{m_2}N_k(T)^2
 \right)^2\right]\\
 & \leqslant &
 \frac{2}{(NT_0)^4}\left(
 m_1\sum_{j = 1}^{m_1}\mathbb E(M_j(T)^4) +
 m_2\sum_{k = 1}^{m_2}\mathbb E(N_k(T)^4)\right)
\end{eqnarray*}
where, for every $j\in\{1,\dots,m_1\}$ and $k\in\{1,\dots,m_2\}$, $M_j$ and $N_k$ are the $\mathbb F$-martingales defined by
\begin{eqnarray*}
 M_j(s) & := &
 \sum_{i = 1}^{N}\int_{t_0}^{s}
 \varphi_j(X_{u}^{i})\sigma(X_{u}^{i})dW_{1}^{i}(u)\\
 & &
 \hspace{1.5cm}
 \quad {\rm and}\quad
 N_k(s) :=
 \sum_{i = 1}^{N}\int_{t_0}^{s}
 \psi_k(Y_{u}^{i})\sigma(X_{u}^{i})dW_{1}^{i}(u)
 \textrm{ $;$ }\forall s\in [t_0,T].
\end{eqnarray*}
Now, by the Burkholder-Davis-Gundy and Jensen's inequalities, there exists a constant $\mathfrak c_1 > 0$, not depending on $\mathbf m$ and $N$, such that for every $j\in\{1,\dots,m_1\}$ and $k\in\{1,\dots,m_2\}$,
\begin{eqnarray*}
 \mathbb E(M_j(T)^4)
 & \leqslant &
 \mathfrak c_1\mathbb E\left[\left(\sum_{i = 1}^{N}\int_{t_0}^{T}
 \varphi_j(X_{u}^{i})^2\sigma(X_{u}^{i})^2du\right)^2\right]\\
 & \leqslant &
 \mathfrak c_1N^2T_0\int_{t_0}^{T}\mathbb E(\varphi_j(X_u)^4\sigma(X_u)^4)du
\end{eqnarray*}
and
\begin{eqnarray*}
 \mathbb E(N_k(T)^4)
 & \leqslant &
 \mathfrak c_1\mathbb E\left[\left(\sum_{i = 1}^{N}\int_{t_0}^{T}
 \psi_k(Y_{u}^{i})^2\sigma(X_{u}^{i})^2du\right)^2\right]\\
 & \leqslant &
 \mathfrak c_1N^2T_0\int_{t_0}^{T}\mathbb E(\psi_k(Y_u)^4\sigma(X_u)^4)du.
\end{eqnarray*}
Therefore,
\begin{eqnarray}
 \mathbb E(\|\mathbf E_{\bf m}\|_{2,\mathbb R^{m_1 + m_2}}^{4})
 & \leqslant &
 \frac{2\mathfrak c_1}{N^2T_{0}^{3}}\left(
 m_1\sum_{j = 1}^{m_1}\int_{t_0}^{T}\mathbb E(\varphi_j(X_u)^4\sigma(X_u)^4)du +
 m_2\sum_{k = 1}^{m_2}\int_{t_0}^{T}\mathbb E(\psi_k(Y_u)^4\sigma(X_u)^4)du\right)
 \nonumber\\
 & \leqslant &
 \frac{2\mathfrak c_1}{N^2T_{0}^{2}}
 (m_1\mathfrak L_{\varphi}(m_1)^2 + m_2\mathfrak L_{\psi}(m_2)^2)
 \nonumber\\
 \label{2nd_order_moment_E_2}
 & &
 \hspace{3.5cm}
 \times\int_{-\infty}^{\infty}\sigma(x)^4f_X(x)dx
 \quad {\rm with}\quad
 f_X(.) =\int_{-\infty}^{\infty}f(.,y)dy.
\end{eqnarray}
The conclusion follows from Inequalities (\ref{2nd_order_moment_E_1}) and (\ref{2nd_order_moment_E_2}), Assumption \ref{stability_condition} and Lemma \ref{bound_Omega}.
%


%
\subsection{Proof of Proposition \ref{nonadaptive_risk_bound_f_weighted}}\label{section_proof_nonadaptive_risk_bound_f_weighted}
First,
\begin{eqnarray}
 \label{nonadaptive_risk_bound_f_weighted_2}
 \mathbb E(\|(\widetilde a_{m_1},\widetilde b_{m_2}) - (a,b)\mathbf 1_A\|_{f}^{2})
 & = &
 \mathbb E(\|(\widehat a_{m_1},\widehat b_{m_2}) - (a,b)\mathbf 1_A\|_{f}^{2}
 \mathbf 1_{\Lambda_{\bf m}}) +
 \mathbb E(\|(a,b)\mathbf 1_A\|_{f}^{2}\mathbf 1_{\Lambda_{\bf m}^{c}})\\
 & =: &
 \mathbb T +\mathbb S.
 \nonumber
\end{eqnarray}
Since $a + b\in\mathbb L^4(A,f(x,y)dxdy)$, and by Inequality (\ref{nonadaptive_risk_bound_2}), there exists a constant $\mathfrak c_1 > 0$, not depending on $\mathbf m$ and $N$, such that
\begin{displaymath}
\mathbb S =
\|(a,b)\mathbf 1_A\|_{f}^{2}\mathbb P(\Lambda_{\bf m}^{c})\leqslant
\frac{\mathfrak c_1}{N}.
\end{displaymath}
Let $(a_{m_1}^{f},b_{m_2}^{f})$ be the minimizer of $(u,v)\mapsto\|(u,v) - (a,b)\mathbf 1_A\|_{f}^{2}$ over $\mathcal S_{\bf m}$. By following the same line as in the proof of Proposition \ref{nonadaptive_risk_bound} (see the beginning of Step 2),
\begin{displaymath}
\mathbb T =
\min_{(\tau,\nu)\in\mathcal S_{\bf m}}
\|(\tau,\nu) - (a,b)\mathbf 1_A\|_{f}^{2} +
\mathbb E(\|(\widehat a_{m_1},\widehat b_{m_2}) -
(a_{m_1}^{f},b_{m_2}^{f})\|_{f}^{2}\mathbf 1_{\Lambda_{\bf m}}).
\end{displaymath}
Now, note that
\begin{eqnarray*}
 \mathbb E(\|(\widehat a_{m_1},\widehat b_{m_2}) -
 (a_{m_1}^{f},b_{m_2}^{f})\|_{f}^{2}\mathbf 1_{\Lambda_{\bf m}})
 & = &
 \mathbb E(\|(\widehat a_{m_1},\widehat b_{m_2}) -
 (a_{m_1}^{f},b_{m_2}^{f})\|_{f}^{2}\mathbf 1_{\Omega_{\bf m}\cap\Lambda_{\bf m}})\\
 & &
 \hspace{2cm} +
 \mathbb E(\|(\widehat a_{m_1},\widehat b_{m_2}) -
 (a_{m_1}^{f},b_{m_2}^{f})\|_{f}^{2}\mathbf 1_{\Omega_{\bf m}^{c}\cap\Lambda_{\bf m}})\\
 & =: &
 \mathbb T_1 +\mathbb T_2,
\end{eqnarray*}
and let us provide suitable controls of $\mathbb T_1$ and $\mathbb T_2$.
\begin{itemize}
 \item {\bf Control of $\mathbb T_1$.} By the definition of $\Omega_{\bf m}$, and by Proposition \ref{nonadaptive_risk_bound},
 \begin{eqnarray*}
  \mathbb T_1 & \leqslant &
  2\mathbb E(\|(\widehat a_{m_1},\widehat b_{m_2}) -
  (a_{m_1}^{f},b_{m_2}^{f})\|_{N}^{2}\mathbf 1_{\Omega_{\bf m}\cap\Lambda_{\bf m}})\\
  & \leqslant &
  4\left(\min_{(\tau,\nu)\in\mathcal S_{\bf m}}\|(\tau,\nu) - (a,b)\mathbf 1_A\|_{f}^{2} +
  \frac{2\|\sigma\|_{\infty}^{2}}{T_0}\cdot\frac{m_1 + m_2}{N} +
  \frac{\mathfrak c_{\ref{nonadaptive_risk_bound}}}{N}\right) +
  4\mathbb E(\|(a_{m_1}^{f},b_{m_2}^{f}) - (a,b)\mathbf 1_A\|_{N}^{2}).
 \end{eqnarray*}
 Moreover,
 \begin{displaymath}
 \mathbb E(\|(a_{m_1}^{f},b_{m_2}^{f}) - (a,b)\mathbf 1_A\|_{N}^{2}) =
 \|(a_{m_1}^{f},b_{m_2}^{f}) - (a,b)\mathbf 1_A\|_{f}^{2},
 \end{displaymath}
 leading to
 \begin{equation}\label{nonadaptive_risk_bound_f_weighted_3}
 \mathbb T_1\leqslant
 8\min_{(\tau,\nu)\in\mathcal S_{\bf m}}\|(\tau,\nu) - (a,b)\mathbf 1_A\|_{f}^{2} +
 \frac{8\|\sigma\|_{\infty}^{2}}{T_0}\cdot\frac{m_1 + m_2}{N} +
 \frac{4\mathfrak c_{\ref{nonadaptive_risk_bound}}}{N}.
 \end{equation}
 \item {\bf Control of $\mathbb T_2$.} By Lemma \ref{bound_Omega},
 \begin{displaymath}
 \|(a_{m_1},b_{m_2})\|_{f}^{2}\mathbb P(\Omega_{\bf m}^{c})
 \leqslant
 2\min_{(\tau,\nu)\in\mathcal S_{\bf m}}\|(\tau,\nu) - (a,b)\mathbf 1_A\|_{f}^{2} +
 2\mathfrak c_{\ref{bound_Omega}}\|(a,b)\mathbf 1_A\|_{f}^{2}\frac{1}{N},
 \end{displaymath}
 leading to
 \begin{eqnarray*}
  \mathbb T_2
  & \leqslant &
  2\mathbb E(\|(\widehat a_{m_1},\widehat b_{m_2})\|_{f}^{2}
  \mathbf 1_{\Omega_{\bf m}^{c}\cap\Lambda_{\bf m}}) +
  2\|(a_{m_1},b_{m_2})\|_{f}^{2}\mathbb P(\Omega_{\bf m}^{c})\\
  & \leqslant &
  4\min_{(\tau,\nu)\in\mathcal S_{\bf m}}\|(\tau,\nu) - (a,b)\mathbf 1_A\|_{f}^{2} +
  2\mathbb E(\|(\widehat a_{m_1},\widehat b_{m_2})\|_{f}^{4}
  \mathbf 1_{\Lambda_{\bf m}})^{\frac{1}{2}}\mathbb P(\Omega_{\bf m}^{c})^{\frac{1}{2}} +
  4\mathfrak c_{\ref{bound_Omega}}\|(a,b)\mathbf 1_A\|_{f}^{2}\frac{1}{N}.
 \end{eqnarray*}
 Moreover,
 \begin{eqnarray*}
  \|(\widehat a_{m_1},\widehat b_{m_2})\|_{f}^{2}
  & = &
  \widehat\theta^*\Psi_{\bf m}\widehat\theta =
  \left(\widehat{\bf Z}_{\bf m} -
  \frac{\mathbf d_{\bf m}^{*}\widehat\Psi_{\bf m}^{-1}\widehat{\bf Z}_{\bf m}}{
  \mathbf d_{\bf m}^{*}\widehat\Psi_{\bf m}^{-1}\mathbf d_{\bf m}}\cdot
  \mathbf d_{\bf m}\right)^*
  \widehat\Psi_{\bf m}^{-1}\Psi_{\bf m}\widehat\Psi_{\bf m}^{-1}
  \left(\widehat{\bf Z}_{\bf m} -
  \frac{\mathbf d_{\bf m}^{*}\widehat\Psi_{\bf m}^{-1}\widehat{\bf Z}_{\bf m}}{
  \mathbf d_{\bf m}^{*}\widehat\Psi_{\bf m}^{-1}\mathbf d_{\bf m}}\cdot
  \mathbf d_{\bf m}\right)\\
  & = &
  \left\|\Psi_{\bf m}^{\frac{1}{2}}\widehat\Psi_{\bf m}^{-1}
  \left(\widehat{\bf Z}_{\bf m} -
  \frac{\mathbf d_{\bf m}^{*}\widehat\Psi_{\bf m}^{-1}\widehat{\bf Z}_{\bf m}}{
  \mathbf d_{\bf m}^{*}\widehat\Psi_{\bf m}^{-1}\mathbf d_{\bf m}}\cdot
  \mathbf d_{\bf m}\right)\right\|_{2,\mathbb R^{m_1 + m_2}}^{2}\\
  & \leqslant &
  2\|\Psi_{\bf m}^{\frac{1}{2}}
  \widehat\Psi_{\bf m}^{-1}\widehat{\bf Z}_{\bf m}\|_{2,\mathbb R^{m_1 + m_2}}^{2} +
  2\left(\frac{\mathbf d_{\bf m}^{*}\widehat\Psi_{\bf m}^{-1}\widehat{\bf Z}_{\bf m}}{
  \mathbf d_{\bf m}^{*}\widehat\Psi_{\bf m}^{-1}\mathbf d_{\bf m}}\right)^2
  \|\Psi_{\bf m}^{\frac{1}{2}}\widehat\Psi_{\bf m}^{-1}
  \mathbf d_{\bf m}\|_{2,\mathbb R^{m_1 + m_2}}^{2}.
 \end{eqnarray*}
 By the Cauchy-Schwarz inequality,
 \begin{eqnarray*}
  (\mathbf d_{\bf m}^{*}\widehat\Psi_{\bf m}^{-1}\widehat{\bf Z}_{\bf m})^2
  & = &
  \langle\widehat\Psi_{\bf m}^{-\frac{1}{2}}\mathbf d_{\bf m},
  \widehat\Psi_{\bf m}^{-\frac{1}{2}}\widehat{\bf Z}_{\bf m}\rangle_{2,\mathbb R^{m_1 + m_2}}^{2}\\
  & \leqslant &
  \|\widehat\Psi_{\bf m}^{-\frac{1}{2}}\mathbf d_{\bf m}\|_{2,\mathbb R^{m_1 + m_2}}^{2}
  \|\widehat\Psi_{\bf m}^{-\frac{1}{2}}\widehat{\bf Z}_{\bf m}\|_{2,\mathbb R^{m_1 + m_2}}^{2} =
  \mathbf d_{\bf m}^{*}\widehat\Psi_{\bf m}^{-1}\mathbf d_{\bf m}
  \widehat{\bf Z}_{\bf m}^{*}\widehat\Psi_{\bf m}^{-1}\widehat{\bf Z}_{\bf m},
 \end{eqnarray*}
 leading to
 \begin{eqnarray*}
  \left(\frac{\mathbf d_{\bf m}^{*}\widehat\Psi_{\bf m}^{-1}\widehat{\bf Z}_{\bf m}}{
  \mathbf d_{\bf m}^{*}\widehat\Psi_{\bf m}^{-1}\mathbf d_{\bf m}}\right)^2
  \|\Psi_{\bf m}^{\frac{1}{2}}\widehat\Psi_{\bf m}^{-1}
  \mathbf d_{\bf m}\|_{2,\mathbb R^{m_1 + m_2}}^{2}
  & \leqslant &
  \frac{\widehat{\bf Z}_{\bf m}^{*}\widehat\Psi_{\bf m}^{-1}
  \widehat{\bf Z}_{\bf m}}{\mathbf d_{\bf m}^{*}\widehat\Psi_{\bf m}^{-1}\mathbf d_{\bf m}}
  \|\Psi_{\bf m}\|_{\rm op}
  \|\widehat\Psi_{\bf m}^{-1}\mathbf d_{\bf m}\|_{2,\mathbb R^{m_1 + m_2}}^{2}\\
  & = &
  \|\Psi_{\bf m}\|_{\rm op}
  \|\widehat\Psi_{\bf m}^{-1}\widehat{\bf Z}_{\bf m}\|_{2,\mathbb R^{m_1 + m_2}}^{2}
  \leqslant
  \|\Psi_{\bf m}\|_{\rm op}
  \|\widehat\Psi_{\bf m}^{-1}\|_{\rm op}^{2}
  \|\widehat{\bf Z}_{\bf m}\|_{2,\mathbb R^{m_1 + m_2}}^{2}.
 \end{eqnarray*}
 Thus, by the definition of $\Lambda_{\bf m}$, and since $\|\Psi_{\bf m}\|_{\rm op}\leqslant 2(\mathfrak L_{\varphi}(m_1) +\mathfrak L_{\psi}(m_2))$,
 \begin{eqnarray*}
  \|(\widehat a_{m_1},\widehat b_{m_2})\|_{f}^{2}\mathbf 1_{\Lambda_{\bf m}}
  & \leqslant &
  4\|\Psi_{\bf m}\|_{\rm op}
  \|\widehat\Psi_{\bf m}^{-1}\|_{\rm op}^{2}
  \|\widehat{\bf Z}_{\bf m}\|_{2,\mathbb R^{m_1 + m_2}}^{2}\\
  & \leqslant &
  8\mathfrak c_{r}^{2}\frac{N^2}{(\mathfrak L_{\varphi}(m_1) +
  \mathfrak L_{\psi}(m_2))\log(N)^2}\|\widehat{\bf Z}_{\bf m}\|_{2,\mathbb R^{m_1 + m_2}}^{2}.
 \end{eqnarray*}
 By Inequality (\ref{2nd_order_moment_E_2}), there exists a constant $\mathfrak c_2 > 0$, not depending on $\mathbf m$ and $N$, such that
 \begin{eqnarray*}
  \mathbb E(\|\widehat{\bf Z}_{\bf m}\|_{2,\mathbb R^{m_1 + m_2}}^{4})
  & \leqslant &
  8(\mathbb E(\|\mathbf Z_{\bf m}\|_{2,\mathbb R^{m_1 + m_2}}^{4}) +
  \mathbb E(\|\mathbf E_{\bf m}\|_{2,\mathbb R^{m_1 + m_2}}^{4}))\\
  & \leqslant &
  \mathfrak c_2N(\mathfrak L_{\varphi}(m_1)^2 +\mathfrak L_{\psi}(m_2)^2)
  \left(\int_A(a(x) + b(y))^4f(x,y)dxdy +
  \int_{-\infty}^{\infty}\sigma(x)^4f_X(x)dx\right).
 \end{eqnarray*}
 So, by Lemma \ref{bound_Omega}, and since $r\geqslant 7$, there exists a constant $\mathfrak c_3 > 0$, not depending on $\mathbf m$ and $N$, such that
 \begin{displaymath}
 \mathbb E(\|(\widehat a_{m_1},\widehat b_{m_2})\|_{f}^{4}
 \mathbf 1_{\Lambda_{\bf m}})^{\frac{1}{2}}\mathbb P(\Omega_{\bf m}^{c})^{\frac{1}{2}}
 \leqslant
 \frac{\mathfrak c_3}{N}.
 \end{displaymath}
 Therefore,
 \begin{equation}\label{nonadaptive_risk_bound_f_weighted_4}
 \mathbb T_2\leqslant
 4\min_{(\tau,\nu)\in\mathcal S_{\bf m}}\|(\tau,\nu) - (a,b)\mathbf 1_A\|_{f}^{2} +
 2(\mathfrak c_3 + 2\mathfrak c_{\ref{bound_Omega}}\|(a,b)\mathbf 1_A\|_{f}^{2})\frac{1}{N}.
 \end{equation}
\end{itemize}
In conclusion, by Inequalities (\ref{nonadaptive_risk_bound_f_weighted_2}), (\ref{nonadaptive_risk_bound_f_weighted_3}) and (\ref{nonadaptive_risk_bound_f_weighted_4}), there exists a constant $\mathfrak c_4 > 0$, not depending on $\mathbf m$ and $N$, such that
\begin{displaymath}
\mathbb E(\|(\widetilde a_{m_1},\widetilde b_{m_2}) - (a,b)\mathbf 1_A\|_{f}^{2})
\leqslant
13\min_{(\tau,\nu)\in\mathcal S_{\bf m}}\|(\tau,\nu) - (a,b)\mathbf 1_A\|_{f}^{2} +
\frac{8\|\sigma\|_{\infty}^{2}}{T_0}\cdot\frac{m_1 + m_2}{N} +
\frac{\mathfrak c_4}{N}.
\end{displaymath}
%


%
\subsection{Proof of Lemma \ref{control_bias_term_non_degenerate_constraint}}\label{section_proof_control_bias_term_non_degenerate_constraint}
Recall that
\begin{displaymath}
b_{m_2} =\sum_{k = 1}^{m_2}\langle b,\psi_k\rangle\psi_k,
\end{displaymath}
and let $\overline b_{m_2}$ be the orthogonal projection of $b$ on $\mathbb S_{m_2}$ for the usual inner product in $\mathbb L^2(A_2)$. Precisely,
\begin{displaymath}
\overline b_{m_2} =\sum_{k = 1}^{m_2}\overline n_k\psi_k
\end{displaymath}
with
\begin{equation}\label{control_bias_term_non_degenerate_constraint_2}
\overline n =
\underset{n\in\mathbb R^{m_2} :\eta(n) = 0}{\rm argmin}J_{m_2}(n),
\end{equation}
where
\begin{displaymath}
J_{m_2}(n) :=\int_{A_2}\left(\sum_{k = 1}^{m_2}n_k\psi_k(y) - b(y)\right)^2dy
\quad {\rm and}\quad
\eta(n) :=\langle n,\delta_{m_2}\rangle_{2,\mathbb R^{m_2}}.
\end{displaymath}
Consider
\begin{displaymath}
\mathbf z_{m_2} :=
(\langle b,\psi_1\rangle,\dots,\langle b,\psi_{m_2}\rangle),
\end{displaymath}
and let $L_{m_2}$ be the Lagrangian for Problem (\ref{control_bias_term_non_degenerate_constraint_2}):
\begin{displaymath}
L_{m_2}(n,\lambda) := J_{m_2}(n) -\lambda\eta(n)
\textrm{ $;$ }(n,\lambda)\in\mathbb R^{m_2}\times\mathbb R.
\end{displaymath}
Necessarily,
\begin{displaymath}
\nabla L_{m_2}(\overline n,\overline\lambda) =
\begin{pmatrix}
 2(\overline n -\mathbf z_{m_2}) -\overline\lambda\delta_{m_2}\\
 -\eta(\overline n)
\end{pmatrix} = 0,
\end{displaymath}
leading to
\begin{displaymath}
\overline n =\mathbf z_{m_2} +\frac{\overline\lambda}{2}\delta_{m_2},
\quad\textrm{and then}\quad
\overline\lambda = -2\cdot\frac{\langle\mathbf z_{m_2},
\delta_{m_2}\rangle_{2,\mathbb R^{m_2}}}{\|\delta_{m_2}\|_{2,\mathbb R^{m_2}}^{2}}.
\end{displaymath}
So,
\begin{displaymath}
\overline n =\mathbf z_{m_2} -\frac{\langle\mathbf z_{m_2},
\delta_{m_2}\rangle_{2,\mathbb R^{m_2}}}{\|\delta_{m_2}\|_{2,\mathbb R^{m_2}}^{2}}\delta_{m_2}
\end{displaymath}
and, by the condition (\ref{identifiability_condition}) on $b$,
\begin{eqnarray*}
 \|\overline b_{m_2} - b_{m_2}\|^2 & = &
 \left\|\frac{\langle\mathbf z_{m_2},
 \delta_{m_2}\rangle_{2,\mathbb R^{m_2}}}{\|\delta_{m_2}\|_{2,\mathbb R^{m_2}}^{2}}
 \sum_{k = 1}^{m_2}\delta_{m_2,k}\psi_k\right\|^2\\
 & = &
 \frac{\langle\mathbf z_{m_2},
 \delta_{m_2}\rangle_{2,\mathbb R^{m_2}}^{2}}{\|\delta_{m_2}\|_{2,\mathbb R^{m_2}}^{2}} =
 \frac{1}{\|\delta_{m_2}\|_{2,\mathbb R^{m_2}}^{2}}\left(
 \int_{A_2}(b_{m_2}(y) - b(y))dy\right)^2.
\end{eqnarray*}
Since $b_{m_2}$ is the orthogonal projection of $b$ on $\Sigma_{m_2}$, and since $\overline b_{m_2} - b_{m_2}\in\Sigma_{m_2}$,
\begin{displaymath}
\min_{\nu\in\mathbb S_{m_2}}\|\nu - b\mathbf 1_{A_2}\|^2 =
\|b_{m_2} - b\mathbf 1_{A_2}\|^2 +
\frac{1}{\|\delta_{m_2}\|_{2,\mathbb R^{m_2}}^{2}}\left(
\int_{A_2}(b_{m_2}(y) - b(y))dy\right)^2.
\end{displaymath}
Then, Inequality (\ref{bias_term_bound}) leads to (\ref{control_bias_term_non_degenerate_constraint_1}), which ends the proof of Lemma \ref{control_bias_term_non_degenerate_constraint}.
%


%
\subsection{Proof of Proposition \ref{rate_Sobolev}}\label{section_proof_rate_Sobolev}
By Proposition \ref{nonadaptive_risk_bound} together with Lemma \ref{control_bias_term_non_degenerate_constraint},
\begin{eqnarray}
 & &
 \mathbb E(\|(\widetilde a_{m_1},\widetilde b_{m_2}) - (a,b)\mathbf 1_A\|_{N}^{2})
 \nonumber\\
 \label{rate_Sobolev_3}
 & &
 \hspace{2cm}\leqslant
 2\mathfrak c_{f,1}\|a_{m_1} - a\mathbf 1_{A_1}\|^2 +
 2\mathfrak c_{f,2}\|b_{m_2} - b\mathbf 1_{A_2}\|^2 + R(m_2) +
 \frac{2\|\sigma\|_{\infty}^{2}}{T_0}\cdot\frac{m_1 + m_2}{N} +
 \frac{\mathfrak c_{\ref{nonadaptive_risk_bound}}}{N},
\end{eqnarray}
where
\begin{displaymath}
R(m_2) :=
\frac{2\mathfrak c_{f,2}}{\|\delta_{m_2}\|_{2,\mathbb R^{m_2}}^{2}}\left(
\int_{A_2}(b_{m_2}(y) - b(y))dy\right)^2.
\end{displaymath}
First, since $a\in\mathbb W_{\varphi}^{\alpha}(A_1,L_1)$ and $b\in\mathbb W_{\psi}^{\beta}(A_2,L_2)$,
\begin{equation}\label{rate_Sobolev_4}
\|a_{m_1} - a\mathbf 1_{A_1}\|^2\leqslant L_1m_{1}^{-\alpha}
\quad {\rm and}\quad
\|b_{m_2} - b\mathbf 1_{A_2}\|^2\leqslant L_2m_{2}^{-\beta}.
\end{equation}
Now, by the conditions (\ref{rate_Sobolev_1}) and (\ref{rate_Sobolev_2}),
\begin{eqnarray*}
 R(m_2) & \leqslant &
 \mathfrak c_1m_{2}^{-(\omega + 1)}\left(
 \sum_{k > m_2}\langle b,\psi_k\rangle\int_{A_2}\psi_k(y)dy\right)^2
 \quad {\rm with}\quad
 \mathfrak c_1 =\frac{2\mathfrak c_{f,2}}{\mathfrak c_{\psi,1}}\\
 & \leqslant &
 \mathfrak c_1m_{2}^{-(\omega + 1)}
 \left(\sum_{k > m_2}k^{\beta}\langle b,\psi_k\rangle^2\right)
 \left(\sum_{k > m_2}k^{-\beta}\left|\int_{A_2}\psi_k(y)dy\right|^2\right)
 \leqslant
 \mathfrak c_2
 m_{2}^{-(\omega + 1)}
 \sum_{k > m_2}k^{\omega -\beta}
\end{eqnarray*}
with $\mathfrak c_2 =\mathfrak c_1L_2\mathfrak c_{\psi,2}$. Moreover,
\begin{displaymath}
\sum_{k > m_2}k^{\omega -\beta}
\leqslant
\int_{m_2}^{\infty}y^{\omega -\beta}dy =
\frac{m_{2}^{1 +\omega -\beta}}{\beta - (1 +\omega)},
\end{displaymath}
leading to
\begin{equation}\label{rate_Sobolev_5}
R(m_2)\leqslant
\frac{\mathfrak c_2}{\beta - (1 +\omega)}m_{2}^{-\beta}.
\end{equation}
Therefore, by plugging (\ref{rate_Sobolev_4}) and (\ref{rate_Sobolev_5}) in Inequality (\ref{rate_Sobolev_3}),
\begin{displaymath}
\mathbb E(\|(\widetilde a_{m_1},\widetilde b_{m_2}) - (a,b)\mathbf 1_A\|_{N}^{2})
\lesssim
m_{1}^{-\alpha} + m_{2}^{-\beta} +
\frac{m_1 + m_2}{N}.
\end{displaymath}
The conclusion comes by taking $m_1 = m_{1}^{\star}\asymp N^{-1/(\alpha + 1)}$ and $m_2 = m_{2}^{\star}\asymp N^{-1/(\beta + 1)}$.
%


%
\subsection{Proof of Theorem \ref{adaptive_risk_bound}}\label{section_proof_adaptive_risk_bound}
Consider
\begin{displaymath}
\mathcal M_{N}^{+} :=
\left\{
\mathbf m = (m_1,m_2)\in\{1,\dots,N\}^2 :
(\mathfrak L_{\varphi}(m_1) +\mathfrak L_{\psi}(m_2))
(\|\Psi_{\bf m}^{-1}\|_{\rm op}\vee 1)\leqslant
\frac{3\mathfrak c_r}{2}\cdot\frac{N}{\log(N)}\right\}
\end{displaymath}
and
\begin{displaymath}
\Omega_N :=
\bigcap_{\mathbf m\in\mathcal M_{N}^{+}}\Omega_{\bf m}.
\end{displaymath}
The proof of Theorem \ref{adaptive_risk_bound} relies on the following lemma.
%


%
\begin{lemma}\label{relationship_collections_models}
Under the assumptions of Theorem \ref{adaptive_risk_bound},
\begin{displaymath}
\Omega_N\subset\Xi_N :=
\{\mathcal M_N\subset\widehat{\mathcal M}_N\subset\mathcal M_{N}^{+}\}.
\end{displaymath}
\end{lemma}
\noindent
The proof of Lemma \ref{relationship_collections_models} is postponed to Section \ref{section_proof_relationship_collections_models}.
%


%
\subsubsection{Steps of the proof}
Let $(a_{\widehat m_1}^{N},b_{\widehat m_2}^{N})$ be the minimizer of $(u,v)\mapsto \|(u,v) - (a,b)\mathbf 1_A\|_{N}^{2}$ over $\mathcal S_{\widehat{\bf m}}$. Then,
\begin{displaymath}
\|(\widehat a_{\widehat m_1},\widehat b_{\widehat m_2}) - (a,b)\mathbf 1_A\|_{N}^{2} =
\min_{(\tau,\nu)\in\mathcal S_{\widehat{\bf m}}}\|(\tau,\nu) - (a,b)\mathbf 1_A\|_{N}^{2} +
\|(\widehat a_{\widehat m_1},\widehat b_{\widehat m_2}) - (a_{\widehat m_1}^{N},b_{\widehat m_2}^{N})\|_{N}^{2}.
\end{displaymath}
Moreover, since $0\in\mathcal S_{\widehat{\bf m}}$,
\begin{displaymath}
\min_{(\tau,\nu)\in\mathcal S_{\widehat{\bf m}}}\|(\tau,\nu) - (a,b)\mathbf 1_A\|_{N}^{2}\leqslant
\|(a,b)\mathbf 1_A\|_{N}^{2}.
\end{displaymath}
So,
\begin{equation}\label{adaptive_risk_bound_1}
\|(\widehat a_{\widehat m_1},\widehat b_{\widehat m_2}) - (a,b)\mathbf 1_A\|_{N}^{2}\leqslant
\|(\widehat a_{\widehat m_1},\widehat b_{\widehat m_2}) - (a,b)\mathbf 1_A\|_{N}^{2}\mathbf 1_{\Omega_N} +
R(\widehat{\bf m})\mathbf 1_{\Omega_{N}^{c}},
\end{equation}
where
\begin{displaymath}
R(\widehat{\bf m}) :=\|(a,b)\mathbf 1_A\|_{N}^{2} +
\|(\widehat a_{\widehat m_1},\widehat b_{\widehat m_2}) - (a_{\widehat m_1}^{N},b_{\widehat m_2}^{N})\|_{N}^{2}.
\end{displaymath}
The proof of Theorem \ref{adaptive_risk_bound} is dissected in three steps. The first step provides a preliminary risk bound on our adaptive projection least squares estimator, which is improved in Step 3 thanks to the bound established in Step 2 on
\begin{displaymath}
\rho(\mathbf m) :=
\mathbb E\left(\left(\left[\sup_{(\tau,\nu)\in\mathcal B_{\mathbf m,\widehat{\bf m}}}|\zeta_N(\tau,\nu)|\right]^2 -
p(\mathbf m,\widehat{\bf m})\right)_+\mathbf 1_{\Omega_N}\right)\textrm{ $;$ }
\mathbf m = (m_1,m_2)\in\mathcal M_N
\end{displaymath}
where, for every $\mathbf m' = (m_1',m_2')\in\mathcal M_N$,
\begin{displaymath}
\mathcal B_{\mathbf m,\mathbf m'} :=
\{(\tau,\nu)\in\mathcal S_{(m_1\vee m_1',m_2\vee m_2')} :\|(\tau,\nu)\|_f = 1\}
\end{displaymath}
and
\begin{displaymath}
p(\mathbf m,\widehat{\bf m}) := p(m_1\vee m_1') + p(m_2\vee m_2')
\quad {\rm with}\quad
p(m) :=
\frac{\kappa\|\sigma\|_{\infty}^{2}}{8T_0}\cdot\frac{m}{N}.
\end{displaymath}
{\bf Step 1.} On the one hand, by Inequality (\ref{nonadaptive_risk_bound_6}), and by the definition of $\widehat{\mathcal M}_N$,
\begin{eqnarray*}
 \|(\widehat a_{\widehat m_1},\widehat b_{\widehat m_2}) - (a_{\widehat m_1}^{N},b_{\widehat m_2}^{N})\|_{N}^{2}
 & \leqslant &
 \|\widehat\Psi_{\widehat{\bf m}}^{-1}\|_{\rm op}
 \|\mathbf E_{\widehat{\bf m}}\|_{2,\mathbb R^{\widehat m_1 +\widehat m_2}}^{2}\\
 & \leqslant &
 \mathfrak c_rN\|\mathbf E_{\bf N}\|_{2,\mathbb R^{2N}}^{2}
 \quad {\rm with}\quad
 \mathbf N = (N,N).
\end{eqnarray*}
Moreover, by Lemma \ref{bound_Omega},
\begin{displaymath}
\mathbb P(\Omega_{N}^{c})\leqslant
\sum_{\mathbf m\in\mathcal M_{N}^{+}}
\mathbb P(\Omega_{\bf m}^{c})\leqslant
\frac{\mathfrak c_{\ref{bound_Omega}}}{N^{r - 2}}.
\end{displaymath}
Then, by Inequality (\ref{2nd_order_moment_E_2}),
\begin{eqnarray*}
 \mathbb E(R(\widehat{\bf m})\mathbf 1_{\Omega_{N}^{c}})
 & \leqslant &
 (\mathbb E(\|(a,b)\|_{N}^{4})^{\frac{1}{2}} +
 \mathfrak c_rN\mathbb E(\|\mathbf E_{\bf N}^{4}\|_{2,\mathbb R^{2N}}^{4})^{\frac{1}{2}})
 \mathbb P(\Omega_{N}^{c})^{\frac{1}{2}}\\
 & \leqslant &
 \mathfrak c_1(1 + N^{\frac{3}{2}})N^{1 -\frac{r}{2}}
 \leqslant 2\mathfrak c_1N^{-\frac{r - 5}{2}},
\end{eqnarray*}
where $\mathfrak c_1$ is a positive constant not depending on $N$. On the other hand, for every $(\tau,\nu),(\overline\tau,\overline\nu)\in\mathcal S_{\bf N}$,
\begin{displaymath}
\gamma_N(\overline\tau,\overline\nu) -\gamma_N(\tau,\nu) =
\|(\overline\tau,\overline\nu) - (a,b)\mathbf 1_A\|_{N}^{2} -
\|(\tau,\nu) - (a,b)\mathbf 1_A\|_{N}^{2} - 2\zeta_N(\overline\tau -\tau,\overline\nu -\nu).
\end{displaymath}
Moreover, by the definition of $\widehat{\bf m}$, for every $\mathbf m = (m_1,m_2)\in\widehat{\mathcal M}_N$,
\begin{equation}\label{adaptive_risk_bound_2}
\gamma_N(\widehat a_{\widehat m_1},\widehat b_{\widehat m_2}) +
{\rm pen}(\widehat{\bf m})\leqslant
\gamma_N(\widehat a_{m_1},\widehat b_{m_2}) +
{\rm pen}(\mathbf m).
\end{equation}
On the event $\Xi_N$, Inequality (\ref{adaptive_risk_bound_2}) remains true for every $\mathbf m\in\mathcal M_N$. Then, on the event $\Omega_N$ (which is contained in $\Xi_N$ by Lemma \ref{relationship_collections_models}), for any $\mathbf m = (m_1,m_2)\in\mathcal M_N$, since $\mathcal S_{\bf m} +\mathcal S_{\widehat{\bf m}}\subset\mathcal S_{(m_1\vee\widehat m_1,m_2\vee\widehat m_2)}$ by (\ref{consequence_nested}),
\begin{eqnarray*}
 \|(\widehat a_{\widehat m_1},\widehat b_{\widehat m_2}) - (a,b)\mathbf 1_A\|_{N}^{2}
 & \leqslant &
 \|(\widehat a_{m_1},\widehat b_{m_2}) - (a,b)\mathbf 1_A\|_{N}^{2}\\
 & &
 \hspace{1cm} +
 2\zeta_N(\widehat a_{\widehat m_1} -\widehat a_{m_1},\widehat b_{\widehat m_2} -\widehat b_{m_2}) +
 {\rm pen}(\mathbf m) - {\rm pen}(\widehat{\bf m})\\
 & \leqslant &
 \|(\widehat a_{m_1},\widehat b_{m_2}) - (a,b)\mathbf 1_A\|_{N}^{2} +
 \frac{1}{8}\|(\widehat a_{\widehat m_1},\widehat b_{\widehat m_2}) -
 (\widehat a_{m_1},\widehat b_{m_2})\|_{f}^{2}\\
 & &
 \hspace{1cm} +
 8\mathcal Z(\mathbf m,\widehat{\bf m}) +
 {\rm pen}(\mathbf m) + 8p(\mathbf m,\widehat{\bf m}) - {\rm pen}(\widehat{\bf m}),
\end{eqnarray*}
where
\begin{displaymath}
\mathcal Z(\mathbf m,\widehat{\bf m}) :=
\left(\left[\sup_{(\tau,\nu)\in\mathcal B_{\mathbf m,\widehat{\bf m}}}|\zeta_N(\tau,\nu)|\right]^2 -
p(\mathbf m,\widehat{\bf m})\right)_+.
\end{displaymath}
Since $\|(\tau,\nu)\|_{f}^{2}\mathbf 1_{\Omega_N}\leqslant 2\|(\tau,\nu)\|_{N}^{2}\mathbf 1_{\Omega_N}$ for every $(\tau,\nu)\in\mathcal S_{\bf N}$, and since $8p(\mathbf m,\widehat{\bf m})\leqslant {\rm pen}(\mathbf m) + {\rm pen}(\widehat{\bf m})$,
\begin{eqnarray*}
 \|(\widehat a_{\widehat m_1},\widehat b_{\widehat m_2}) - (a,b)\mathbf 1_A\|_{N}^{2}
 & \leqslant &
 3\|(\widehat a_{m_1},\widehat b_{m_2}) - (a,b)\mathbf 1_A\|_{N}^{2}\\
 & &
 \hspace{2cm} +
 4{\rm pen}(\mathbf m) + 16\mathcal Z(\mathbf m,\widehat{\bf m})
 \quad {\rm on}\quad\Omega_N.
\end{eqnarray*}
So,
\begin{displaymath}
\mathbb E(\|(\widehat a_{\widehat m_1},\widehat b_{\widehat m_2}) -
(a,b)\mathbf 1_A\|_{N}^{2}\mathbf 1_{\Omega_N})
\leqslant
\min_{\mathbf m\in\mathcal M_N}\{
3\mathbb E(\|(\widehat a_{m_1},\widehat b_{m_2}) -
(a,b)\mathbf 1_A\|_{N}^{2}\mathbf 1_{\Omega_N}) +
4{\rm pen}(\mathbf m) + 16\rho(\mathbf m)\}.
\end{displaymath}
Therefore, to conclude this first step, by Inequality (\ref{adaptive_risk_bound_1}),
\begin{eqnarray}
 \label{adaptive_risk_bound_3}
 \mathbb E(\|(\widehat a_{\widehat m_1},\widehat b_{\widehat m_2}) - (a,b)\mathbf 1_A\|_{N}^{2})
 & \leqslant &
 \min_{\mathbf m\in\mathcal M_N}\{
 3\mathbb E(\|(\widehat a_{m_1},\widehat b_{m_2}) - (a,b)\mathbf 1_A\|_{N}^{2}\mathbf 1_{\Omega_N})\\
 & &
 \hspace{3.5cm} +
 4{\rm pen}(\mathbf m) + 16\rho(\mathbf m)\} + 2\mathfrak c_1N^{-\frac{r - 5}{2}}.
 \nonumber
\end{eqnarray}
{\bf Step 2.} First, consider $(\tau,\nu)\in\mathcal S_{\bf N}$, and let $M(\tau,\nu)$ be the $\mathbb F$-martingale defined by
\begin{displaymath}
M_s(\tau,\nu) :=
\sum_{i = 1}^{N}\int_{t_0}^{t}(\tau(X_{s}^{i}) +\nu(Y_{s}^{i}))\sigma(X_{s}^{i})dW_{1}^{i}(s)
\textrm{ $;$ }\forall t\in [t_0,T].
\end{displaymath}
Since $W^1,\dots,W^N$ are independent Brownian motions,
\begin{displaymath}
\langle M(\tau,\nu)\rangle_T =
\sum_{i = 1}^{N}\int_{t_0}^{T}(\tau(X_{s}^{i}) +\nu(Y_{s}^{i}))^2\sigma(X_{s}^{i})^2ds
\leqslant NT_0\|\sigma\|_{\infty}^{2}\|(\tau,\nu)\|_{N}^{2}.
\end{displaymath}
Then, by the Bernstein inequality for continuous local martingales (see Revuz and Yor \cite{RY99}, p. 153), for any $\varepsilon,\upsilon > 0$,
\begin{eqnarray*}
 \mathbb P(\zeta_N(\tau,\nu)\geqslant\varepsilon,\|(\tau,\nu)\|_{N}^{2}\leqslant\upsilon^2)
 & \leqslant &
 \mathbb P(M_T(\tau,\nu)^*\geqslant NT_0\varepsilon,
 \langle M(\tau,\nu)\rangle_T\leqslant NT_0\upsilon^2\|\sigma\|_{\infty}^{2})\\
 & \leqslant &
 \exp\left(-\frac{NT_0\varepsilon^2}{2\upsilon^2\|\sigma\|_{\infty}^{2}}\right).
\end{eqnarray*}
So, since this bound remains true by replacing $(\tau,\nu)$ by $-(\tau,\nu)$,
\begin{equation}\label{adaptive_risk_bound_4}
\mathbb P(|\zeta_N(\tau,\nu)|\geqslant\varepsilon,\|(\tau,\nu)\|_{N}^{2}\leqslant\upsilon^2)\leqslant
2\exp\left(-\frac{NT_0\varepsilon^2}{2\upsilon^2\|\sigma\|_{\infty}^{2}}\right).
\end{equation}
Now, for any $\mathbf m = (m_1,m_2)$ and $\mathbf m' = (m_1',m_2')$ belonging to $\mathcal M_N$, note that
\begin{displaymath}
\sup_{(\tau,\nu)\in\mathcal B_{\mathbf m,\mathbf m'}}|\zeta_N(\tau,\nu)|\leqslant
\sup_{h\in\mathbb B_{\mathbf m,\mathbf m'}}|\mathbb Z_N(h)|,
\end{displaymath}
where
\begin{displaymath}
\mathbb Z_N(h) :=
\frac{1}{NT_0}\sum_{i = 1}^{N}\int_{t_0}^{T}h(X_{s}^{i},Y_{s}^{i})\sigma(X_{s}^{i})dW_{1}^{i}(s),
\end{displaymath}
and
\begin{displaymath}
\mathbb B_{\mathbf m,\mathbf m'} :=
\left\{h\in\mathbb K_{\mathbf m,\mathbf m'} :
\int_{\mathbb R^2}h(x,y)^2f(x,y)dxdy\leqslant 1\right\}
\end{displaymath}
with
\begin{eqnarray*}
 \mathbb K_{\mathbf m,\mathbf m'} & = &
 \{h :\mathbb R^2\rightarrow\mathbb R\textrm{, measurable} :\\
 & &
 \hspace{2.5cm}
 \exists (\tau,\nu)\in S_{m_1\vee m_1'}\times\Sigma_{m_2\vee m_2'}\textrm{, }
 \forall (x,y)\in\mathbb R^2\textrm{, }h(x,y) =\tau(x) +\nu(y)\}.
\end{eqnarray*}
Since $\mathbb K_{\mathbf m,\mathbf m'}$ is a vector subspace of $\mathbb L^2(\mathbb R^2,f(x,y)dxdy)$ of dimension $D(\mathbf m,\mathbf m') := m_1\vee m_1' + m_2\vee m_2'$, by Lorentz et al. \cite{LGM96}, Chapter 15, Proposition 1.3, for every $n\in\mathbb N$, there exists $T_n\subset\mathbb B_{\mathbf m,\mathbf m'}$ such that $|T_n|\leqslant (3/\delta_n)^{D(\mathbf m,\mathbf m')}$ with $\delta_n =\alpha 2^{-n}$ ($\alpha\in (0,1)$), and for every $h\in\mathbb B_{\mathbf m,\mathbf m'}$,
\begin{equation}\label{adaptive_risk_bound_5}
\exists h_n\in T_n :\int_{\mathbb R^2}(h(x,y) - h_n(x,y))^2f(x,y)dxdy\leqslant\delta_{n}^{2}.
\end{equation}
Thanks to (\ref{adaptive_risk_bound_4}) and (\ref{adaptive_risk_bound_5}), and since $p(\mathbf m,\mathbf m')\asymp D(\mathbf m,\mathbf m')/N$, the third step of the proof of Marie \cite{MARIE25}, Theorem 3.2 - based on the chaining technique - extends from the case $b(.) = 0$ to the case $b(.)\neq 0$. Therefore, there exists a constant $\kappa_0 > 0$, not depending on $N$, such that for every $\kappa\geqslant\kappa_0$,
\begin{equation}\label{adaptive_risk_bound_6}
\rho(\mathbf m)\leqslant
\mathbb E\left(\left(\left[\sup_{h\in\mathbb B_{\mathbf m,\widehat{\bf m}}}|\mathbb Z_N(h)|\right]^2 -
p(\mathbf m,\widehat{\bf m})\right)_+\mathbf 1_{\Omega_N}\right)
\lesssim\frac{1}{N}\textrm{ $;$ }\forall\mathbf m\in\mathcal M_N.
\end{equation}
{\bf Step 3.} By pugging Inequality (\ref{adaptive_risk_bound_6}) in Inequality (\ref{adaptive_risk_bound_3}), there exists a constant $\mathfrak c_2 > 0$, not depending on $N$, such that
\begin{displaymath}
\mathbb E(\|(\widehat a_{\widehat m_1},\widehat b_{\widehat m_2}) - (a,b)\mathbf 1_A\|_{N}^{2})
\leqslant
\mathfrak c_2\left(
\min_{\mathbf m\in\mathcal M_N}\{
3\mathbb E(\|(\widehat a_{m_1},\widehat b_{m_2}) - (a,b)\mathbf 1_A\|_{N}^{2}\mathbf 1_{\Omega_N}) +
{\rm pen}(\mathbf m)\} +\frac{1}{N}\right).
\end{displaymath}
%


%
\subsubsection{Proof of Lemma \ref{relationship_collections_models}}\label{section_proof_relationship_collections_models} We acknowledge Huang \cite{HUANG25} for a decisive improvement of the result obtained in this lemma. Consider $\omega\in\Omega_N$ and
\begin{displaymath}
\widehat G_{\bf m}(\omega) :=
\Psi_{\bf m}^{-\frac{1}{2}}\widehat\Psi_{\bf m}(\omega)\Psi_{\bf m}^{-\frac{1}{2}}
\textrm{ $;$ }\forall\mathbf m\in\mathcal M_{N}^{+}.
\end{displaymath}
For any $\mathbf m\in\mathcal M_{N}^{+}$, since $\omega\in\Omega_{\bf m}$,
\begin{displaymath}
{\rm Sp}(\widehat G_{\bf m}(\omega))
\subset\left[\frac{1}{2},\frac{3}{2}\right],
\quad\textrm{and then}\quad
{\rm Sp}(\widehat G_{\bf m}^{-1}(\omega))
\subset\left[\frac{2}{3},2\right].
\end{displaymath}
Moreover,
\begin{displaymath}
\widehat\Psi_{\bf m}^{-1}(\omega) =
\Psi_{\bf m}^{-\frac{1}{2}}\widehat G_{\bf m}^{-1}(\omega)\Psi_{\bf m}^{-\frac{1}{2}}.
\end{displaymath}
So, thanks to a well-known property of the Loewner order,
\begin{displaymath}
\frac{2}{3}\mathbf x^*\Psi_{\bf m}^{-1}\mathbf x
\leqslant\mathbf x^*\widehat\Psi_{\bf m}^{-1}(\omega)\mathbf x\leqslant
2\mathbf x^*\Psi_{\bf m}^{-1}\mathbf x
\textrm{ $;$ }\forall\mathbf x\in\mathbb R^{m_1 + m_2},
\end{displaymath}
leading to
\begin{displaymath}
\|\Psi_{\bf m}^{-1}\|_{\rm op}\leqslant
\frac{3}{2}\|\widehat\Psi_{\bf m}^{-1}(\omega)\|_{\rm op}
\quad {\rm and}\quad
\|\widehat\Psi_{\bf m}^{-1}(\omega)\|_{\rm op}\leqslant
2\|\Psi_{\bf m}^{-1}\|_{\rm op}.
\end{displaymath}
Therefore, $\omega\in\Xi_N$.
%


%

%
\end{document}